\documentclass[12pt, reqno]{amsart}
\usepackage{graphicx}
\usepackage[letterpaper, margin=1in]{geometry}
\usepackage{times}
\usepackage{setspace}
\usepackage{soul}
 \usepackage{color}
\usepackage{tikz}
\usepackage{tikz-cd}
\usepackage{bm}
\usepackage{amsmath, amsfonts, amssymb, amsthm}
\usepackage{listings}
\usepackage{color}
\usepackage{mathtools}
\definecolor{dkgreen}{rgb}{0,0.6,0}
\definecolor{gray}{rgb}{0.5,0.5,0.5}
\definecolor{mauve}{rgb}{0.58,0,0.82}
\usepackage{braket}
\usepackage{pgfplots}\pgfplotsset{compat=1.7}
\usepackage{bbold}
\usepackage{ragged2e}
\usepackage{physics}
\usepackage{tcolorbox}
\tcbuselibrary{theorems}
\newtheorem{theorem}{Theorem}[section]  
\newtheorem{lemma}[theorem]{Lemma}      
\newtheorem{proposition}[theorem]{Proposition}
\newtheorem{corollary}[theorem]{Corollary}                
\theoremstyle{definition}
\newtheorem{definition}[theorem]{Definition}
\newtheorem{example}[theorem]{Example}

\newcommand{\al}{\aleph}
\newcommand{\mbb}{\mathbb}
\newcommand{\mcal}{\mathcal}
\usepackage{pgfplots}
\pgfplotsset{width=10cm,compat=1.9}

\usepackage{hyperref}
\hypersetup{
    colorlinks=false,
    linkbordercolor={1 0 0}, 
    citebordercolor={0 1 0},
    urlbordercolor={1 0 0},
    pdfborder={0 0 1}
}
\title{A Structural Analysis of Infinity in Set Theory and Modern Algebra}
\author{Noah Betz}
\date{}
\begin{document}
\begin{abstract}
\justifying
We present a self-contained analysis of infinity from two mathematical perspectives: set theory and algebra. We begin with cardinal and ordinal numbers, examining deep questions such as the continuum hypothesis, along with foundational results such as the Schröder–Bernstein theorem, multiple proofs of the well-ordering of cardinals, and various properties of infinite cardinals and ordinals. Transitioning to algebra, we analyze the interplay between finite and infinite algebraic structures, including groups, rings, and $R$-modules. Major results, such as the fundamental theorem of finitely generated abelian groups, Krull's Theorem, Hilbert's basis theorem, and the equivalence of free and projective modules over principal ideal domains, highlight the connections and differences between finite and infinite structures, as well as demonstrating the relationship between set-theoretic and algebraic treatments of infinity. Through this approach, we provide insights into how key results about infinity interact with and inform one another across set-theoretic and algebraic mathematics.
\end{abstract}
\maketitle
\tableofcontents
\begin{center}
\section*{Introduction}
\end{center}

The concept of infinity has fascinated mathematicians and philosophers alike for centuries. Unlike finite quantities, infinity represents something abstract, without bound, and often counterintuitive. Yet, infinity plays a fundamental role in modern mathematics. 

In set theory, infinity is viewed as a number corresponding to the size of a set, called a cardinal number. The properties of infinite cardinals in contrast to finite cardinals are fundamentally distinct, particularly in their arithmetic. Ordinal numbers provide a framework for distinguishing between two sets of the same cardinality. The well-ordering theorem, equivalent to the axiom of choice \cite{grabczewski, hewitt, howard, moore}, allows for a well-ordering of any non-empty set, and a unique ordinal number associated with that ordering.

Infinity with algebraic structures could either extend or cause natural generalizations, or cause a structural breakdowns. With groups, the structural differences between the direct product and the weak direct product illustrates challenges with taking the product of infinitely many elements. Classification and structural theorems, such as the fundamental theorem of finitely generated abelian groups and the classification of cyclic groups, further display the structural differences between finite and infinite groups. Ring structures, such as finite and infinite integral domains, polynomial rings, formal power series, and Noetherian Rings showcase natural extensions between finite and infinite structures as well as differences. In module theory, the differences between the theory of finite and infinite dimensional vector spaces, along with the structure of projective modules, further highlights the nuanced role infinity plays in algebraic mathematics.

The study of infinity leads to surprising results that challenge one's intuition. By presenting infinity through the lenses set theory and algebra, we aim to provide a self-contained, rigorous, and accessible treatment of how infinity interacts between the two disciplines.
\section{Set-Theoretic Perspective}\label{section1}
We assume the axioms of Zermelo–Fraenkel set theory (ZFC) (see \cite{shoenfieldaxioms}) with the notion of class left informal, denote the empty set as $\emptyset$, and assume, as in ZF \cite[Section 3]{shoenfieldaxioms}
 \begin{align*}
\exists X(\exists \emptyset \in X \hspace{1mm} \forall z(z \notin \emptyset) \land \forall y \in X \hspace{1mm} \exists z \in X \hspace{1mm} \forall w(w \in z \leftrightarrow w\ \in y \lor w = y)).
\end{align*} 
The statement above is called the $\textit{axiom of infinity}$. More compactly, with set-builder notation, it can be stated as \cite[Section 7.4]{stoll} \begin{align*}
\exists X(\emptyset \in X \land \forall x(x \in X \implies (x \cup \{x\}) \in X)).
\end{align*} 
A set $X$ satisfying the above statement is called an \textit{inductive set} \cite[Section 2.7]{roitman}. We define the intersection of all inductive sets as the set $\mbb{N}$ of \textit{natural numbers}, denoting its elements by
\begin{align*}
0  & \coloneqq \emptyset \\
1 & \coloneqq \{\emptyset\} \\
2 & \coloneqq \{\emptyset, \{\emptyset\}\} \\
3 & \coloneqq \{\emptyset, \{\emptyset\},\{\emptyset, \{\emptyset\}\}\} \\
\vdots
\end{align*}
We assume knowledge of the sets $\mbb{Z}$, $\mbb{Q}$, and $\mbb{R}$ of integer, rational, and real numbers, respectively. For constructions of $\mbb{Z}$ and $\mbb{Q}$, see \cite[Sections 3.3-4]{stoll}. For the dedekind cuts construction of $\mbb{R}$, see \cite[Chapter 1 Appendix]{babyrudin} and \cite[Section B]{axlermeasuresupplement}. For the equivalence classes of Cauchy sequences construction, see \cite[Section 2.1]{ziemer}. Additionally, under the respective isomorphisms if necessary, we write
\begin{align*}
\mbb{N} \subsetneq \mbb{Z} \subsetneq \mbb{Q} \subsetneq \mbb{R}.
\end{align*} 
We define $\mbb{Z}^{+} \coloneqq \{z \in \mbb{Z}|z>0\}$, calling it the set of \textit{positive integers}. Under the appropriate isomorphisms if necessary, we have $\mbb{Z}^{+} = \mbb{N} \setminus \{0\}$, and we do not distinguish between the two sets. The \textit{axiom of choice}, which states \cite[Section 6]{shoenfieldaxioms}
\begin{align*}
\forall X\left[\emptyset \notin X \implies \exists f:X \to \bigcup_{A \in X} A \hspace{2mm} \forall A \in X (f(A) \in A)\right]
\end{align*}
will be denoted in parentheses as (AC). The weaker countable choice and dependent choice will be stated as ($\text{AC}_{\omega}$) and (DC), respectively. Note that finite choice is provable in ZF \cite[Theorem 8.2]{suppes}, but countable choice is not \cite[Section 9.4]{potter}.
\subsection{Cardinal Numbers}\label{sub1.1}
In set theory, infinity is often viewed as a number corresponding to the size of a set, called a $\textit{cardinal number}.$ For a finite set $A$, the cardinal number corresponding to $A$ is simply the number of elements in the set $A$. However, for an infinite set, a more general definition needs to be made. For completeness, we recall some fundamental definitions. Given any sets $A$ and $B$, for $a \in A$, $b \in B$, define $(a,b) \coloneqq \{\{a\},\{a,b\}\}$ and \begin{align}\label{l1}
A \times B \coloneqq \left\{(a,b) \;\middle|\; a \in A, b \in B\right\}
\end{align}
The set (\ref{l1}) is called the $\textit{Cartesian product}$ of $A$ and $B$. A subset $R$ of $A \times B$ is called a  $\textit{relation}$ on $A \times B$. A $\textit{function}$ $f:A \to B$, is a relation $R$ on $A \times B$ with the following property: for every $a \in A$, there exists one and only one $b \in B$ such that $(a, b) \in R$. We write $f(a) \coloneqq b$.  Given functions $f:A \to B$ and $g:B \to C$, the function $x \mapsto g(f(x))$ ($x \in A$) is called the \textit{composition} of $g$ with $f$ and denoted $gf$. It is easy to verify function composition is associative ($(fg)h = f(gh)$).
\begin{definition}\label{def1.1}
Let $A$ be a set. An equivalence relation $R$ on the set $A$ is a relation on $A \times A$ with the following properties \\
(i) Reflectivity: for every $a \in A$,  $(a,a) \in R$. \\
(ii) Transitivity: If $(a,b) \in R$ and $(b,c) \in R$, then $(a,c) \in R$. \\
(iii) Symmetry: If $(a,b) \in R$, then so is $(b,a) \in R$. \\
If $R$ is an equivalence relation, then it is common to write $(a,b) \in R$ by $a \sim b$. An equivalence class of an element $a \in A$ is the set $$\bar{a} = \left\{b \in A \;\middle|\; a \sim b\right\}.$$
\end{definition}
A function $f:A \to B$ is $\textit{bijective}$ if it has the following properties: for every $b \in B$, there exists $a \in A$ such that $f(a) = b$, and if $ a\neq b$ for $a,b \in A$, then $f(a) \neq f(b)$. A function satisfying the first property is called $\textit{surjective}$, and a function satisfying the second property is called $\textit{injective}$. Given two sets $A$ and $B$, if there exists a bijection $f:A \to B$, then $A$ and $B$ are called $\textit{equipollent}$ sets, denoted $A \sim B$. The map $1_{A}:A \to A$, defined by $a \mapsto a$, is called the identity map on $A$. For a function $f:A \to B$, $S \subset B$, $C \subset A$, and $b \in B$, define the following sets \begin{align}
f(C) \coloneqq \left\{f(a) \;\middle|\; a \in C\right\}\label{l2} \\
f^{-1}(b) \coloneqq \left\{a \in A \;\middle|\; f(a) = b\right\}\label{l3} \\
f^{-1}(S) \coloneqq \left\{a \in A \;\middle|\; f(a) \in S\right\}\label{l4}.
\end{align}
The set (\ref{l2}) is the $\textit{image}$ of $C$ under $f$, (\ref{l3}) is the $\textit{pre-image}$ of $b \in B$ under $f$, and (\ref{l4}) the pre-image of $S$ under $f$. For sets $A,B$ and function $f:A \to B$, a function $g:B \to A$ such that $gf = 1_{A}$ is called a $\textit{left inverse}$ of $f$ and a function $h:B \to A$ such that $fh = 1_{B}$ a $\textit{right inverse}$ of $f$. A function that's both a left and right inverse of $f$ is called a $\textit{two-sided inverse}$ of $f$. 

Suppose $f:A \to B$ has a left-inverse $g$ and right-inverse $h$. Then \begin{align}\label{5} g = g1_{B} = g(fh) = (gf)h = 1_{A}h = h. \end{align}
Thus if $f$ has a left inverse and a right inverse, the inverses are the same map. Moreover, this shows that any two-sided inverse of $f$ is unique. The following theorem is a weaker version of \cite[Theorem 3J]{enderton}; however, it differs fundamentally, as it does not require the axiom of choice.
\begin{theorem}\label{thm1.1}
Let $f:A \to B$ be a function with $A$ non-empty. Then $f$ is bijective if and only if $f$ has a unique two-sided inverse, denoted $f^{-1}$.
\end{theorem}
\noindent $\textbf{Proof}.$ If $f:A \to B$ is bijective, consider the map $f^{-1}:B \to A$ by $b \mapsto a$, where $a \in f^{-1}(b)$. By injectivity, such element is unique, and by surjectivity, such element always exists, thus $f^{-1}$ is well-defined. Moreover, for any $a \in A$, $f^{-1}f(a) = f^{-1}(f(a)) = a$, and for any $b \in B$, $ff^{-1}(b) = f\qty(f^{-1}(b)) = b$, thus $f^{-1}$ is a two-sided inverse, and is unique by (\ref{5}). Conversely, if $f$ has a unique two-sided inverse $f^{-1}$, then for any $b \in B$, $f\qty(f^{-1}(b)) = b$ and $f$ is surjective. If $a,c \in A$ are such that $f(a) = f(c)$, then $a = f^{-1}\qty(f(a)) = f^{-1}\qty(f(c)) = c$ and $f$ is injective. Thus, $f$ is bijective, as desired. $\blacksquare$

\noindent Theorem \ref{thm1.1} can be extended to the following statement: \\
(i) $f$ is injective if and only if there is a map $g:B \to A$ such that $gf = 1_{A}$ \\
(ii) $f$ is surjective if and only if there is a map $h:B \to A$ such that $fh = 1_{B}$. \\
The proof of (i) does not require AC. The backward direction follows immediately from the fact that $a = gf(a) = gf(b) = b$ if $f(a) = f(b)$ ($a,b \in A$), and the forward direction by fixing some $a_0 \in A$ and defining
\begin{align*}
g(b) = \begin{cases}
a, \text{ if } b \in f(A) \land f(a) = b \\
a_0, \text{ if } b \notin f(A).
\end{cases}
\end{align*}
The proof of (ii), however, does require AC. In the forward direction, we observe $f^{-1}(b) \subset A$ is nonempty for each $b \in B$, then for each $b \in B$, choose $a_b \in f^{-1}(b)$ (AC). Then the map $b \mapsto a_b$ ($b \in B$) becomes a right inverse of $f$.

This illustrates how seemingly elementary properties of functions can hinge upon powerful assumptions like AC. While the existence of unique inverses for arbitrary bijective functions does not require choice due to the uniqueness of the pre-image of each element, the extension to arbitrary surjective functions depends crucially on AC. This reliance only emerges in the context of infinite sets. In the finite case, only finite choices are required to construct a right inverse. Therefore, it is precisely with infinite domains that the subtlety of choice of pre-images becomes logically and philosophically important, and highlights an important notion of how infinity intertwines deeply within the foundations of mathematics. Not in grand results like the continuum hypothesis, but through basic properties of functions. We now prove a few theorems to motivate the definition of cardinality on arbitrary sets.
\begin{theorem}\label{thm1.2}
Equipollence is an equivalence relation on the class of all sets
\end{theorem}
\noindent $\textbf{Proof}$. First, notice that $\emptyset \subset \emptyset \times \emptyset$ is (vacuously) a bijective function and $\emptyset \sim \emptyset$. For $A \neq \emptyset$, the identity function $1_{A}:A \to A$ defined by $1_A(a) = a$ is a bijection, and thus $A \sim A$ for every set $A$. Now suppose $A \sim B$. Then let $f:A \to B$ be a bijection. By Theorem \ref{thm1.1}, $f$ has two-sided inverse $f^{-1}:B \to A$. Moreover, $f^{-1}$ is bijective, as it has a two-sided inverse $f$, and thus $B \sim A$. Now suppose $A \sim B$ and $B \sim C$. Let $f:A \to B$ be a bijection and $g:B \to C$ a bijection. Consider the composition $gf:A \to C$. If $a,b \in A$ are such that $gf(a) = gf(b)$, then \begin{align}\label{l6} g(f(a)) = g(f(b)) \implies f(a) = f(b) \implies a = b \end{align} and $gf$ is injective. For any $c \in C$, there exists $b \in B$ such that $g(b) = c$, and for any $b \in B$, there exists $a \in A$ such that $f(a) = b$, whence $gf(a) = c$ and $gf$ is surjective. Thus, $gf$ is bijective and $A \sim C$, as desired. $\blacksquare$

The preceding theorem is in \cite[Theorem 0.8.1]{hungerford} (as an exercise). Let $\textit{I}_{0} = \emptyset$ and for each $n \in \mathbb{Z}^{+}$, define $\textit{I}_{n} = \left\{1,2,...,n\right\}$. The following proposition is in \cite[Exercise 0.8.1]{hungerford}.
\begin{proposition}\label{prop1.3}
Let $n,m \in \mathbb{N}$. \\
(i) $\textit{I}_{n}$ is not equipollent to any of its proper subsets \\
(ii) $\textit{I}_{n}$ is equipollent to $\textit{I}_{m}$ if and only if $n = m$.
\end{proposition}
\noindent $\textbf{Proof}$. [(i)] Let $S$ be the set of all $n \in \mathbb{N}$ such that $\textit{I}_{n}$ is equipollent to one of its proper subsets. Suppose, by way of contradiction, that $S \neq \emptyset$. By the well-ordering principle, $S$ has a least element, $k \in S$. Since $\textit{I}_{k}$ has a proper subset, $k \geq 1$. If $k = 1$, then the only proper subset of $\textit{I}_{k}$ is $\emptyset$, and the only function $\emptyset \to \textit{I}_{k}$ is the empty function, which clearly isn't bijective. Thus $k > 1$ and  Let $f:A \to \textit{I}_{k}$ be a bijection with $k \notin A$ and $A \subsetneq \textit{I}_{k}$ (if $k \in A$, then choose some $j \in \textit{I}_k \setminus A$ and replace $k$ with $j$, having it map to the same element as $k$). Let $B = A \setminus \left\{f^{-1}(k)\right\}$ and $C = \textit{I}_{k} \setminus \left\{k\right\} = \textit{I}_{k-1}$. The restriction $f|_{B}:B \to C$ is a bijection, contradicting the minimality of $k$.

[(ii)] If $\textit{I}_{n} \sim \textit{I}_{m}$, then by (i), $m$ cannot be less than $n$ (or else $\textit{I}_{m}$ is a proper subset of $\textit{I}_{n}$) and likewise $n$ cannot be less than $m$, and thus $n = m$. If $n = m$, then the identity map is a bijection and $\textit{I}_{n} \sim \textit{I}_{m}$, as desired. $\blacksquare$ 

Proposition \ref{prop1.3} highlights a fundamental characterization of finite sets by showing there exists no bijection between a finite set and any of its proper subsets, and moreover that (as a natural consequence of the first statement) two finite sets have a bijection between them if and only if they have the same number of elements. In the context of infinite sets, however, equipollence with proper subsets becomes a defining feature (e.g., $\mbb{N}$ and $2\mbb{N}$ are equipollent with $\mbb{N} \to 2\mbb{N}$ by $n \mapsto 2n$). Therefore, this result serves as a fundamental contrast between finite and infinite sets, to then highlight crucial differences between finite and infinite cardinal numbers, and illustrating crucial structural differences between finite and infinite set-theoretic universes.

The cardinal number of a set $A$ attempts to answer how many elements are in the set $A$. If $A$ has $n$ elements, say $A = \left\{a_1,...,a_n\right\}$ for distinct elements $a_1,...,a_n$, then the map given by $a_j \mapsto j$ ($1 \leq j \leq n$) is clearly a bijection and $A \sim \textit{I}_{n}$. Moreover, Proposition \ref{prop1.3} and Theorem \ref{thm1.2} show that $A$ is not equipollent to any other $\textit{I}_n$. Thus, to say $A$ has $n$ elements means precisely that $A \sim \textit{I}_{n}$, or in other words, $A$ and $\textit{I}_n$ are in the same equivalence class under the equivalence relation of equipollence. Therefore, for a finite set $A$ (that is, a set equipollent to some $\textit{I}_n$), the equivalence class of $A$ under equipollence provides the necessary information to determine the number of elements in $A$. Extending this notion to infinite sets gives the definition 
\begin{definition}\label{def3}
The cardinal number (or cardinality) of a set $A$ is the equivalence class of $A$ under the equivalence relation of equipollence, denoted $|A|$. If $A$ is a finite set, then $|A|$ is a finite cardinal, and if $A$ is an infinite set, then $|A|$ is an infinite cardinal.
\end{definition}
The definition of cardinality presented in this paper follows the presentation of \cite[Section 0.8]{hungerford}, where cardinality is a proper class rather than a set, as in \cite[Section 6.2]{enderton}. This definition is a prime example of how the concept of infinity tends to complicate concepts. Without infinity, it becomes clear how to define cardinality, simply count the number of elements in the set (which need be finite), and that is the set's cardinality. But infinity complicates this matter, for how can one count an infinite number of elements and make that the set's cardinality? 

Naturally, we identify the natural number $n$ with the cardinal number $|\textit{I}_n|$ and write $|\textit{I}_n| = n$, so that the cardinal number of a finite set is precisely the number of elements in the set. For $|A| = \alpha$ and $|B| = \beta$ cardinal numbers, and $A,B$ disjoint, the sum $\alpha + \beta$ is $|A \cup B|$ and product $\alpha \beta$ is $|A \times B|$. It is easy to see this definition is independent of the sets $A,B$ chosen. Moreover, for $\alpha \beta$, the sets need not be disjoint.
\begin{definition}\label{def4}
Let $A$ be a set. A $\textit{partial ordering}$ on $A$ is a relation $\leq$ on $A \times A$ that is \\
(i) Reflective: for every $a \in A$, $a \leq a$ \\
(ii) Transitive: if $a \leq b$ and $b \leq c$, then $a \leq c$ ($a,b,c \in A$) \\
(iii) Antisymmetric: if $a \leq b$ and $b \leq a$, then $a = b$ ($a,b \in A$). \\
Two elements $a,b \in A$ are comparable if $a \leq b$ or $b \leq a$. A partial ordering $R$ is a $\textit{linear ordering}$ if every $a,b \in A$ are comparable with respect to $R$. $R$ is a $\textit{well-ordering}$ if every subset of $A$ has a least element with respect to $R$.
\end{definition}
\begin{example}\label{ex4}
Given a set $X$, ordering by inclusion is a partial ordering that need not be a linear ordering. For example, ordering by inclusion on subsets of $\{1,2\} \subset \mbb{N}$ is not linearly ordered since $\{1\}$ and $\{2\}$ are not comparable. In contract, ordering by inclusion on only the subsets $\mathcal{A} = \{\emptyset, \{1\}, \{1,2\}\}$ is linearly ordered, since $\emptyset \subset \{1\} \subset \{1,2\}$.
\end{example}
\begin{example}\label{ex5}
Consider the set $A \times B$ for partially ordered sets $(A, \leq_{A})$ and $(B, \leq _{B})$. Then the \textit{dictionary ordering} on $A \times B$ is the partial ordering defined by $(a,b) \leq (a',b')$ if and only if $a <_{A} a'$ (i.e., $a \leq_{A} a'$ and $a \neq a'$) or $a = a'$ and $b \leq_{B} b'$. If $A$ and $B$ are linearly ordered, then so is $A \times B$ with respect to this ordering.
\end{example}
This definition aligns with \cite[Section 14]{halmos}, which includes reflexivity, in contrast to non-reflexivity like in \cite[Section 1.3]{munkres}. Every well-ordered set $(A, \leq)$ is also linearly ordered, for any subset $\{a,b\}$ of $A$ must have a least element, so $a \leq b$ or $b \leq a$. Given a partially ordered set $(A, \leq)$, a $\textit{maximal element}$ $a \in A$ is an element such that for any $c \in A$ comparable to $a$, $c \leq a$. An $\textit{upper bound}$ of a subset $B$ of $A$ is an element $x \in A$ such that $b \leq x$ for every $b \in B$. A nonempty subset $B$ of $A$ that is linearly ordered is called a $\textit{chain}$ in $A$. $\textit{Zorn's lemma}$, which is equivalent to the axiom of choice (see \cite[Section 1.3]{hewitt}, \cite{grabczewski}, \cite{moore}, or \cite{howard}), states that if $A$ is a non-empty partially ordered set such that every chain in $A$ has an upper bound in $A$, then $A$ contains a maximal element. For two sets $A,B$, $|A| \leq |B|$ if there exists an injection $f:A \to B$, and $|A| = |B|$ if there is a bijection $f:A \to B$ (since $A$ and $B$ would have the same equivalence class under equipollence). If $|A| \leq |B|$ and $|A| \neq |B|$, then it is denoted $|A| < |B|$. Sometimes, $|A| \leq |B|$ or $|A| < |B|$ will be denoted instead as $|B| \geq |A|$ or $|B| > |A|$, respectively.
\begin{example}\label{ex1}
The sets $\mbb{N}$, $\mbb{Z}$, $\mbb{Q}$, and $\mbb{R}$ all have infinite cardinalities. Under inclusion, clearly $|\mbb{N}| \leq |\mbb{Z}| \leq |\mbb{Q}| \leq |\mbb{R}|$. However unintuitive, it holds that $|\mbb{N}| = |\mbb{Z}| = |\mbb{Q}|$ (Proposition \ref{prop1.10}). For the real numbers however, it holds that $|\mbb{R}| = |\mathcal{P}(\mbb{N})| > \mbb{N}$ (Theorem \ref{thm1.15}, \ref{thm1.11}).
\end{example}
\begin{lemma}[Schröder–Bernstein]\label{lem1.4}
If $A,B$ are sets such that $|A| \leq |B|$ and $|B| \leq |A|$, then $|A| = |B|$. 
\end{lemma}
\noindent $\textbf{Proof}$. Let $f:A \to B$ be an injection and $g:B \to A$ an injection. Let $E_0 = A \setminus g(B)$ and define recursively $E_{n+1} = g\qty(f\qty(E_n))$. Let $E = \bigcup_{n \in \mbb{N}} E_n$ and define $h:A \to B$ by $$h(x) = \begin{cases} 
f(x), \text{ if } x \in E \\
g^{-1}(x), \text{ if } x \in A-E.
\end{cases}$$
First off, notice that $$x \in A-E \implies x \notin E_{0} \implies x \in g(B).$$
Thus $h$ is well-defined. If $x,y \in E$ and $h(x) = h(y)$, then $f(x) = f(y)$ and $x = y$, so $h$ is injective on $E$. If $x,y \in A-E$ and $h(x) = h(y)$, then $g^{-1}(x) = g^{-1}(y)$, and applying $g$ to both sides gives $x = y$, so $h$ is injective on $A-E$. Now if $x \in E$ and $y \in A-E$, and $h(x) = h(y)$, then $f(x) = g^{-1}(y)$, or $y = gf(x)$. Since $x \in E$, for some $n \in \mbb{N}$, $x \in E_n$. Thus $y = gf(x) \in gf(E_n) = E_{n+1} \subset E$, a contradiction. It follows that $h$ is injective. Now let $b \in B$ and consider $x = g(b)$. If $x \in E$, then $x \in E_m$ for some $m \in \mbb{N}$, and $m > 0$ since $x \in g(B)$.  Write $x = gf(x')$ with $x' \in E_{m-1}$, and notice $$b = g^{-1}(x) = g^{-1}\qty(gf(x')) = f(x') = h(x').$$ 
Finally, if $x \notin E$, then $$h(x) = g^{-1}(x) = g^{-1}\qty(g(b)) = b.$$ Thus $h$ is surjective, and $|A| = |B|$, as desired. $\blacksquare$ 

For further proofs of the preceding theorem, see \cite[Section 6.4]{enderton}, \cite[Theorem 2.7]{kaplansky}, and \cite[Theorem 2.3.1]{stoll}. Lemma \ref{lem1.4} shows that the ordering $\leq$ on cardinal numbers is anti-symmetric. While this result is natural and intuitively clear for finite sets and follows from the pigeonhole principle, for infinite sets such as $\mbb{N}$ or $\mbb{R}$, neither the pigeonhole principle nor counting-based arguments are applicable.

The importance of Schröder–Bernstein becomes particularly important when studying infinite sets. With finite sets, bijections are generally easily constructed or shown to not exist, however with infinite sets, bijections are generally non-trivial to explicitly construct and are much more subtle. With infinite sets, it is often considerably simpler to construct two injections than to construct one bijection, particularly in cases when, for example, the sets in questions are embedded in one another. This is illustrated by Example \ref{ex2} and the several key results throughout this paper that use a ``cardinality squeeze" in its proof (e.g., Proposition \ref{prop1.10}, Theorem \ref{thm1.15}, Theorem \ref{thm2.11}). It is Schröder–Bernstein that allows one to introduce that cardinality squeeze. 

Therefore, Lemma \ref{lem1.4} not only demonstrates a foundational result in the theory of cardinal numbers, but also reflects the richness and intricacy of infinite set theory, where seemingly intuitive and straightforward assumptions like it require much more deeper machinery. In that sense, it illustrates a fundamental distinction between the more well-behaved structure of finite sets and the more nuanced structure of the infinite.
\begin{example}\label{ex2}
Consider the real numbers $\mbb{R}$ and the interval $[0,1] \subset \mbb{R}$. By inclusion $\abs{[0,1]} \leq |\mbb{R}|$. Moreover, the map $f(x) = \frac{1}{2}+\frac{\arctan x}{\pi}$ is injective since $\arctan x$ is, and $-\frac{\pi}{2} \leq \arctan x \leq \frac{\pi}{2}$ so $\Im f \subset [0,1]$. Thus $|\mbb{R}| \leq |[0,1]|$ and by Lemma \ref{lem1.4} $|\mbb{R}| = \abs{[0,1]}$.
\end{example}
The concept of the Cartesian product of two sets in (\ref{l1}) can be generalized to an arbitrary family of sets as follows. Let $\left\{A_{i}\right\}_{i \in \textit{I}}$ be a family of sets. Then the $\textit{Cartesian product}$ of $\left\{A_{i}\right\}_{i \in \textit{I}}$ is the set of all functions $f:\textit{I} \to \bigcup_{i \in \textit{I}} A_i$ such that $f(i) \in A_{i}$ for every $i \in \textit{I}$, and is denoted $\prod_{i \in \textit{I}} A_{i}$. We denote an element of the Cartesian Product by $\left\{x_{i}\right\}_{i \in \textit{I}} \in \prod_{i \in \textit{I}} A_{i}$, where each $x_j$ is the output $f(j) \in A_{j}$. The axiom of choice guarantees that the Cartesian product of a family of non-empty sets over a non-empty index is non-empty, or in other words, that one can ``choose" an element from each set. For $\left\{A_{i}\right\}_{i \in \textit{I}}$ a family of sets, define the map $\pi_{j}:\prod_{i \in \textit{I}} A_{i} \to A_{j}$ by \begin{align}\label{l7}
    \pi_{j}\qty(f) = f(j).
\end{align}
The map (\ref{l7}) is called the $\textit{canonical projection}$ from $\prod_{i \in \textit{I}} A_{i}$ onto $A_{j}$.

The following result is foundational to the theory of cardinal numbers, for it ensures that given sets $A$ and $B$, their cardinalities can be compared, and also that given a family of cardinals, there is a smallest one. Since cardinal numbers generalize natural numbers (the natural numbers are finite cardinals), it is natural to seek properties that mirror that of $\mbb{N}$, particularly the well-ordering of $\mbb{N}$ (and hence follows the comparability of any two elements).

Cardinality intuitively reflects the ``size" of a set, and one may (naively) expect that from this any two cardinals may be compared and that a family of sizes contains a smallest element. Infinity, however, complicates this matter substantially. The naive approach of comparing sizes by some counting argument fails, and one cannot simply examine all elements from a family to determine the least one. These difficulties demonstrate the subtleties of infinite sets, highlights the complexity of infinity in set-theoretic results, and motives the use of nontrivial methods such as Zorn's Lemma. The following proof follows a structure similar to Hönig \cite{honig}, but with added details and different notation for clarity.
\begin{theorem}\label{thm1.5}
The class of all cardinal numbers is well-ordered by $\leq$.
\end{theorem}
\noindent $\textbf{Proof}$. Given any set $A$, the identity map $1_A$ is clearly injective and $|A| \leq |A|$. If $A,B,C$ are sets such that $|A| \leq |B|$ and $|B| \leq |C|$, then let $f:A \to B$ and $g:B \to C$ be injections. Then (\ref{l6}) shows the composition $gf$ is an injection and $|A| \leq |C|$. Lemma \ref{lem1.4} shows $\leq$ is antisymmetric and therefore $\leq$ is a partial ordering.  

Let $A = \left\{A_{i}\right\}_{i \in \textit{I}}$ be a family of sets with cardinalities $|A_{i}| = \alpha_{i}$ for each $i \in \textit{I}$. Let $S = \left\{\alpha_{i} \mid \alpha_{i} \text{ is finite}\right\}$. If $S \neq \emptyset$, then by the well-ordering principle, $S$ has a least element $\alpha$. Clearly $\alpha \leq \alpha_{i}$ for any infinite $\alpha_{i}$ as well, and we are done. So assume $S = \emptyset$. It suffices to show that there exists $i_{0} \in \textit{I}$ such that for every $i \in \textit{I}$, there exists an injection $f:A_{i_{0}} \to A_{i}$. Consider the class $\mcal{B}$ of subsets $B$ of $\prod_{i \in \textit{I}} A_{i}$ such that \begin{align}\label{l8}
    x = \{x_{i}\}_{i \in \textit{I}} \in B, y = \{y_{i}\}_{i \in \textit{I}} \in B \land x \neq y \implies x_{i} \neq y_{i}, \forall i \in \textit{I}.
\end{align}  
Take a subset $B = \left\{\{x_{i}\}_{i \in \textit{I}}\right\}$ such that each $x_{i}$ is chosen from each $A_{i}$ (AC). Clearly $B \in \mcal{B}$ (vacuously) and $\mcal{B} \neq \emptyset$. Partially order $\mcal{B}$ by inclusion and let $C = \{B_{j} \mid j \in J\}$ be a chain in $\mcal{B}$. Consider $\bigcup_{j \in J} B_{j}$ . If $x,y \in \bigcup_{j \in J} B_{j}$, then $x \in B_{i}$ and $y \in B_{j}$ for some $i, j \in \textit{J}$. Without loss of generality, assume $B_{i} \subset B_{j}$. Then $x,y \in B_{j} \in \mcal{B}$ and thus satisfy (\ref{l8}), whence $\bigcup_{j \in J} B_{j}$ satisfies (\ref{l8}) as well and $\bigcup_{j \in J} B_{j} \in \mcal{B}$. Clearly $B_{j} \subset  \bigcup_{j \in J} B_{j}$ for every $j \in J$ and thus $\bigcup_{j \in J} B_{j} \in \mcal{B}$  is an upper bound of $C$ and $\mcal{B}$ fits the conditions for Zorn's lemma. Let $M$ be a maximal element. We claim there exists $i_{0} \in \textit{I}$ such that $\pi_{i_{0}}\qty(M) = A_{i_{0}}$. So suppose, by way of contradiction, that $\pi_{i}\qty(M) \neq A_{i}$ for every $i \in \textit{I}$. Then choose $a_{i} \in A_{i}-\pi_{i}\qty(M)$ for each $i \in \textit{I}$ (AC). Then $M \cup \left\{a_{i}\right\}_{i \in \textit{I}}$ strictly contains $M$, but fits (\ref{l8}), contradicting the maximality of $M$. Thus, let $i_{0} \in \textit{I}$ be such that $\pi_{i_{0}}\qty(M) = A_{i_{0}}$. Define $f_{i}:A_{i_{0}} \to A_{i}$ by $x_{i_{0}} \mapsto \pi_{i}(x)$ for each $i \in \textit{I}$, where $x = \{x_{i}\}_{i \in \textit{I}} \in M$ is such that $\pi_{i_{0}}(x) = x_{i_{0}}$. By (\ref{l8}), $x$ is unique and $f_{i}$ well-defined. Moreover, if $x_{i_{0}} \neq y_{i_{0}}$, then $x \neq y$ (by definition) and $x_{i} \neq y_{i}$ (by (\ref{l8})). Thus, $f_{i}$ is injective, as desired. $\blacksquare$ 

The cardinality of the natural numbers $|\mathbb{N}|$ is denoted $\bm{\aleph}_{0}$.  A set $S$ equipollent to $\mathbb{N}$ is called a $\textit{denumerable}$ set. For further proofs and discussion of the following proposition, see \cite[Section 2.1]{kaplansky} and \cite[Theorem 0.8.8]{hungerford}.
\begin{proposition}\label{prop1.6}
Every infinite set has a denumerable subset
\end{proposition}
\noindent $\textbf{Proof}$. Let $S$ be an infinite set. Choose $x_0 \in S$ and recursively define $x_n \in S \setminus \{x_0,...,x_{n-1}\}$ ($n \in \mbb{Z}^{+}$), which is possible since $S$ is infinite. Consider the set $X = \{x_n \mid n \in \mbb{N}\} \subset S$ ($\text{AC}_{\omega}$). The map $x_n \mapsto n$ ($n \in \mbb{N}$) is clearly bijective hence $X$ is a denumerable subset of $S$, as desired. $\blacksquare$

Proposition \ref{prop1.6} requires, and is equivalent to, the axiom of countable choice \cite{howard}. Moreover, it shows that $\bm{\al}_{0}$ is in fact the smallest infinite cardinal. For given any infinite set $A$ with $|A| = \alpha$, and a denumerable subset $B$, the inclusion map $B \to A$ is injective and so $\bm{\al}_{0} = |B| \leq |A| = \alpha$. The following Proposition follows the same proof strategy as in \cite[Lemma 0.8.9]{hungerford}, but with additional details for completeness and clarity.
\begin{proposition}\label{prop1.7}
If $A$ is an infinite set and $F$ a finite set, then $|A \cup F| = |A|$.
\end{proposition}
\noindent $\textbf{Proof}.$ First, suppose $A \cap F = \emptyset$. Let $F = \{c_0,...,c_{n-1}\}$ ($c_i$ distinct) and $D = \{x_{i} \mid i \in \mbb{N}\}$ ($x_i$ distinct) a denumerable subset of $A$ by Proposition \ref{prop1.6}. Consider the map $f:A \to A \cup F$ given by, $$f(x) = \begin{cases}
    c_{i}, \text{ if } x = x_{i}, 0 \leq i \leq n-1 \\
    x_{i-n} \text{ if } x = x_{i}, i>n-1 \\
    x,  \text{ if } x \in A-D.
\end{cases}$$
Clearly $f$ is surjective. Let $B = \{x_{i} \hspace{1mm}| \hspace{1mm} 0 \leq i \leq n-1\}$. If $x,y$ are such that $x \neq y$ and $x,y \in B$, or $x,y \in D-B$, or $x,y \in A-D$, then clearly $f(x) \neq f(y)$. If $x$ and $y$ are in different sets, then $f(x) \neq f(y)$ because $A$ is disjoint to $F$ and $D \cap (A-D) = \emptyset$. Thus $f$ is bijective and $|A \cup F| = |A|$. If $A \cap F \neq \emptyset$, then $A \cup F = A \cup (F-A)$ and $|A \cup F| = |A \cup (F-A)| = |A|$ (since $F-A$ is disjoint to $A$), as desired. $\blacksquare$ 

Proposition \ref{prop1.7} illustrates a key divergence between finite and infinite cardinalities. With finite sets, adding even a single element to a set increases its cardinality, but in the case of infinite sets, such intuitions break down and the proposition shows that adjoining any finite number of elements maintains the same cardinality. Thus the proposition highlights a fundamental discontinuity between infinite and finite set-theoretic mathematics; the way arithmetic operates within them.

To see this more concretely, at least in the case where $|A| = \bm{\al}_{0}$, imagine a hotel with countably infinitely many rooms, all rooms completely filled. If $n$ new guests show up, then shifting each guest $n$ rooms down leaves precisely $n$ rooms open for the new guests to fill. This hotel analogy is often called Hilbert's hotel, introduced by Hilbert in his lecture notes \cite[Über das Unendliche]{hilbertlect}, and captures the intuition for the statement $\bm{\al}_{0}+n = \bm{\al}_{0}$. We need to prove two more strong and counterintuitive results about the arithmetic of infinite cardinals. For further proofs and discussion of the following Lemma, see \cite[Lemma 0.8.10]{hungerford} and \cite{whitehead}.
\begin{lemma}\label{lem1.8}
If $\alpha, \beta$ are cardinal numbers such that $\beta \leq \alpha$ and $\alpha$ is infinite, then $\alpha + \beta = \alpha$. 
\end{lemma}
\noindent $\textbf{Proof}.$ Let $A,B$ be sets with $|A| = \alpha$, $|B| = \beta$, and $A \cap B = \emptyset$. Clearly $\alpha \leq \alpha+\beta$ (inclusion map). Let $C$ be such that $|C| = \alpha$ and $A \cap C = \emptyset$. Consider an injection $f:B \to C$ and define the map $h:A \cup B \to A \cup C$ by $h|_{A} = 1_{A}$ and $h|_{B} = f$. Since $B$ and $C$ are disjoint to $A$, it follows that $h$ is a well-defined, injective map, and $\alpha+\beta = |A \cup B| \leq |A \cup C| = \alpha+\alpha$. Therefore, by Lemma \ref{lem1.4} it suffices to show $\alpha+\alpha = \alpha$. Let $\mcal{F}$ be the set of all pairs $(f,X)$, where $X \subset A$ and $f:X \times \{0,1\} \to X$ is a bijection. Define an order $\leq$ on $\mcal{F}$ by \begin{align}\label{l9}
    \qty(f_{1},X_{1}) \leq \qty(f_{2},X_{2}) \iff X_{1} \subset X_{2} \land f_{2}|_{X_{1} \times \{0,1\}} = f_{1}.
\end{align}
First off, since $A$ is infinite, by Proposition \ref{prop1.6}, it has a denumerable subset $D$. Notice the map $f:\mbb{N} \times \{0,1\} \mapsto \mbb{N}$ by $(n,0) \mapsto 2n$ and $(n,1) \mapsto 2n+1$ is a bijection. Let $g:D \to \mbb{N}$ be a bijection. Then clearly the map $h:D \times \{0,1\} \to \mbb{N} \times \{0,1\}$ by $(d,0) \mapsto (g(d),0)$ and $(d,1) \mapsto (g(d),1)$ for $d \in D$ is bijective, and moreover the composition \begin{align}\label{l10}
    D \times \{0,1\} \xrightarrow{h} \mbb{N} \times \{0,1\} \xrightarrow{f} \mbb{N} \xrightarrow{g^{-1}} D 
\end{align}
is a bijection and thus $\mcal{F} \neq \emptyset$. Clearly, $(f,X) \leq (f,X)$ for any $(f,X) \in \mcal{F}$. If $(f_{1},X_{1}) \leq (f_{2},X_{2})$ and $(f_{2},X_{2}) \leq (f_{3},X_{3})$, then $X_{1} \subset X_{2} \subset X_{3}$ and $f_{3}|_{X_{1} \times \{0,1\}} = f_{2}|_{X_{1} \times \{0,1\}} = f_{1}$ (since $X_{1} \times \{0,1\} \subset X_{2} \times \{0,1\}$) and thus $(f_{1},X_{1}) \leq (f_{3},X_{3})$. If $(f_{1},X_{1}) \leq (f_{2},X_{2})$ and $(f_{2},X_{2}) \leq (f_{1},X_{1})$, then $X_{1} \subset X_{2}$ and $X_{2} \subset X_{1}$, or $X_{1} = X_{2}$, and $f_{2} = f_{2}|_{X_{1}} = f_{1}$, and so $(f_{1},X_{1}) = (f_{2},X_{2})$. Therefore, $\leq$ is a partial ordering on $\mcal{F}$. 

 Let $B = \left\{(f_{i},X_{i}) \mid i \in \textit{I} \right\}$ be a chain in $\mcal{F}$. Let $X = \bigcup_{i \in \textit{I}} X_{i}$ and define $f:X \times \{0,1\} \to X$ by $f(x) = f_{i}(x)$ for $x \in X_{i} \times \{0,1\}$. Suppose $x \in  X_{i} \times \{0,1\}$ and $y \in X_{j} \times \{0,1\}$ are such that $x = y$. Without loss of generality, assume  $(f_{i},X_{i}) \leq (f_{j},X_{j})$. Then, $f(y) = f_{j}(y) = f_{j}(x) = f_{i}(x) = f(x)$ and $f$ is well-defined. Moreover, if $x$ and $y$ are such that $f(x) = f(y)$, then $f_{j}(y) = f_{i}(x) = f_{j}(x)$ and $x = y$. Since each $f_{i}$ is surjective, it immediately follows that $f$ is surjective and thus $(f,X) \in \mcal{F}$. Finally, given any $(f_{i},X_{i}) \in B$, clearly $X_{i} \subset X$ and $f|_{X_{i} \times \{0,1\}} = f_{i}$ so $(f_{i},X_{i}) \leq (f,X)$, and $(f,X)$ is an upper bound of $B$, whence $\mcal{F}$ fits the conditions for Zorn's lemma. 

Let $(g,Y) \in \mcal{F}$ be a maximal element. Clearly, $Y_{0} = \{(y,0) \mid y \in Y\}$ and $Y_{1} = \{(y,1) \mid y \in Y\}$ are disjoint sets such that $|Y| = |Y_{0}| = |Y_{1}|$ and $Y \times \{0,1\} = Y_{0} \cup Y_{1}$. Therefore, by definition \begin{align}\label{l11}
    |Y| = |Y \times \{0,1\}| = |Y_{0} \cup Y_{1}| = |Y_{0}|+|Y_{1}| = |Y|+|Y|.
\end{align}
We now show that $|Y| = \alpha$, completing the proof. If $A-Y$ is infinite, it would contain a denumerable subset $E$ by Proposition \ref{prop1.6}, and, as in (\ref{l10}), there exists a bijection $\varphi:E \times \{0,1\} \to E$. By combining $\varphi$ with $g$, we can construct a bijection $h:\qty(Y \cup E) \times \{0,1\} \to Y \cup E$ such that $(g,Y) \leq \qty(h, Y \cup E)$ but $(g,Y) \neq \qty(h, Y \cup E)$, contradicting the maximality of $(g,Y)$. Therefore, $A-Y$ must be finite. Since $A$ is infinite, and $A = Y \cup (A-Y)$, $Y$ must also be infinite. Thus, by Proposition \ref{prop1.7}, $|Y| = |Y \cup (A-Y)| = |A| = \alpha$, as desired. $\blacksquare$ 

As with Proposition \ref{prop1.7}, this highlights a fundamental distinction in the arithmetic of infinite cardinals from their finite counterparts. Continuing with the hotel analogy, if $\alpha = \bm{\al}_{0}$ and $\beta$ is finite, the above analogy covers it. If $\beta = \bm{\al}_{0}$, that's like if $\bm{\al}_{0}$ guests arrived at the full hotel. Then, having each of the guests move to double their room number gives room for the $\bm{\al}_{0}$ new guests (they each fill the odd numbered rooms). This then captures the intuitive nature of the property $\bm{\al}_{0}+\bm{\al}_{0} = \bm{\al}_{0}$. The partial ordering (\ref{l9}) is called an $\textit{ordering by extension}$ (see \cite[Section 0.8]{hungerford}). For further proofs and discussion of the following Lemma, see \cite[Theorem 0.8.11]{hungerford} and \cite{whitehead}.
\begin{theorem}\label{thm1.9}
If $\alpha$ and $\beta$ are cardinal numbers such that $0 \neq \beta \leq \alpha$ and $\alpha$ is infinite, then $\alpha \beta = \alpha$.
\end{theorem}
\noindent $\textbf{Proof}.$ Clearly $\alpha \leq \alpha \beta$. Let $A,B$ be sets such that $|A| = \alpha$ and $|B| = \beta$. Let $f:B \to A$ be an injection, and then the map $A \times B \to A \times A$ by $(a,b) \mapsto (a,f(b))$ is clearly injective and $\alpha \beta \leq \alpha \alpha$. Thus, by Lemma \ref{lem1.4} it suffices to show $\alpha = \alpha \alpha$. Let $\mcal{F}$ be the set of all bijections $f:X \times X \to X$, where $X$ is an infinite subset of $A$. Clearly $|\mbb{N}| \leq |\mbb{N} \times \mbb{N}|$ ($n \mapsto (n,0)$). The map $\mbb{N} \times \mbb{N} \to \mbb{N}$ by $(n,m) \mapsto 2^{n} 3^{m}$ is injective by the fundamental theorem of arithmetic \cite[Theorem 1.14]{niven}, and thus $|\mbb{N}| = |\mbb{N} \times \mbb{N}|$ so there exists bijection $f:\mbb{N} \to \mbb{N} \times \mbb{N}$ (Lemma \ref{lem1.4}). By Proposition \ref{prop1.6}, $A$ has denumerable subset $D$. Let $g:D \to N$ be a bijection. Clearly the map $h: D \times D \to \mbb{N} \times \mbb{N}$ by $h(a,b) = (g(a),g(b))$ is bijective and moreover the composition \begin{align}\label{l12}
D \times D \xrightarrow{h} \mbb{N} \times \mbb{N} \xrightarrow{f^{-1}} \mbb{N} \xrightarrow{g^{-1}} D
\end{align}
is bijective and $\mcal{F} \neq \emptyset$. Partially order $\mcal{F}$ by extension (as in (\ref{l9})) and by Zorn's lemma there exists maximal element $g:M \times M \to M$. Since $g$ is bijective, $|M||M| = |M \times M| = |M|$. 

We now show $|M| = |A| = \alpha$, completing the proof. Suppose $|M| < |A-M|$. Then by definition there is a proper subset $C$ of $A-M$ such that $|C| = |M|$. It is easy to see that \begin{align}\label{l13} |C| = |M| = |M \times M| = |M \times C| = |C \times M| = |C \times C| \end{align} and each of these sets are mutually disjoint. Consequently, by Lemma \ref{lem1.8}, $|\qty(M \times C) \cup \qty(C \times M) \cup \qty(C \times C)| = |M \times C| + |C \times M| + |C \times C| = |C|+|C|+|C| = |C|$ and there is a bijection $h:\qty(M \times C) \cup \qty(C \times M) \cup \qty(C \times C) \to C$. Let $E = \qty(M \times C) \cup \qty(C \times M) \cup \qty(C \times C)$ and define $\varphi: \qty(M \cup C) \times \qty(M \cup C) \to M \cup C$ by $\varphi|_{M \times M} = g$ and $\varphi|_{E} = h$. Since all sets in (\ref{l13}) are mutually disjoint, $\varphi$ is a well-defined bijection with $(g,M \times M) < \qty(\varphi,  \qty(M \cup C) \times \qty(M \cup C))$, contradicting the maximality of $g$. Thus, by Theorem \ref{thm1.5} and Lemma \ref{lem1.8}, $|A-M| \leq |M|$ and $|M| = |A-M|+|M| = |\qty(A-M) \cup M| = |A| = \alpha$, as desired. $\blacksquare$

Proposition \ref{prop1.7}-Theorem \ref{thm1.9} show a clear disconnect between infinity and finiteness in set-theory. In Proposition \ref{prop1.7}, if $A$ is finite and $A-F \neq \emptyset$, then clearly $|A \cup F| > |A|$. In Lemma \ref{lem1.8}, if $\alpha$ is finite the result is absurd if $\beta \neq 0$, and in Theorem \ref{thm1.9}, if $\alpha$ is finite then again the result if absurd if $\beta \geq 2$. This difference is important as it highlights that while infinite cardinals were defined as an extension of finite cardinals, they have fundamentally different properties, and properties true (or even trivial) for finite cardinals need not be true for infinite cardinals, and if they are true, a proof may not be nearly as trivial as in the finite case. Thus the three theorems demonstrate that infinite cardinals constitute a structure fundamentally distinct from that of finite cardinals, in spite of the fact of how they were defined, requiring new proof techniques such as Zorn's Lemma and AC along with refined intuition grounded in infinity.

The next result shows that, in fact, $\mbb{Q}$ and $\mbb{Z}$ are denumerable sets, which is quite unintuitive for $\mbb{N} \subsetneq \mbb{Z} \subsetneq \mbb{Q}.$ For a further proof of the following theorem, see \cite[Theorems 2.2-3]{kaplansky}.
\begin{proposition}\label{prop1.10}
$\mathbb{Q}$ and $\mathbb{Z}$ are denumberable.
\end{proposition}
\noindent $\textbf{Proof}$. Consider the map $f:\mathbb{Z} \to \mathbb{N}$ defined by $$f(n) = \begin{cases} 2n-1, \text{ if } n>0 \\ 
-2n, \text{ if } n \leq 0.
\end{cases}$$
Since the positive integers map to all odd naturals, and the non-positive integers to all even naturals, and clearly $f$ is injective, $f$ is a bijection and $|\mathbb{Z}| = \bm{\aleph}_{0}$. The map $\mathbb{Z} \to \mathbb{Q}$ by $n \mapsto \frac{n}{1}$ is clearly injective, and so $|\mathbb{Z}| \leq |\mathbb{Q}|$. Likewise, the map $\mbb{Q} \to \mbb{Z} \times \mbb{Z}$ by $\frac{n}{m} \mapsto \qty(\frac{n}{\gcd(n,m)}, \frac{m}{\gcd(n,m)})$ is clearly injective, and thus $|\mathbb{Z}| \leq |\mathbb{Q}| \leq |\mathbb{Z} \times \mathbb{Z}|$. Since $\mbb{Z}$ is denumerable, and the product of denumerable sets is denumerable by Theorem \ref{thm1.9}, $|\mbb{Z} \times \mbb{Z}| = \bm{\al_{0}}$. It then follows from Lemma \ref{lem1.4} that $|\mbb{Q}| = \bm{\al_{0}}$, as desired. $\blacksquare$ 
\begin{example}\label{ex3}
By the Archimedean principle (see the paragraphs preceding Theorem \ref{thm1.15}), there exists a rational numbers between any two real numbers. Due to this property, we say that $\mbb{Q}$ is \textit{dense} in $\mbb{R}$. Despite this, Proposition \ref{prop1.10} still shows $\mbb{Q}$ is denumerable. Since $\mbb{R}$ isn't denumerable by Theorem \ref{thm1.13}, one may expect that dense subsets of $\mbb{R}$ also aren't denumerable, but the rational numbers $\mbb{Q}$ show this isn't the case.
\end{example}
To this point, it has not been show that there are cardinals strictly larger than $\bm{\al}_{0}$. In fact, given any infinite cardinal $|A| = \alpha$, $\mcal{P}(A)$ (the set of all subsets of $A$) has a strictly larger cardinality, as is shown below. The result likely goes against ones natural intuition about infinity, for if infinity represents a process that goes on forever, how can one go past an infinity? If we instead think in terms of bijections instead of processes, it makes more sense that one can adjoin enough elements to a set in such a way that there is no bijection from the first set into the new set. Note, however, that by Proposition \ref{prop1.7}, in the infinite case, one cannot simply adjoin a finite number of elements, or even the cardinality of the original set by Lemma \ref{lem1.8}. For further proofs of the following Theorem, see \cite[Theorem 2.3.8]{jechintro}, \cite[Theorem 1.6.2]{abbott}, and \cite[Theorem 0.8.5]{hungerford}.
\begin{theorem}\label{thm1.11}
Given any set $A$, $\mcal{P}(A)$ has a strictly larger cardinality.
\end{theorem}
\noindent $\textbf{Proof}.$ The map $a \mapsto \{a\}$ is an injective map $A \to \mcal{P}(A)$ and thus $|A| \leq |\mcal{P}(A)|$. If there were a bijective map $f:A \to \mcal{P}(A)$, then for some $a_{0} \in A$, $f\qty(a_{0}) = B$, where $B = \{a \in A \mid a \notin f(a)\} \subset A$. But this yields a contradiction, for if $a_{0} \in B$, then $a_{0} \notin f(a_0) = B$, and if $a_{0} \notin B$, then $a_{0} \in f(a_0) = B$. Therefore, $|A| \neq |\mcal{P}(A)|$ and $|A| < |\mcal{P}(A)|$, as desired. $\blacksquare$ 

The real numbers $\mbb{R}$ are an example of an infinite set that isn't denumerable. In order to prove this result, a basic result from analysis, the nested interval property, will be used. First, a $\textit{sequence}$ in set $A$ is a function $f:\mbb{N} \to A$, denoted $\{a_{n}\}_{n \in \mbb{N}}$ ($n \mapsto a_{n}$). A sequence of intervals $\{\textit{I}_{n}\}_{n \in \mbb{N}}$ with $\textit{I}_{n} = [a_{n},b_{n}]$ and $a_{n},b_{n} \in \mbb{R}$ is $\textit{nested}$ if $\textit{I}_{n} \supset \textit{I}_{n+1}$ for every $n \in \mbb{N}$. Given any non-empty subset $A$ of $\mbb{R}$ with an upper bound, the $\textit{least upper bound property}$ guarantees the existence of a least upper bound of $A$, denoted $\sup{A}$. For further discussion of the following Lemma, see \cite[Theorem 1.4.1]{abbott}. The result relates to infinity through the use of intersections over countably infinite sets.
\begin{lemma}\label{lem1.12}
Let $\{\textit{I}_{n}\}_{n \in \mbb{N}}$ be a sequence of nested intervals $\textit{I}_{n} = [a_{n},b_{n}]$ with $a_{n}, b_{n} \in \mbb{R}$. Then $$\bigcap_{n \in \mbb{N}} \textit{I}_{n} \neq \emptyset.$$
\end{lemma}
\noindent $\textbf{Proof}.$ Consider the set $A = \{a_{j} \mid j \in \mbb{N}\}$. If there exists $j,n \in \mbb{N}$ such that $a_{j} > b_{n}$, then $\textit{I}_{j}$ and  $\textit{I}_{n}$ would be disjoint, so $a_{j} \leq b_{n}$ for every $j,n \in \mbb{N}$. Let $x = \sup{A}$. Since each $b_{n}$ is also an upper bound of $A$, $x \leq b_{n}$ for all $n \in \mbb{N}$, and by definition $a_{n} \leq x$ for every  $n \in \mbb{N}$. Thus $x \in \bigcap_{n \in \mbb{N}} \textit{I}_{n}$ (since it's in each $\textit{I}_{n}$) and $\bigcap_{n \in \mbb{N}} \textit{I}_{n} \neq \emptyset$, as desired. $\blacksquare$ 
\begin{theorem}\label{thm1.13}
The real numbers $\mbb{R}$ have a cardinality strictly larger than $\bm{\al}_{0}$.
\end{theorem}
\noindent $\textbf{Proof}.$ Clearly $|\mbb{N}| \leq |\mbb{R}|$. Suppose, by way of contradiction, there exists bijection $f:\mbb{N} \to \mbb{R}$. Let $x_{n} = f(n)$ for each $n$ and then we have $\mbb{R} = \{x_{n}\}_{n \in \mbb{N}}$. Let $\textit{I}_0$ be a (non-singleton) closed interval that does not contain $x_{0}$.  Define recursively $\textit{I}_{n}$ such that $\textit{I}_{n-1} \supset  \textit{I}_{n}$ and $x_{n} \notin \textit{I}_{n}$. By definition, given any $x_{j}\in \mbb{R}$, $x_{j} \notin \textit{I}_{j}$ and thus $x_{j} \notin \bigcap_{n \in \mbb{N}} \textit{I}_{n}$. Since all elements of $\mbb{R}$ are some $x_{j}$ ($j \in \mbb{N}$), it follows that $\bigcap_{n \in \mbb{N}} \textit{I}_{n} = \emptyset$, contradicting Lemma \ref{lem1.12}. Thus $|\mbb{N}| < |\mbb{R}|$, as desired. $\blacksquare$ 

For further discussion of the preceding theorem, see \cite[Theorem 1.5.6]{abbott} and \cite[Theorem 4.6.1]{jechintro}. The real numbers $\mbb{R}$ are often viewed by analysts as the completion of the rational numbers, in the sense that the real numbers are an ordered field with the least upper bound property and $\mbb{Q}$ as a subfield. A subset of $\mbb{Q}$ with an upper bound need not have a least upper bound, but the real numbers ``complete" the rational numbers by adjoining elements corresponding the least upper bounds of these sets. This ``completion" causes the real numbers to have a cardinality strictly larger than $|\mbb{Q}| = \bm{\al}_{0}$, and serves as a concrete, natural example of an infinity larger than $\bm{\al}_{0}$. This illustrates that infinity is not a single concept, but has many levels; $\mbb{R}$ marks the transition from a countable to an uncountable world. A rich question is then to ask: Is the cardinality of $\mbb{R}$ the first cardinal to leave this countable world? 

To prove the next result, again some basic analysis knowledge is used. The reader is likely familiar with these results, but they are given for completeness. A sequence of real numbers $\{a_{n}\}_{n \in \mbb{N}}$ $\textit{converges}$ to $L$, denoted $\lim_{n \to \infty} a_{n} = L$ or $a_n \to L$ if for every $\varepsilon > 0$, there exists $N \in \mbb{N}$ such that for every $n > N$, $|a_n-L| < \varepsilon$. We assume knowledge of basic properties of convergent sequences (see \cite[Theorem 3.3]{babyrudin}). Let $\textit{I}^{*}_{n} = \{0,1,2,...,n\}$ for each $n \in \mbb{N}$. A map $f:\textit{I}^{*}_{n} \to A$ for some $n \in \mbb{N}$ is called a $\textit{finite sequence}$ in $A$. Say a finite sequence $f:\textit{I}^{*}_{n} \to \mbb{R}$ is of the form $f(n) = ar^{n}$ ($a,r \in \mbb{R}, r \neq 1$). Let $S_{k} = \sum_{i = 0}^{k} f(i)$. Then $rS_{k} = \sum_{i = 1}^{k} f(i) + ar^{k+1}$ and by subtracting the two equations $S_{k} = \frac{a-ar^{k+1}}{1-r}$ ($0 \leq k \leq n$). The sum of a real sequence $\{a_{n}\}_{n \in \mbb{N}}$ is defined as $\lim_{k \to \infty} S_k$, where $S_k = \sum_{i = 0}^{k} a_i$. Now suppose $\{a_{n}\}_{n \in \mbb{N}}$ is a sequence of the form $a_n = ar^{n}$ with $|r| < 1$ ($a,r \in \mbb{R}$). Notice if $\varepsilon > 0$ and $r \neq 0$ ($r = 0$ is trivial), for all $k \in \mbb{Z}^{+}$ with $k > \lceil \frac{\ln \varepsilon}{\ln |r|} \rceil$, we have $k \ln|r| < \ln \varepsilon$ ($\ln|r| < 1$), whence $|r^k-0| < \varepsilon$ and $\lim_{k \to \infty} r^k = 0$. Thus \begin{align}\label{l14}
\sum_{n \in \mbb{N}} a_n = \lim_{k \to \infty} S_k = \lim_{k \to \infty} \frac{a-ar^{k+1}}{1-r} = \frac{a}{1-r}.
\end{align}
A sequence $\{a_{n}\}_{n \in \mbb{N}}$ is $\textit{bounded above}$ [resp. below] if it has an upper [resp. lower] bound. $\{a_{n}\}_{n \in \mbb{N}}$ is $\textit{increasing}$ [resp. \textit{decreasing}] if $a_{n} \leq a_{n+1}$ [resp. $a_{n} \geq a_{n+1}$]  for every $n \in \mbb{N}$. For further proofs of Lemma \ref{lem1.14}, see \cite[Theorem 3.14]{babyrudin} or \cite[Theorem 2.4.2]{abbott}.
\begin{lemma}\label{lem1.14}
Let  $\{a_{n}\}_{n \in \mbb{N}}$ be an increasing sequence of real numbers bounded from above. Then $\{a_{n}\}_{n \in \mbb{N}}$ converges.
\end{lemma}
\noindent $\textbf{Proof}.$ Let $\varepsilon > 0$ and $a = \sup\{a_{n}\}_{n \in \mbb{N}}$. Then there exists $N \in \mbb{N}$ such that $a_N > a-\varepsilon$ (otherwise $a-\varepsilon < a$ would be an upper bound). Since $\{a_{n}\}_{n \in \mbb{N}}$ is increasing, for every $n > N$, $a_n \geq a_N > a-\varepsilon$. Rearranging  gives $|a_n-a| = a-a_n < \varepsilon$, as desired. $\blacksquare$

Given $x,y \in \mbb{R}$ with $x > 0$, the $\textit{Archimedean principle}$ says there exists $n \in \mbb{N}$ such that $nx>y$. This follows by contradiction. Suppose there is $x,y \in \mbb{R}$ such that $y \geq nx$ for every $n \in \mbb{N}$. Then $A = \{nx \mid n \in \mbb{N} \}$ is bounded and has $\sup{A} = a$. Then $a-x < a$ so $a-x < mx$ for some $m \in \mbb{N}$, whence $a< (m+1)x \in A$, a contradiction. Using this, given $a < b$ in $\mbb{R}$ ($b<a$ case is analogous), $b-a > 0$ and choose $n \in \mbb{N}$ with $n(b-a) >1$. So there is $m \in \mbb{Z}$ with $na<m<nb$, yielding  $a<\frac{m}{n} < b$, and there exists a rational number between any two real numbers.

Currently, the theory is building towards the usual statement of the continuum hypothesis. To do this, it becomes necessary to show that $|\mbb{R}|$ (also called the cardinality of the continuum) is equal to $|\mcal{P}(\mbb{N})|$. To get some intuition for why this is true, imagine representing a real number $r \in \mbb{R}$ in binary, or more precisely, as a sum \begin{align}\label{l15}
r = \sum_{n \in \mbb{Z}} \delta_{r}(n)2^{n}
\end{align}
where $\delta_{r}(n)$ is a function taking on either one or zero depending on $n$ so that the sum is $r$. Then any real number can be constructed by making $|\mbb{Z}| = \bm{\al}_{0}$ choices as to whether $\delta_{r}(n)$ is one or zero. Likewise, a subset of $\mbb{Z}$ can constructed by making $|\mbb{Z}|$ choices as to whether an element is in the subset or not. So naturally one would (naively) expect $|\mbb{R}| = |\mcal{P}(\mbb{Z})| = |\mcal{P}(\mbb{N})|$. The issue with this intuition is that the binary representation of a real number need not be unique, so a more formal proof is given below. Overall, however, this connection between $\mbb{R}$ and $\mcal{P}(\mbb{N})$ illuminates the deep set-theoretic nature of the continuum. It highlights that the uncountability of $\mbb{R}$ is not only a result of its completion of the countable set $\mbb{Q}$, but also a manifestation of the higher cardinality that arises from considering sets of subsets of the countable set $\mbb{N}$. For further discussion of the following theorem (and Proposition \ref{prop1.16}), see \cite[Theorem 4.6.3]{jechintro}. Before the theorem, we define for a set $X$ and $A \subset X$ the \textit{characteristic function} $\chi_{A}:X \to \{0,1\}$ by \begin{align}\label{l105}
\chi_{A}(x) = \begin{cases}
1, \text{ if } x \in A \\
0, \text{ if } x \notin A.
\end{cases}
\end{align}
\begin{theorem}\label{thm1.15}
$\mcal{P}(\mbb{N})$ and $\mbb{R}$ have the same cardinality.
\end{theorem}
\noindent $\textbf{Proof}.$ Define $f: \mcal{P}\qty(\mbb{N}) \to\mbb{R}$ by $S \mapsto \sum_{n \in S} 10^{-n}$. Using (\ref{l14}) for $S = \mbb{N}$ gives $\sum_{n \in S} 10^{-n} = \frac{1}{1-\frac{1}{10}} = \frac{10}{9}$. If $S \subset \mbb{N}$ is finite, then clearly $\sum_{n \in S} 10^{-n}$ is a finite real number. If $S$ is infinite, then the sequence of partial sums is increasing and bounded above by $\frac{10}{9}$, and thus convergent by Lemma \ref{lem1.14}. So again, $\sum_{n \in S} 10^{-n}$ is a defined real number and $f$ is well-defined. If $S_{1},S_{2} \in \mcal{P}\qty(\mbb{N})$ with $S_{1} \neq S_{2}$, then there exists $n \in \mbb{N}$ in only one subset. Let $k \in \mbb{N}$ be the least such element. Without loss of generality, assume $k \in S_{1}$ and $k \notin S_{2}$. \newpage We have \begin{align}
f(S_1)-f(S_2) & = \sum_{n \in S_1} 10^{-n}-\sum_{n \in S_2} 10^{-n} = 10^{-k}+\sum_{n > k} (\chi_{S_1}(n)-\chi_{S_2}(n))10^{-n}\label{l101} \\
& \geq 10^{-k}-\sum_{n > k} 10^{-n}\label{l102} \\ & = 10^{-k}-\frac{10^{-(k+1)}}{1-\frac{1}{10}}\label{l103} \\ & = \frac{8}{9} 10^{-k}\label{l104}
\end{align} 
where (\ref{l103}) follows from (\ref{l14}). Thus $f(S_1)-f(S_2) > 0$ and $f(S_1) \neq f(S_2)$, so $f$ is injective and $|\mcal{P}\qty(\mbb{N})| \leq |\mbb{R}|$. By Proposition \ref{prop1.10}, there exists a bijection $f: \mbb{Q} \to \mbb{N}$. Then the map $g:\mcal{P}\qty(\mbb{Q}) \to \mcal{P}\qty(\mbb{N})$ by $S \mapsto f(S)$  is clearly bijective and $\abs{\mcal{P}\qty(\mbb{Q})} = \abs{\mcal{P}\qty(\mbb{N})}$. Define $h:\mbb{R} \to \mcal{P}\qty(\mbb{Q})$ by $r \mapsto \{q \in \mbb{Q} \mid q < r\}$. Clearly $h$ is well-defined. Moreover, if $r_1, r_2 \in \mbb{R}$ are such that $r_1 \neq r_2$, without loss of generality, assume $r_1 < r_2$. Then there exists a rational number $s$ such that $r_1 < s < r_2$, and $s \in h(r_2)$ but $s \notin h(r_1)$. Thus $h(r_1) \neq h(r_2)$ and $h$ is injective, so $\mbb{R} \xrightarrow{h} \mcal{P}\qty(\mbb{Q}) \xrightarrow{g} \mcal{P}\qty(\mbb{N})$ is an injection and $|\mbb{R}| \leq \abs{\mcal{P}\qty(\mbb{N})}$. Therefore, by Lemma \ref{lem1.4}, $\abs{\mbb{R}} = \abs{\mcal{P}\qty(\mbb{N})}$, as desired. $\blacksquare$ 

Given two cardinal numbers $|A| = \alpha$ and $|B| = \beta$, the cardinal $\alpha^{\beta}$ is defined to be the cardinality of the set of all functions $B \to A$. It is easy to see this is independent of the sets $A$ and $B$ chosen. This should make some intuitive sense, for given any $b \in B$, there are $|A| = \alpha$ ``choices" to map $b$ to in a function $B \to A$. Given a set $A$, we could view a subset $S$ of $A$ as a map $A \to \{0,1\}$, where $a \mapsto 1$ means $a \in S$ and $a \mapsto 0$ means $a \notin S$ (i.e., the characteristic function $\chi_{S}$). So intuitively, we have 
\begin{proposition}\label{prop1.16}
For any set $A$, $|\mcal{P}(A)| = 2^{|A|}$.
\end{proposition}
\noindent $\textbf{Proof}.$ Clearly $|\{0,1\}| = 2$ and so $2^{|A|}$ is the cardinality of the set $F\qty(A,\{0,1\})$ of all functions $A \to \{0,1\}$. Consider the map $f:\mcal{P}(A) \to F\qty(A,\{0,1\})$ by $S \mapsto \chi_{S}$ ($S \in \mcal{P}(A)$). Clearly $f$ is surjective as any map in $F\qty(A,\{0,1\})$ will have some subset of $A$ map to $1$ and the rest $0$. Moreover, if $B \neq C$ with $B,C \subset A$, then there is some $x$ in only one subset, say, without loss of generality, $x \in B$ and $x \notin C$. Then $\chi_{B}(x) = 1 \neq 0 = \chi_{C}(x)$ and thus $\chi_{B} \neq \chi_{C}$, whence $f$ is injective. Thus $f$ is bijective and $\abs{\mcal{P}(A)} = \abs{F\qty(A,\{0,1\})} = 2^{|A|}$, as desired. $\blacksquare$ 

The theory of cardinals for this paper is complete. By Theorem \ref{thm1.5}, the class of all cardinals $\mathfrak{C}$ is well-ordered by the ordering defined above. Let $\bar{\alpha} = \{\beta \in \mathfrak{C} \mid \beta \leq \alpha\}$ for any cardinal $\alpha$. Theorem \ref{thm1.11} shows that $\mathfrak{C} \setminus \bar{\alpha}$ is non-empty for any $\alpha \in \mathfrak{C}$. Moreover, by Theorem \ref{thm1.5}, each subset $\mathfrak{C} \setminus \overline{\alpha}$ must have a least element. Let $\{a_{n}\}_{n \in \mbb{N}}$ be a sequence in $\mathfrak{C}$ defined by $a_{0} = \bm{\al}_{0}$ and $a_{n} = \min\{\mathfrak{C} \setminus \overline{a_{n-1}}\}$. Then  $\{a_{n}\}_{n \in \mbb{N}}$ is a strictly increasing sequence of infinite cardinal numbers. Let $a_{j} = \bm{\al}_{j}$ for each $j \in \mbb{N}$ and we have \begin{align}\label{l16}
\bm{\al}_{0}<\bm{\al}_{1}<\bm{\al}_{2}<\bm{\al}_{3}<...
\end{align}
Theorem \ref{thm1.5} shows there exists no cardinal number $\beta$ such that $\bm{\al}_{j} < \beta < \bm{\al}_{j+1}$ ($j \in \mbb{N}$), for if there was, $\beta < \bm{\al}_{j+1} = \min\{\mathfrak{C} \setminus \overline{\bm{\al}_{j}}\}$ and $\beta \in \mathfrak{C} \setminus \overline{\bm{\al}_{j}}$, contradicting the minimality of $\bm{\al}_{j+1}$. The $\textit{continuum hypothesis}$ is the statement that there is no cardinal number $\beta$ such that $|\mbb{N}| < \beta < |\mbb{R}|$. Theorem \ref{thm1.15} shows $|\mcal{P}(\mbb{N})| = |\mbb{R}|$ and Proposition \ref{prop1.16} shows that $|\mcal{P}(\mbb{N})| = 2^{|\mbb{N}|}$. Moreover, for there to exist no such $\beta$, we must have $|\mbb{R}| = \bm{\al}_{1}$ (otherwise $\beta = \bm{\al}_{1}$ holds). So the continuum hypothesis (CH) can be rewritten in the familiar form \begin{align}\label{l17}
2^{\bm{\al}_{0}} = \bm{\al}_{1}.
\end{align}
Or equivalently, define $\{\bm{\beth}_{n}\}_{n \in \mbb{N}}$ by $\bm{\beth}_{0} = \bm{\al}_{0}$ and $\bm{\beth}_{n} = 2^{\bm{\beth}_{n-1}}$. Then \begin{align}\label{l18}
\bm{\beth}_{1} = \bm{\al}_{1}.
\end{align}
Gödel \cite{godel1} showed that the continuum hypothesis is consistent with the other Zermelo-Frankel axioms (ZFC). Later, Cohen \cite{cohen} showed that the continuum hypothesis is independent of them, meaning the usual axioms of ZFC aren't enough to determine the truth value of the statement. That is, the statement is undecidable. This, however, does not mean the problem is resolved. The following excerpt is from Gödel \cite{godel2}, on the matter of a proof of the undecidability of CH.
\begin{center}
\begin{quote}
\emph{"Only someone who (like the intuitionist) denies that the concepts and axioms of
classical set theory have any meaning (or any well-defined meaning) could be
satisfied with such a solution, not someone who believes them to describe some
well-determined reality. For in this reality [CH] must be either
true or false, and its undecidability from the axioms as known today can only
mean that these axioms do not contain a complete description of this reality;
and such a belief is by no means chimerical, since it is possible to point out ways
in which a decision of the question, even if it is undecidable from the axioms in
their present form, might nevertheless be obtained."}
\end{quote}
\end{center}
This excerpt reflects a deep philosophical stance on the limitations of mathematical truth. Gödel argues that the undecidability of CH within ZFC is inherently unsatisfying to the many mathematicians who view set-theoretic axioms as descriptions of an objective, well-defined mathematical reality. According to Gödel, if one believes that the axioms, and therefore the sets they govern, refer to a real, determinate universe, then CH must have a definite truth value, regardless of the current inability of our axiomatic systems to prove it. In this view, it is not a reflection of indeterminacy in mathematics itself, but merely the incompleteness of existing axioms. Gödel contrasts this point of view with Brouwer's intuitionism, where the theory of infinite ordinals (see section \ref{sub1.2}) greater than $\omega_{1}$ is rejected as meaningless. Intuitionists may be satisfied with the undecidability of CH from the inherently destructive viewpoint. For them, undecidability is not an issue of incompleteness of axioms, unlike perhaps a Platonist or realist, but a reflection of the limitations of formalized mathematics.

This discussion of the continuum hypothesis illustrates one of the most profound intersections between infinity and the foundations of mathematics. The undecidability of CH in ZFC highlights the limitations of the structure of infinite sets within them---particularly with the existence of sets with cardinalities between $\mbb{N}$ and $\mbb{R}$. If they exist, then the axioms fail to capture such sets; if they don't, the axioms fail to prove otherwise. If we wish to keep ZFC (or ZF) as the core of set theory, unless one (like the intuitionist) rejects the notion that mathematics describes an objective, well-determined mathematical reality and holds the position above, or one rejects the law of excluded middle, this can only be resolved in one of two ways (similar to AC). We accept CH, or we reject it. Cohen himself expresses ``A point of view which the author feels may eventually come to be
accepted is that CH is obviously false" \cite[Section 13]{cohen}. The author holds this position as well.

Regardless of whether one accepts CH, rejects it, or believes that it is truly undecidable, it highlights that infinity in this context is not merely a mathematical tool, but a philosophical boundary. It demonstrates the limits of formal reasoning and compels one to confront the nature of mathematical truth itself. The existence of undecidable statements like CH in ZFC ultimately emphasizes the depth, complexity, and richness of infinity within the foundations of mathematics.
\subsection{Ordinal Numbers}\label{sub1.2}
Proposition \ref{prop1.10} showed that $|\mbb{N}| = |\mbb{Z}| = \bm{\al}_{0}$, yet $\mbb{N} \subsetneq \mbb{Z}$. So what is the difference between these sets? One major difference is their $\textit{order type}$. To define this term, some theory needs to be developed. Given two linearly ordered sets $(A, \leq_{A})$ and $(B,\leq_{B})$, a mapping $f:A \to B$ is an $\textit{order isomorphism}$ if $f$ is a bijection and for every $a_1,a_2 \in A$ \begin{align}\label{l19}
    a_1 \leq_{A} a_2 \iff f(a_1) \leq_{B} f(a_2).
\end{align}
Two linearly ordered sets are $\textit{order isomorphic}$ (denoted $(A, \leq_{A}) \cong (B,\leq_{B})$)  if there exists an order isomorphism between them.
\begin{example}\label{ex6}
The map $f(x) = \tan(\pi x-\pi/2)$ is a strictly increasing bijection from $(0,1)$ onto $\mbb{R}$, and hence $(0,1) \cong \mbb{R}$ (under usual orderings).
\end{example}
\begin{example}\label{ex7}
The set $\{0,1,2\} \subset \mbb{N}$ ordered by $a \leq b$ if and only if $a = b$ or $a \in b$ is order isomorphic to $\{x,y,z\}$ with $x < y < z$. More generally, any finite well-ordered set of cardinality $n+1$ is order isomorphic to $\{0,1,2,...,n\}$
\end{example}
The following theorem is stated in \cite[Section 2.6]{stoll}, but a proof is given below.
\begin{theorem}\label{thm1.17}
The relation defined by $(A, \leq_{A}) \sim (B,\leq_{B})$ if and only if there exists an order isomorphism $f:A \to B$ is an equivalence relation on the class of all linearly ordered sets.
\end{theorem}    
\noindent $\textbf{Proof}.$ Clearly $(A, \leq_{A}) \leq (A, \leq_{A})$ for any linearly ordered set $(A, \leq_{A})$, for the identity map $f:A \to A$ is clearly an order isomorphism. Now suppose $(A, \leq_{A}) \sim (B,\leq_{B})$ with order isomorphism $f:A \to B$. Since $f:A \to B$ is a bijection, $f^{-1}:B \to A$ is a bijection (Theorem \ref{thm1.1}). Write $b_1 = f(a_1)$ and $b_2 = f(a_2)$, and notice \begin{align}\label{l20}
    b_1 \leq_{B} b_2 \iff f(a_1) \leq_{B} f(a_2) \iff a_1 \leq_{A} a_2 \iff f^{-1}(b_1) \leq_{A} f^{-1}(b_2).
\end{align}
Thus $f^{-1}$ is an order isomorphism and $(B,\leq_{B}) \sim (A, \leq_{A})$. If $(A, \leq_{A}) \sim (B,\leq_{B})$ and $(B,\leq_{B}) \sim (C, \leq_{C})$ with order isomorphisms $f:A \to B$ and $g:B \to C$, then the proof of Theorem \ref{thm1.2} shows that the composition $gf:A \to C$ is a bijection. Since $f$ and $g$ are order isomorphisms, \begin{align}\label{l21}
    a_1 \leq_{A} a_2 \iff f(a_1) \leq_{B} f(a_2) \iff gf(a_1) \leq_{C} gf(a_2)
\end{align}
 and $gf$ is an order isomorphism, as desired. $\blacksquare$
 \begin{definition}\label{def2.1}
The order type of a linearly ordered set $(A, \leq_{A})$, denoted $\text{ord}(A, \leq_{A})$, is the equivalence class of $(A, \leq_{A})$ under the equivalence relation defined in Theorem \ref{thm1.17}. When the ordering is clear by context, the order type of $A$ will be denoted simply as $\text{ord} A$
 \end{definition}
\begin{example}\label{ex8}
The set $\{0,1,2,...,n-1\} \subset \mbb{N}$ has order type $n$ (under usual ordering). The natural numbers $\mbb{N}$ have the same order type as the set $2\mbb{N}$, and more generally, as a clear consequence of the definition, two sets have the same order type if and only if they are order isomorphic (see (\ref{l35})).
\end{example}
Theorem \ref{thm1.17} shows this definition is well-defined. Order types can be used to show the differences between sets of the same cardinality. For example, say $A = \{a_1,...,a_n\}$ ($a_i$ distinct) is a set with $A \cap \mbb{N} = \emptyset$. Then Proposition \ref{prop1.7} shows that $|A \cup \mbb{N}| = |\mbb{N}|$. However, consider the ordering on $A \cup \mbb{N}$ defined by $n < a_j$ for every $n \in \mbb{N}$ and $a_j \in A$, $a_1<...<a_n$ the ordering in $A$, and the usual ordering in $\mbb{N}$. Then an order isomorphism must map $0 \mapsto 0$, as it is the least elements of both sets. Suppose then that $n > 0$ and an order isomorphism must map $k \mapsto k$ ($k <n$). Then $n$ cannot map to a smaller element or the map isn't injective, and it cannot map to a larger element or some other element larger than $n$ must map to $n$, so by induction $n \mapsto n$ for every $n \in \mbb{N}$. However, then each $a_j$ would map to an element already mapped to, and the map would not be a bijection. Thus $A \cup \mbb{N}$ and $\mbb{N}$ are not order isomorphic. So while these two infinite sets may have the same cardinality, there is an important difference between them, primarily their order type. This analysis demonstrates how infinity in set-theoretic mathematics manifests in more than size alone, order plays a crucial role. Even when two sets share the same cardinality, their differing order types highlights a richer structure underlying infinite sets.

We now define the ordinal number of a well-ordered set. One definition is to define the order type of a well-ordered set to be its ordinal number. This definition works for some purposes, but its main downfall is that, as with the definition of cardinality above, the ordinal number of a set would not be a set, but rather a proper class. Instead, a singular set from the order type of a well-ordered set will be chosen as its ordinal number. Recall from above the construction of $\mbb{N}$
\begin{align*}
0  & \coloneqq \emptyset \\
1 & \coloneqq \{\emptyset\} \\
2 & \coloneqq \{\emptyset, \{\emptyset\}\} \\
3 & \coloneqq \{\emptyset, \{\emptyset\},\{\emptyset, \{\emptyset\}\}\} \\
\vdots
\end{align*}
The sets above correspond to the respective finite ordinal numbers. Notice that in each of the ordinal numbers above, the set is well ordered with respect to set membership (i.e., $a \leq b$ if and only if $a = b$ or $a \in b$), and that every element of the ordinal is also a subset of the ordinal, which motivates 
\begin{definition}\label{def2.2}
A set $\alpha$ is an ordinal number if $\alpha$ is well-ordered with respect to $\leq$ defined by $$a \leq b \iff a = b \lor a \in b$$
and every element of $\alpha$ is also a subset of $\alpha$.
\end{definition}

The concept of ordinal numbers and order types illustrate a relationship between how finite and infinite sets operate, rather than a disconnect. Cardinal numbers of infinite sets and their properties as shown above differed quite substantially from what is expected in a finite world. Ordinal numbers and order types capture structural distinctions that cardinality alone cannot. For instance, consider a well-ordered infinite set $A$, and $A$ with finitely many elements adjoined to the end of it, or even $|A|$ elements, or $|A|$ elements adjoined to the beginning of $A$. In the finite case, such modifications produce fundamentally different sets, but in the infinite case, these differences are obscured by the fact that the sets share the same cardinality. The order type and order isomorphic ordinal numbers (or lack thereof) of such sets when order is defined in a natural way are quite different however, showcasing the clear differences in the structure of the above sets from $A$ that are hidden when only considering cardinality, and highlighting how orders behave in finite versus infinite contexts. This underscores that infinity is not monolithic. While cardinal numbers view all sets of the same cardinality as equipollent, order types and numbers reveal a deeper nature of infinity, one concerning not just size, but also the structural behavior of underlying sets.
 
It turns out that every well-ordered set is order isomorphic to precisely one ordinal number. To prove this, however, requires some theory. First off, the class of all ordinal numbers is denoted by $\text{Ord}$. Given $x,y \in \text{Ord}$, define $x \leq y$ if and only if $x = y$ or $x \in y$. For further discussion of the following Lemma, see \cite[Lemma 2.11]{jech}. Note that the $\textit{axiom of regularity}$ \cite[Section 2.9]{suppes} guarantees that any nonempty set must contain an element disjoint to it, and applying it to $\{A\}$ shows that $A \notin A$ for any non-empty set $A$ (otherwise, $A \cap \{A\}$ would contain $A$, which isn't empty). Note that $\emptyset \notin \emptyset$ and so this is true for any set.
\begin{lemma}\label{lem1.18}
Let $\alpha$ be an ordinal. \\
(i) If $\beta \in \alpha$, then $\beta$ is an ordinal. \\
(ii) If $\beta$ is an ordinal with $\alpha \neq \beta$ and $\alpha \subset \beta$, then $\alpha \in \beta$. \\
(iii) If $\beta$ is an ordinal, then either $\alpha \subset \beta$ or $\beta \subset \alpha$. \\
(iv) The class of all ordinals is linearly ordered by $\leq$.
\end{lemma}
\noindent $\textbf{Proof}.$ [(i)] Consider $\beta \in \alpha$. If $\beta = \emptyset$, it is (vacuously) an ordinal. If not, by definition $\beta \subset \alpha$, and any nonempty subset of $\beta$ is a subset of $\alpha$ and thus has a least element, so $\beta$ is well-ordered. Now consider $\gamma \in \beta$. Let $\zeta \in \gamma$ (if $\gamma = \emptyset$ by above it's an ordinal). Then $\beta \subset \alpha$ and so $\gamma \in \alpha$, and likewise $\zeta \in \alpha$. Then $\zeta \in \gamma \land \gamma \in \beta$, and by transitivity $\zeta \in \beta$, so $\gamma \subset \beta$ and $\beta$ is an ordinal.

[(ii)] If $\alpha \subsetneq \beta$, let $\gamma$  be the least element of the set $\beta \setminus \alpha$. Then, all elements smaller than $\gamma$ are elements of $\alpha$. If there is some element $\zeta \in \alpha$ with $\gamma \leq \zeta$, then either $\gamma = \zeta$ or $\gamma \in \zeta \land \zeta \subset \alpha$, so $\gamma \in \alpha$. Either case is a contradiction and so $\alpha = \{\varphi \in \beta \mid \varphi < \gamma\} =  \{\varphi \in \beta \mid \varphi \in \gamma\} = \gamma$, thus $\alpha \in \beta$.

[(iii)] Consider $\alpha \cap \beta$. If $\alpha \cap \beta = \emptyset$, then it is (vacuously) an ordinal. Otherwise, since it is a subset of $\alpha$, it's well-ordered by memberhsip, and if $\gamma \in \alpha \cap \beta$, then $\gamma \in \alpha \land \gamma \in \beta$ and so $\gamma \subset \alpha \land \gamma \subset \beta$, thus $\gamma \subset \alpha \cap \beta$ and $\alpha \cap \beta$ is an ordinal. Let $\zeta = \alpha \cap \beta$. If $\zeta \neq \alpha$ and $\zeta \neq \beta$, then by (ii) $\zeta \in \alpha$ and $\zeta \in \beta$, so $\zeta \in \zeta$, contradicting the axiom of regularity. So $\zeta = \alpha$ or $\zeta = \beta$, and $\alpha \subset \beta$ or $\beta \subset \alpha$. 

[(iv)] Clearly $x \leq x$ for any $x \in \text{Ord}$. If $x,y,z \in \text{Ord}$ are such that $x \leq y$ and $y \leq z$ and $x = y$ or $y = z$, then transitivity is immediate. Otherwise, $x < y$ and $y < z$ and so $(x \in y) \land (y \subset z)$, which implies $x \in z$ and $x < z$. If $x,y \in \text{Ord}$ are such that $x \leq y$ and $y \leq x$, then $x \subset y$ and $y \subset x$ and so $x = y$. If $x,y \in \text{Ord}$, then by (iii), $x \subset y$ or $y \subset x$, and then by (ii) $x = y$ or $(x \in y) \lor (y \in x)$, in any case $x \leq y$ or $y \leq x$, as desired. $\blacksquare$.
\begin{theorem}\label{thm1.19}
The class of all ordinals is well-ordered by $\leq$.
\end{theorem}
\noindent $\textbf{Proof}.$ Lemma \ref{lem1.18} (iv) already shows $\leq$ is a linear ordering. Let $\{\alpha_{i}\}_{i \in \textit{I}} \subset \text{Ord}$ be a non-empty family of ordinals. Consider $\bigcap_{i \in \textit{I}} \alpha_{i}$. Given any $\alpha \in \{\alpha_{i}\}_{i \in \textit{I}}$, $\bigcap_{i \in \textit{I}} \alpha_{i} \subset \alpha$ and $\bigcap_{i \in \textit{I}} \alpha_{i}$ is well-ordered. If $\alpha \in \bigcap_{i \in \textit{I}} \alpha_{i}$,then $\alpha \in \alpha_{i}$ for every $i \in \textit{I}$, and so $\alpha \subset \alpha_{i}$ for every $i \in \textit{I}$ and thus $\alpha \subset \bigcap_{i \in \textit{I}} \alpha_{i}$, so $\bigcap_{i \in \textit{I}} \alpha_{i}$ is an ordinal. For any $j \in \textit{I}$, $\bigcap_{i \in \textit{I}} \alpha_{i} \subset \alpha_{j}$, and then by Lemma \ref{lem1.18} (ii), $\bigcap_{i \in \textit{I}} \alpha_{i} \leq \alpha_{j}$ for any $j \in \textit{I}$. Now all that's left to show is that $\bigcap_{i \in \textit{I}} \alpha_{i} \in  \{\alpha_{i}\}_{i \in \textit{I}}$.

Suppose, by way of contradiction, that $\bigcap_{i \in \textit{I}} \alpha_{i} \notin  \{\alpha_{i}\}_{i \in \textit{I}}$. Then for every $\alpha_{i}$ ($i \in \textit{I}$), $\bigcap_{i \in \textit{I}} \alpha_{i} \neq \alpha_{i}$. Thus, $\bigcap_{i \in \textit{I}} \alpha_{i} < \alpha_{i}$ for every $i \in \textit{I}$. Note that \begin{align}\label{l22}
\bigcap_{i \in \textit{I}} \alpha_{i} = \left\{x\;\middle|\;x \in \alpha_{i} \hspace{2mm} \forall i \in \textit{I}\right\} = \left\{x\;\middle|\;x< \alpha_{i}\forall i \in \textit{I}\right\}.
\end{align}
Thus it follows that $\bigcap_{i \in \textit{I}} \alpha_{i} \in \bigcap_{i \in \textit{I}} \alpha_{i}$, a contradiction of the axiom of regularity. $\blacksquare$

Theorem \ref{thm1.19} is the ordinal analogue of Theorem \ref{thm1.5}. In contrast to Theorem \ref{thm1.5}, however, it does not require AC. This distinction underscores the order-theoretic nature of ordinals in comparison to cardinals. While cardinals and ordinals both share this property, the statement for cardinals requires a notion of uncountably infinitely many choices, while at most the ordinal statement would require finitely many choices. This highlights the nuanced nature of infinity, where in some contexts it is well-behaved, but in others, dependent on strong axioms like AC.

Given a well ordered set $(S, \leq)$, define the $\textit{segment}$ of $x \in S$ by $S_{x} = \{a \in S \mid a<x\}$. Consider $\qty(S_{y})_{x}$ with $x < y$. Then clearly $\qty(S_{y})_{x} \subset S_{x}$. Moreover, any $z \in S_{x}$ has $z<x<y$ and so $z \in S_{y}$ with $z < x$ and $S_{x} \subset \qty(S_{y})_{x}$, whence $\qty(S_{y})_{x} = S_{x}$.  Let $(S', \leq')$ be well-ordered and $(S, \leq) \cong (S', \leq')$ with order isomorphism $f:S \to S'$. Then 
\begin{align}\label{l23}
f(S_{x}) = \left\{f(a) \in S' \;\middle|\; a<x\right\} = \left\{y \in S'\;\middle|\; y<f(x)\right\} =  S'_{f(x)}.
\end{align}
\begin{proposition}\label{prop1.20}
Let $\alpha$ and $\beta$ be ordinals. \\
(i) If $\alpha < \beta$, then $\beta_{\alpha} = \alpha$ \\
(ii) If $\alpha \cong \beta$, then $\alpha = \beta$.
\end{proposition}
\noindent $\textbf{Proof}.$ [(i)] Since all ordinals smaller than $\alpha$ are also smaller than $\beta$, \begin{align}\label{l24}
\beta_{\alpha} = \left\{x \in \beta \;\middle|\; x<\alpha\right\} = \left\{x \in \text{Ord} \;\middle|\; x < \alpha\right\} = \alpha,
\end{align}
where the last equality follows from Lemma \ref{lem1.18} (i). 

[(ii)] Suppose, by way of contradiction, that $\alpha \neq \beta$. By Lemma \ref{lem1.18} (iv) either $\alpha \in \beta$ or $\beta \in \alpha$. Without loss of generality, assume $\alpha \in \beta$. Let $f:\beta \to \alpha$ be an order isomorphism. Then $f(\alpha) \in \alpha$ and $f(\alpha) < \alpha$. Define $f^{n}$ to be $f$ iterated $n$ times and $f^{0} = 1_{\beta}$, then a simple induction argument shows $f^{n+1}(\alpha) < f^{n}(\alpha)$ for any $n \in \mbb{N}$. Thus $\{f^{n}(\alpha) \mid n \in \mbb{N}\}$ is a strictly decreasing sequence in $\alpha$, which cannot have a least element, contradicting the well-ordering of $\alpha$. Thus $\alpha = \beta$, as desired. $\blacksquare$
\begin{example}\label{ex9}
As a consequence of Lemma \ref{lem1.18}(iii) and Proposition \ref{prop1.20}(ii), we see that for each finite cardinal $n \in \mbb{N}$, the set $A = \{0,1,2,...,n-1\}$ is the unique ordinal number of cardinality $n$. To prove this, notice if $\beta$ is an ordinal of cardinality $n$, we have $\beta \subset A$ or $A \subset \beta$. Without loss of generality, assume $\beta \subset A$. Then the inclusion $\beta \xrightarrow{\subset} A$ is an injective order preserving map with $|\beta| = |A|$, so by the pigeonhole principle (see the proof of Theorem \ref{thm2.18}), the map is bijective and an order isomorphism, so $\beta \cong A$ and thus by Proposition \ref{prop1.20}(ii), $\beta = A$. It may seem one could create another ordinal with replacement, but the properties that the set must not only be well-ordered by membership, but also having each element be a subset, makes this not possible. This illustrates how finite ordinals behave in a rigid, well-determined way. In contrast, infinite ordinals are more complex, as infinite orders can differ in structure despite having the same cardinality, and highlights the more subtle nature of infinity in the theory of ordinals.
\end{example}
\begin{lemma}\label{lem1.21}
Let $(S, \leq)$ be a well-ordered set with $T \subset S$ and $f:S \to T$ an order isomorphism. Then for every $x \in S$, $x \leq f(x)$.
\end{lemma}
\noindent $\textbf{Proof}.$ Let $A = \{x \in S \mid f(x) < x\}$. Suppose, by way of contradiction, that $A \neq \emptyset$. Then $A$ has a minimal element, say $x_{0}$, so $f(x_{0}) < x_{0}$. Let $x_{1} = f(x_{0})$, and since $f$ is order preserving, $f(x_{1}) < f(x_{0}) = x_{1}$, and $x_{1} < x_{0} \in A$, contradicting the minimality of $x_0$. Thus $A = \emptyset$ and $x \leq f(x)$ for every $x \in S$, as desired. $\blacksquare$
\begin{theorem}\label{thm1.22}
Let $(S, \leq)$ and $(S', \leq')$ be well-ordered sets with $(S, \leq) \cong (S', \leq')$. Then there exists a unique order isomorphism $f:S \to S'$.
\end{theorem}
\noindent $\textbf{Proof}.$ Let $f:S \to S'$ and $g: S \to S'$ be order isomorphisms. Consider $h = f^{-1} g$, which is an order isomorphism by the proof of Theorem \ref{thm1.17}. Then by Lemma \ref{lem1.21}, for every $x \in S$, $x \leq h(x)$. Apply $f$ to both sides to obtain $f(x) \leq' g(x)$ for every $x \in S$. Similarly, consider $\varphi = g^{-1} f$. Likewise, $x \leq \varphi(x)$ for every $x \in S$, and applying $g$ to both sides gives $g(x) \leq' f(x)$ for every $x \in S$. Thus $f = g$, as desired. $\blacksquare$
\begin{lemma}\label{lem1.23}
Let $(S, \leq)$ be a well-ordered set. If for every $x \in S$, $S_{x}$ is isomorphic to an ordinal number, then $S$ is isomorphic to an ordinal number.
\end{lemma}
\noindent $\textbf{Proof}.$ For every $x \in S$, let $f_{x}:S_{x} \to f(x) \in \text{Ord}$ be an order isomorphism (AC not required, map is unique by Theorem \ref{thm1.22}). Consider $\alpha = \{f(x) \mid x \in S\}$. Since each $f(x) \in \alpha$ is an ordinal, by Theorem \ref{thm1.19} $\alpha$ is well-ordered. Let $\gamma \in f(x)$, and then $\gamma = f_{x}(y)$ for some $y \in S_{x}$. Since $f_{x}(y) < f(x)$, by Proposition \ref{prop1.20} \begin{align}\label{l25}
f_{x}(y) = \qty(f(x))_{f_{x}(y)} = \left\{f_{x}(z)\;\middle|\;f_{x}(z) < f_{x}(y)\right\} = \left\{f_{x}(z)\;\middle|\;z<y\right\}.
\end{align}
The map $f_{x}(y) \to S_{y}$ by $f_{x}(z) \mapsto z$ is then an order isomorphism (since $f_{x}$ is) and 
\begin{align}\label{l26}
f_{x}(y) \cong S_{y} \cong f(y).
\end{align}
Thus $\gamma = f(y)$ by Proposition \ref{prop1.20}(ii) and $f(x) \subset \alpha$, so $\alpha$ is an ordinal. Consider the map $f:S \to \alpha$ by $x \mapsto f(x)$. Suppose $x,y \in S$ with $x<y$. The maps $f_{x}^{-1}: f(x) \to S_{x}$ and $f_y|_{S_x}: \qty(S_{y})_{x} \to \qty(f(y))_{f_y(x)}$  are order isomorphisms. Consider the map $\varphi = f_{y}|_{S_{x}} f_{x}^{-1}:f(x) \to \qty(f(y))_{f_{y}(x)}$. Since $\varphi$ is an order isomorphism (Theorem \ref{thm1.17}), then by Proposition \ref{prop1.20} $f(x) = f_{y}(x) < f(y)$ (by \ref{l26}). So $x \leq y$ implies $f(x) \leq f(y)$. If $x \neq y$, then $x < y$ or $y < x$. Then $f(x) < f(y)$ or $f(y) < f(x)$ and $f(x) \neq f(y)$, so $f$ is injective. Suppose $f(x) \leq f(y)$. Then by comparability, $x \leq y$ or $y \leq x$. If $y < x$, then $f(y) \leq f(x)$. By anti-symmetry, $f(y) = f(x)$, contradicting the injectivity of $f$. Thus $x \leq y$ and $f$ is order preserving. Clearly $f$ is surjective, thus $S \cong \alpha$, as desired. $\blacksquare$
\begin{theorem}\label{thm1.24}
Every well-ordered set is isomorphic to a unique ordinal number.
\end{theorem}
\noindent $\textbf{Proof}.$ Let $(S, \leq)$ be a well-ordered set and $E = \{a \in S \mid S_{a} \ncong \beta \hspace{2mm} \forall \beta \in \text{Ord} \}$. Suppose, by way of contradiction, that $E \neq \emptyset$. Let $\gamma = \min E$. Then for any $x < \gamma$, $S_{x}$ is order isomorphic to some ordinal. Since $S_{x} = \qty(S_{\gamma})_{x}$ for any $x \in S_{\gamma}$, which is necessarily less than $\gamma$, $\qty(S_{\gamma})_{x}$ is isomorphic to an ordinal, and then by Lemma \ref{lem1.23}, $S_{\gamma}$ is isomorphic to some ordinal, a contradiction. Thus $E = \emptyset$, and by Lemma \ref{lem1.23}, $S$ is isomorphic to some ordinal, say $\alpha$. If $S \cong \alpha$ and $S \cong \zeta$, then $\alpha \cong S \cong \zeta$ and $\alpha = \zeta$ by Proposition \ref{prop1.20}(ii), as desired. $\blacksquare$

In light of Theorem \ref{thm1.24}, every well-ordered set $(A, \leq)$ contains some unique ordinal number $\alpha$ in its order type, and so $(A, \leq)$ and $\alpha$ have the same order type. For any ordinal $\alpha$, denote the order type of $\alpha$ as simply $\alpha$, so that the order type of a well-ordered set is a unique ordinal number. The order type of $\mbb{N}$ (i.e., the ordinal number order isomorphic to $\mbb{N}$) under the usual ordering is denoted $\omega$. Take the set $A \cup \mbb{N}$ where $A = \{a_1,...,a_n\}$ ($a_i$ distinct) and $A \cap \mbb{N} = \emptyset$ with the ordering defined by $n < a_j$ for every $n \in \mbb{N}$ and $a_j \in A$, $a_1<...<a_n$ the ordering in $A$, and the usual ordering in $\mbb{N}$. In the discussion following Definition \ref{def2.1}, we showed that $A \cup \mbb{N}$ and $\mbb{N}$ have different order types. The order type of $A \cup \mbb{N}$ (which is clearly well-ordered) is denoted $\omega+n$. If $\omega+n < \omega$, then $\omega+n \subsetneq \omega$, and the map \begin{align}\label{l27}
A \cup \mbb{N} \xrightarrow{\cong} \omega+n \xrightarrow{\subset} \omega \xrightarrow{\cong} \mbb{N}
\end{align}
is an order preserving injection. However, then as above $0 \mapsto 0$ and by induction $k \mapsto k$ for every $k \in \mbb{N}$, but then the map cannot be injective, a contradiction. Thus, $\omega < \omega+n$ for any $n \in \mbb{Z}^{+}$. Moreover, an analogous argument to above shows that \begin{align}\label{l28}
\omega<\omega+1<\omega+2<\omega+3<...
\end{align}
\begin{example}\label{ex10}
The ordinal $1+\omega$ is order isomorphic to the set $\{a\} \cup \{a_1,a_2,...\}$ with $a < a_i$ for every $i \in \mbb{Z}^{+}$ and $a_i \leq a_j$ if and only if $i \leq j$, whereas $\omega+1$ is order isomorphic to the set $\{a_1,a_2,...\} \cup \{a\}$ with the same ordering on $\{a_1,a_2,...\}$ and $a_i < a$ for all $i \in \mbb{Z}^{+}$. The first set is order isomorphic to $\mbb{N}$, and hence $1+\omega = \omega$, but the second is not, so $\omega+1 \neq \omega$. This illustrates that for infinite ordinals, we lose commutativity unlike finite ordinals, and demonstrates a difference between infinity and finiteness with ordinals.
\end{example}
\begin{proposition}\label{prop1.25}
Let $\{\alpha_{i}\}_{i \in \textit{I}}$ be a family of ordinal numbers and $\textit{I}$ a set. Then $\bigcup_{i \in \textit{I}} \alpha_{i}$ is the least upper bound of $\{\alpha_{i}\}_{i \in \textit{I}}$.
\end{proposition}
\noindent $\textbf{Proof}.$ Let $\alpha \in \bigcup_{i \in \textit{I}} \alpha_{i}$. Then there exists $i \in \textit{I}$ such that $\alpha \in \alpha_{i}$ and so $\alpha \subset \alpha_{i} \subset \bigcup_{i \in \textit{I}} \alpha_{i}$. The elements of the union are all ordinals by Lemma \ref{lem1.18} (i), so $\bigcup_{i \in \textit{I}} \alpha_{i}$ is well ordered by Theorem \ref{thm1.19}, and hence an ordinal. Given any $\alpha_{i}$ ($i \in \textit{I}$) by definition $\alpha_{i} \subset \bigcup_{i \in \textit{I}} \alpha_{i}$, and thus by Lemma \ref{lem1.18}(ii), $\alpha_{i} \leq \bigcup_{i \in \textit{I}} \alpha_{i}$, and $\bigcup_{i \in \textit{I}} \alpha_{i}$ is an upper bound of $\{\alpha_{i}\}_{i \in \textit{I}}$. Moreover, if $z < \bigcup_{i \in \textit{I}} \alpha_{i}$, then $z \in  \bigcup_{i \in \textit{I}} \alpha_{i}$ and $z \in \alpha_{j}$ for some $j \in \textit{I}$. So $z < \alpha_{j}$, and $z$ is not an upper bound of  $\{\alpha_{i}\}_{i \in \textit{I}}$. Thus $\bigcup_{i \in \textit{I}} \alpha_{i}$ is the least upper bound of $\{\alpha_{i}\}_{i \in \textit{I}}$, as desired. $\blacksquare$

Using Proposition \ref{prop1.25}, even larger ordinals can be constructed, for example \begin{align}
\omega \cdot 2 & = \sup\{\omega,\omega+1,\omega+2,...\} \\
\omega^{2} & = \sup\{\omega, \omega \cdot 2, \omega \cdot 3,...\} \\
\omega^{\omega} & = \sup\{\omega,\omega^{2},\omega^{3},...\}.
\end{align} 
Note that it is important to distinguish between ordinals and cardinals, for a larger ordinal need not mean a larger cardinal. Repeated application of Lemma \ref{lem1.8} and Theorem \ref{thm1.9} show that each of the ordinals above have a cardinality $\bm{\al}_{0}$, despite their ordinals being larger than $\omega$. In some sense, ordinal numbers act almost like how one might intuitively expect cardinal numbers of infinite sets to act, where adding more elements to the end creates a strictly larger ordinal than the ordinal number of the original well-ordered set. 

So far, in the theory of ordinals, an effort has been made to ensure no usage of the axiom of choice. However, accepting it allows for any non-empty set to have an ordinal associated with it. The $\textit{well-ordering theorem}$, which is equivalent to the axiom of choice, states that for any non-empty set $A$, there exists a linear ordering $\leq$ on $A$, such that $(A, \leq)$ is well-ordered. Thus, any set $A$ has an order $\leq$ such that $(A, \leq)$ is well-ordered, and thus is order isomorphic to a unique ordinal by Theorem \ref{thm1.24}. This allows for a new proof of Theorem \ref{thm1.5}, for let $A = \{A_{i}\}_{i \in \textit{I}}$ be a family of sets with $|A_{i}| = \alpha_{i}$ for each $i \in \textit{I}$ (as in Theorem \ref{thm1.5}). For each $A_{i} \in A$, choose an ordering $\leq_{i}$ such that $(A_{i},\leq_{i})$ is well-ordered (AC). Let $\beta_{i}$ be the ordinal number of the well-ordered sets for each $i \in \textit{I}$. Then $\{\beta_{i}\}_{i \in \textit{I}}$ has a least element by Theorem \ref{thm1.19}, say $\beta_{j}$. Then given any $A_{i} \in A$, the map
 \begin{align}\label{l32}
A_{j} \xrightarrow{\cong} \beta_{j} \xrightarrow{\subset} \beta_{i} \xrightarrow{\cong} A_{i}
\end{align}
is injective. Thus, it follows that $\alpha_{j}$ is the least cardinality of the sets in class $A$. This highlights a deep connection between the theory of cardinals and ordinals. The first proof of Theorem \ref{thm1.5} relied purely on cardinals and the construction an injection, whereas the proof above draws on ordinal numbers and well-orderings. Through AC these two approaches converge, demonstrating how different avenues to infinite sets can ultimately yield the same results, and emphasizing the richness of infinity in set-theoretic mathematics.

The theory of ordinals for this paper is complete. One more result, the principle of transfinite induction, is shown due to its use in the next section, particularly in Theorem \ref{thm2.54}. This statement and proof of the theorem follows the structure of \cite[Theorem 0.7.1]{hungerford}.
\begin{theorem}[Principle of Transfinite Induction]\label{thm1.26}
If $B$ is a subset of a well-ordered set $(A, \leq)$ such that for every $a \in A$, $$\left\{c \in A \;\middle|\; c < a\right\} \subset B \implies a \in B,$$
then $B = A$.
\end{theorem}
\noindent $\textbf{Proof}.$ Suppose, by way of contradiction, that $A-B \neq \emptyset$. Then there exists a least element $a \in A-B$. By definition of least element, there exists no element $b \in A$ with $b<a$ such that $b \notin B$ (otherwise, $b<a \in A-B$). Thus $\{c \in A \mid c<a \} \subset B$. By hypothesis, $a \in B$, so $a \in B \cap (A-B) = \emptyset$, a contradiction. Thus $A-B = \emptyset$ and $A = B$, as desired. $\blacksquare$

While ordinal numbers and order types do not directly answer the counterintuitive nature of cardinal equivalence among infinite sets, they highlight a key structural difference between such sets. To close with the example beginning this subsection, the set $\mbb{Z}$ with its usual ordering isn't even well-ordered, and thus cannot be associated with an ordinal number without changing the order to something non-standard. In contrast, $\mbb{N}$ is already well-ordered by the well-ordering principle, and has an ordinal number $\omega$ associated with it. The order of a set can change its structure drastically, and in part explains why $\mbb{Z}$ and $\mbb{N}$, or even $\mbb{Q}$ (with their usual orderings), do not all have the same properties, yet all have the same cardinality (Proposition \ref{prop1.10}). Ordinals help reinforce the multifaceted nature of infinity. Even in cases when two sets share the same cardinality, their order types and corresponding ordinal numbers (or lack thereof) can offer meaningful distinctions between the sets and their behaviors.
\section{Algebraic Perspectives}\label{section2}

The concept of infinite sets as introduced in the previous section has numerous applications in the discipline of algebra. Although finite algebraic structures already have a rich theory in their own right, infinite sets bring that theory to a whole new level. We now introduce some foundational definitions. A $\textit{n-ary operation}$ ($n \in \mbb{Z}^{+}$) on a nonempty set is a function $f:A^{n} \to A$ for some $n \in \mathbb{Z}^{+}$, where $A^{n} = \prod_{i = 1}^{n} A_{i}$ with $A_{i} = A$ for every $1 \leq i \leq n$. Naturally, a $\textit{binary}$ operation is a 2-ary operation and a $\textit{unary}$ operation a 1-ary operation. A $0$-ary (nullary) operation is simply an element of $A$ (formally a function $f: \{\emptyset\} \to A$). This aligns with viewing constants as 0-ary operations.
\begin{definition}\label{defsec2.1}
An algebraic structure is a nonempty set $A$, together with a family of operations on $A$ (typically binary) and a finite set of identities (e.g., commutativity, associativity, etc.) that those operations must satisfy.
\end{definition}

\subsection{Group Theory} Before delving into theorems about algebraic structures that directly use infinite sets and structures, we naturally must first develop the fundamentals about the most basic of such structures, primarily groups, as we now define.
\begin{definition}\label{defsec2.2}
A semigroup is a nonempty set $A$ together with a binary operation (denoted as addition here) such that for every $a,b,c \in A$ $$(a+b)+c = a+(b+c).$$
The property above is called associativity. A monoid is a semigroup with an element $0 \in A$ such that $a+0 = 0+a = a$ for every $a \in A$ (called an additive identity). A group is a monoid with the property that for every $a \in A$, there exists $b \in A$ such that $a+b = b+a = 0$ (called an additive inverse). A group $G$ is called abelian if $a+b = b+a$ for every $a,b \in G$. 
\end{definition}
\begin{example}\label{ex11}
The set $\mbb{Z}^{+}$ under usual addition is clearly associative and hence a semigroup, but has no identity element and thus is not a monoid. The set $\mbb{N}$ does have an identity of $0$, and thus is a monoid under usual addition. The entire set $\mbb{Z}$ of integers then adjoins additive inverses and thus creates a group under usual addition.
\end{example}
\begin{example}\label{ex12}
The finite set of all bijections $f:\{1,2,3\} \to \{1,2,3\}$ forms a group under function composition, denoted $S_3$, and called the \textit{symmetric group} on $\{1,2,3\}$. However, $S_3$ is not abelian, as one may verify maps defined by $f(1) = 2$, $f(2) = 1$, $f(3) = 3$ and $g(1) = 1$ $g(2) = 3$, $g(3) = 2$, do not commute ($(fg)(1) \neq (gf)(1)$). In contrast, the set $\mbb{Z}_{4} = \{\bar{0},\bar{1},\bar{2},\bar{3}\}$ of integers modulo $4$ is a finite, additive abelian group under usual addition ($\bar{a}+\bar{b} = \overline{a+b}$).
\end{example}
The focus of this paper will be on groups rather than semigroups or monoids. Notice that for a group $G$, if $0_1, 0_2 \in G$ are additive identities, and $b,c \in G$ are additive inverses of $a \in G$, then \begin{align}
0_1 & = 0_1+0_2 = 0_2+0_1 = 0_2\label{l33} \\
b & = b+0 = b+(a+c) = (b+a)+c = 0+c = c. \label{l34}
\end{align}
Thus the additive identity of a group is unique, as well as the additive inverse of an element. In light of (\ref{l34}), we denote the additive inverse of an element $a$ as $-a$. In general, from now on, multiplicative notation will be used when dealing with general groups or non-abelian groups, and additive notation with abelian groups. The multiplicative identity for a monoid $G$ will be denoted either $1_{G}$ or $1$ depending on the context.

We now build the theory necessarily to define the quotient of a group to a normal subgroup (Definition \ref{defsec2.4}). Quotient groups aid in proving numerous crucial results about infinity in group theory (e.g., Theorem \ref{thm2.11}, Theorem \ref{thm2.14}). Let $G$ be a set and $R$ ($\sim$) an equivalence relation. The set of all equivalence classes of $G$ with respect to $R$ is denoted $G/R$. Before proving the next result, we note that if $R$ ($\sim$) is an equivalence relation on a set $A$, \begin{align}\label{l35}
a \sim b \iff \bar{a} = \bar{b}.
\end{align} 
To prove this, notice if $a \sim b$, then if $c \in \bar{a}$, we have $c \sim a \sim b$ and $c \in \bar{b}$, and by symmetry $\bar{a} = \bar{b}$. Conversely, if $\bar{a} = \bar{b}$, then $a \in \bar{b}$, so $a \sim b$. For a further proof of the following theorem, see \cite[Theorem I.1.5]{hungerford}.
\begin{theorem}\label{thm2.1}
Let $R$ ($\sim$) be an equivalence relation on a monoid $G$ such that $a_1 \sim a_2$ and $b_1 \sim b_2$ implies that $a_1b_1 \sim a_2b_2$ for all $a_1,a_2,b_1,b_2 \in G$. Then the set $G/R$ is a monoid under the binary operation defined by $(\bar{a})(\bar{b}) = \overline{ab}$. If $G$ is a group, so is $G/R$, and if $G$ is abelian, then so is $G/R$.
\end{theorem}
\noindent $\textbf{Proof}.$ If $\bar{a_1} = \bar{a_2}$ and $\bar{b_1} = \bar{b_2}$ ($a_1,a_2,b_1,b_2 \in G$), then  $a_1 \sim a_2$, and $b_1 \sim b_2$ by (\ref{l35}). Thus, by hypothesis, $a_1b_1 \sim a_2b_2$, and by (\ref{l35}) again $\overline{a_1b_1} = \overline{a_2b_2}$, so the operation is well-defined. We see for all $a,b,c \in G$, $\bar{a}(\bar{b}\bar{c}) = \bar{a}(\overline{bc}) = \overline{a(bc)} = \overline{(ab)c} = \overline{ab}\bar{c} = (\bar{a}\bar{b})\bar{c},$
hence $G/R$ is associative. We see for all $a \in G$, $\bar{a}\bar{1} = \overline{a1} = \bar{a} = \overline{1a} = \bar{1}\bar{a}$, whence $\bar{1} \in G/R$ is the identity, and $G/R$ is a monoid. If $G$ is a group, then $\bar{a} \in G/R$ naturally has multiplicative inverse $\overline{a^{-1}}$ as $\bar{a}\overline{a^{-1}} = \overline{aa^{-1}} = \bar{1} = \overline{a^{-1}a} = \overline{a^{-1}}\bar{a}$, so that $G/R$ is a group. Finally, if $G$ is abelian and $\bar{a},\bar{b} \in G/R$, then $\bar{a}\bar{b} = \overline{ab} = \overline{ba} = \bar{b}\bar{a},$
and $G/R$ is abelian, as desired. $\blacksquare$ 

We now extend the notation of congruence modulo $m$ in $\mbb{Z}$. Recall for $a,b,m \in \mbb{Z}$ and $m > 0$, $a \equiv b \hspace{1mm} (mod \hspace{1mm} m)$ if $m \vert a-b$, that is $a-b \in m\mbb{Z}$. More generally, we have 
\begin{definition}\label{defsec2.3}
Let $H$ be a subgroup of a group $G$ and $a,b \in G$. a is right congruent to b modulo H, denoted $a \equiv_{r} b \hspace{1mm} (mod \hspace{1mm} H)$ if $ab^{-1} \in H$. a is left congruent to b modulo $H$, denoted $a \equiv_{l} b \hspace{1mm} (mod \hspace{1mm} H),$ if $a^{-1}b \in H$.
\end{definition}

\noindent We note that if $G$ is abelian, then \begin{align}\label{l36}
a \equiv_{r} b \hspace{1mm} (mod \hspace{1mm} H) \iff ab^{-1} \in H \iff (ab^{-1})^{-1} = a^{-1}b \in H \iff a \equiv_{l} b \hspace{1mm} (mod \hspace{1mm} H)
\end{align}
and left and right congruence are equivalent.  A subset $H$ of a group $G$ that is also a group under the same operation is called a subgroup of $G$, denoted $H \leq G$. Notice that the identity $1_{H}$ of a subgroup $H$ of a group $G$ equals the identity $1_{G} \in G$, for $1_{H} = 1_{H}^{-1}(1_{H}1_{H}) = 1_{H}^{-1}(1_{H}) = 1_{H}^{-1}(1_{H}1_{G}) = 1_{G}$.
It then follows that the inverse of $a \in H$ is the same as the inverse of $a \in G$. For further proofs of the following theorem, see \cite[Theorem I.4.2]{hungerford}. Notice that (iii) highlights a relationship between the cardinalities of $aH$, $Ha$, and $H$ even in the case when $H$ is an infinite subgroup of $G$.
\begin{theorem}\label{thm2.2}
Let $H$ be a subgroup of a group $G$ \\
(i) Both right and left congruence modulo $H$ are equivalence relations on $G$. \\
(ii) The equivalence class of $a \in G$ under right [resp. left] congruence modulo $H$ is the set $Ha = \{ha \mid h \in H\}$ [resp. $aH = \{ah \mid h \in H\}$]. \\
(iii) $|Ha| = |H| = |aH|$ for all $a \in G$
\end{theorem} 
\noindent $\textbf{Proof}.$ We prove the theorem for right congruence, left congruence then follows by symmetry.

[(i)]. Clearly $a \equiv_{r} a \hspace{1mm} (mod \hspace{1mm} H)$ as $aa^{-1} = 1 \in H$. If $a \equiv_{r} b \hspace{1mm} (mod \hspace{1mm} H)$, then $ab^{-1} \in H$, so $ba^{-1} = (ab^{-1})^{-1} \in H$ and $b \equiv_{r} a \hspace{1mm} (mod \hspace{1mm} H)$. Finally, if $a \equiv_{r} b \hspace{1mm} (mod \hspace{1mm} H)$ and $b \equiv_{r} c \hspace{1mm} (mod \hspace{1mm} H)$, then $ab^{-1} \in H$ and $bc^{-1} \in H$, so by closure $ac^{-1} = (ab^{-1})(bc^{-1}) \in H$ and $a \equiv_{r} c \hspace{1mm} (mod \hspace{1mm} H)$.

[(ii)] The equivalence class of $a \in G$ is given by \begin{align}
\bar{a} & = \left\{x \in G\;\middle|\;x \equiv_{r} a \hspace{1mm} (mod \hspace{1mm} H)\right\} = \left\{x \in G\;\middle|\; xa^{-1} = h \in H\right\}\label{l37} \\ & = \left\{x \in G\;\middle|\; x = ha, h \in H\right\}\label{l38} \\ & = Ha\label{l39}.
\end{align}

[(iii)] Consider the map $f:H \to Ha$ by $h \mapsto ha$. Clearly $f$ is surjective, and if $h_1a = h_2a$, then $h_1 = (h_1a)a^{-1} = (h_2a)a^{-1} = h_2$, so $f$ is injective and thus bijective, as desired. $\blacksquare$

The sets $Ha$ and $aH$ in Theorem \ref{thm2.2} are called the right and left cosets of $a \in G$, respectively. For a certain class of groups, called normal subgroups, which need not be abelian, left and right congruence are also equivalent. In fact, this statement could be used as a definition of a normal subgroup. For further proofs and discussion of the following Theorem, see \cite[Lemma 2.65]{rotman}, \cite[Theorem I.5.1]{hungerford}, and \cite[Proposition 3.6]{dummit}.
\begin{theorem}\label{thm2.3}
Let $N$ be a subgroup of a group $G$. Then the following conditions are equivalent.

(i) Left and right congruence modulo $N$ are equivalent;

(ii) every left coset of $N$ in $G$ is a right coset of $N$ in $G$;

(iii) $aN = Na$ for all $a \in G$;

(iv) for all $a \in G$, $aNa^{-1} \subset N$, where $aNa^{-1} = \{ana^{-1} \mid n \in N\}$;

(v) for all $a \in G$, $aNa^{-1} = N$.
\end{theorem}
\noindent $\textbf{Proof}.$ [(i) $\Leftrightarrow$ (iii)] Two equivalence relations $\sim_{1}$ and $\sim_{2}$ are equivalent if and only if the equivalence class of each element under $\sim_{1}$ is equal to its equivalence class under $\sim_{2}$, for letting $\bar{a}_{1}$ be the equivalence class of $a \in G$ under $\sim_{1}$, and $\bar{a}_{2}$ under $\sim_{2}$, we have that $\sim_{1}$ and $\sim_{2}$ are equivalent if and only if \begin{align}\label{l40}
b \in \bar{a}_{1} \iff b \sim_{1} a \iff b \sim_{2} a \iff b \in  \bar{a}_{2}. \quad \forall a \in G.
\end{align} 
By Theorem \ref{thm2.2}, the equivalence classes under right and left congruence are the cosets $Na$ and $aN$, respectively.

[(ii) $\Rightarrow$ (iii)] Fix $a \in G$. If $aN = Nb$ for some $b \in G$, then $a \in Na \cap Nb$, so $a \in Nb$ and by Theorem \ref{thm2.2}, $a \equiv_{r} b \hspace{1mm} (mod \hspace{1mm} H)$. Thus, by (\ref{l35}) and Theorem \ref{thm2.2}, $Na = Nb$ and so $aN = Na$. 

[(iii) $\Rightarrow$ (iv)] If $aN = Na$ ($a \in G$), then for all $n_1 \in N$, there exists $n_2 \in N$ such that $an_1 = n_2a$, so $an_1a^{-1} = n_2$ and $aNa^{-1} \subset N$. [(iv) $\Rightarrow$ (v)] By hypothesis $aNa^{-1} \subset N$. Moreover, since $a^{-1} \in G$, $a^{-1}Na \subset N$. Notice for every $n \in N$, $n = a(a^{-1}na)a^{-1} \in aNa^{-1}$ and thus $N \subset aNa^{-1}$. [(v) $\Rightarrow$ (ii)] Fix $an_1 \in aN$ ($a \in G$). Since $aNa^{-1} = N$, there exists $n_2 \in N$ with $an_1a^{-1} = n_2$, so $an_1 = n_2a \in Na$ and $aN \subset Na$. Similarly, if $n_3a \in Na$, since $a^{-1}Na = N$, there exists $n_4 \in N$ such that $a^{-1}n_3a = n_4$, so $n_3a = an_4 \in aN$ and $Na \subset aN$. Thus $aN = Na$, and the left coset $aN$ is also the right coset $Na$, as desired. $\blacksquare$

A subgroup $N$ of a group $G$ satisfying any of the equivalent conditions of Theorem \ref{thm2.3} is called a $\textit{normal}$ subgroup of $G$, and denoted $N \lhd G$. In light of Theorem \ref{thm2.3}, when dealing with  a normal subgroup $N$ of a group $G$, we'll omit the ``$r$" and ``$l$" subscripts when working modulo $N$. 
\begin{definition}\label{defsec2.4}
Let $G$ be a group and $N \lhd G$. The set $G/N$ is defined as the quotient $G/R$, where $R$ is the equivalence relation of left (or equivalently right) congruence modulo $N$. By Theorem \ref{thm2.2}, $G/N$ is the set $$G/N = \left\{aN\;\middle|\;a \in G\right\}.$$ 
\end{definition}
\begin{theorem}\label{thm2.4}
Let $G$ be a group and $N \lhd G$. The quotient $G/N$ is a group under multiplication, defined by $(aN)(bN) = (ab)N$.
\end{theorem}
\noindent $\textbf{Proof}.$ By Theorem \ref{thm2.1}, it suffices to show that $a_1 \equiv a_2 \hspace{1mm} (mod \hspace{1mm} N)$ and $b_1 \equiv b_2 \hspace{1mm} (mod \hspace{1mm} N)$ implies that $a_1b_1 \equiv a_2b_2 \hspace{1mm} (mod \hspace{1mm} N)$ for all $a_1,a_2,b_1,b_2 \in G$. By hypothesis $a_1a_2^{-1} = n_1 \in N$ and $b_1b_2^{-1} = n_2 \in N$. Thus $(a_1b_1)(a_2b_2)^{-1} = a_1b_1b_2^{-1}a_2^{-1} = (a_1n_2)a_2^{-1}$. Since $N$ is normal, $a_1N = Na_1$, and thus $a_1n_2 = n_3a_1$ for some $n_3 \in N$. Thus, $(a_1b_1)(a_2b_2)^{-1} = (a_1n_2)a_2^{-1} = (n_3a_1)a_2^{-1} = n_3n_1 \in N$ and $a_1b_1 \equiv a_2b_2 \hspace{1mm} (mod \hspace{1mm} N)$, as desired. $\blacksquare$

For further proofs of the above theorem, see \cite[Proposition 3.5]{dummit}, \cite[Theorem 2.67]{rotman}, and \cite[Theorem I.5.4]{hungerford}. In light of Theorem \ref{thm2.4}, we call $G/N$ the $\textit{quotient group}$ of $G$ and $N$. Now that we have defined quotients, naturally we define products of groups. The \textit{product group} of two groups $G$ and $H$ is the Cartesian product $G \times H$ (see (\ref{l1})) along with multiplication defined by \begin{align}\label{l41}
(a,a')(b,b') = (ab,a'b') \hspace{2mm} a,b \in G, a',b' \in H.
\end{align}
Naturally the identity is $1 = (1_{G},1_{H})$ and the multiplicative inverse of $(a,b) \in G \times H$ is $(a^{-1},b^{-1}) \in G \times H$. It is then easy to verify $G \times H$ forms a group under (\ref{l41}). Extending this to a (nonempty) family of groups $\{G_{i} \mid i \in \textit{I}\}$, the \textit{direct product} is defined as the set $\prod_{i \in \textit{I}} G_{i}$ together with multiplication $(fg)(i) = f(i)g(i)$ for $f,g \in \prod_{i \in \textit{I}} G_{i}$ and $i \in \textit{I}$. The identity is then naturally the map $i \mapsto 1_{G_{i}}$ ($i \in \textit{I}$) and the multiplicative inverse of $f \in \prod_{i \in \textit{I}} G_{i}$ the map $i \mapsto (f(i))^{-1}$ ($i \in \textit{I}$). The $\textit{weak direct product}$ of $\{G_{i} \mid i \in \textit{I}\}$ is the subset $\prod^{w}_{i \in \textit{I}} G_{i}$ of $\prod_{i \in \textit{I}} G_{i}$ consisting only of elements with finitely many non-identity outputs  (see \cite{hungerford}). Is it easy to verify $\prod^{w}_{i \in \textit{I}} G_{i}$ forms a subgroup of $\prod_{i \in \textit{I}} G_{i}$. When dealing with abelian groups, the weak direct product is denoted $\bigoplus_{i \in \textit{I}} G_{i}$ and called the $\textit{direct sum}$ of $\{G_{i} \mid i \in \textit{I}\}$ (see \cite{lang}).

One way we attempt to study an algebraic structure is to study maps between two such algebraic structures of the same type. While some value could be brought by studying arbitrary functions, it is particularly useful if the function preserves the structures between the sets. More precisely, we have the following definition
\begin{definition}\label{defsec2.5}
Let $A$ and $B$ be semigroups. A map $f:A \to B$ is a semigroup homomorphism if $$f(ab) = f(a)f(b).$$
If $A$ and $B$ are groups, $f$ is called a group homomorphism, and if $A$ and $B$ are monoids and $f(1_{A}) = 1_{B}$, $f$ is called a monoid homomorphism. If $f$ is injective [resp. surjective], $f$ is called a monomorphism [resp. epimorphism].
\end{definition}
\begin{example}\label{ex13}
Given two monoids $G$ and $H$, the map $f:G \to H$ by $f(g) = 1_{H}$ for all $g \in G$ is called the \textit{trivial homomorphism}. If $G$ and $H$ are groups, then $f$ is a group homomorphism. In a similar fashion, given a semigroup, monid, or group $G$, the identity map $1_{G}:G \to G$ by $g \mapsto g$ is a homomorphism.
\end{example}
\begin{example}\label{ex14}
The map $\mbb{Z} \to \mbb{Z}_{n}$ by $k \mapsto \bar{k}$, where $\mbb{Z}_{n} = \{\bar{0},\bar{1},...,\overline{n-1}\}$ is the integers modulo $n$ is a group epimorphism. 
\end{example}
\begin{example}\label{ex15}
The map $\det:GL_{n}(\mbb{R}) \to \mbb{R}$ by $A \mapsto \det(A)$, where $GL_{n}(\mbb{R})$ is the set of all invertible $n$-by-$n$ matrices with elements in $\mbb{R}$, and $\det A$ the determinant of matrix $A$, is a homomorphism of groups under multiplication as $\det(AB) = (\det A)(\det B)$ \cite[Theorem 10.40]{axlerlinear}.
\end{example}
Note that in the case of groups, $f(1) = f(1(1)) = f(1)f(1)$ so $f(1) = 1$, and $f(aa^{-1}) = f(a)f(a^{-1}) = 1$ so $f(a^{-1}) = (f(a))^{-1}$. As expected, it naturally follows from the definition that the composition of two homomorphisms is a homomorphism. Let $f:G \to H$ be a homomorphism of groups. The $\textit{kernel}$ of $f$ is the set of all $a \in G$ such that $f(a) = 1_{H}$, denoted $\ker f$. Since given $a,b \in \ker f$, $f(ab) = f(a)f(b) = 1$, it follows naturally that $\ker f \leq G$. Moreover, \begin{align}\label{l42}
f \hspace{2mm} \text{is injective} \iff \ker f = \{1\}.
\end{align} 
To prove this, observe if $f$ is injective, then if $a \in \ker f$, $f(a) = 1 = f(1)$ so $a = 1$ and $\ker f = \{1\}$. Conversely, if $\ker f = \{1\}$ and $f(a) = f(b)$ ($a,b \in G$), then $f(ab^{-1}) = 1$ and $a = b$. For further proofs and discussion of the following theorem, see  \cite[Proposition 3.7]{dummit} and \cite[Theorem I.5.5]{hungerford}. Notice here that, even in the case where $N$ is an infinite subgroup of $G$, it is still the kernel of home homomorphism. This highlights a continuity between normal subgroups, kernels, and infinity.
\begin{proposition}\label{prop2.5}
A subgroup $N$ of a group $G$ is normal if and only if it is the kernel of some homomorphism.
\end{proposition}
\noindent $\textbf{Proof}.$ If $N$ is normal, consider the map $\pi:G \to G/N$ defined by $\pi(a) = aN$. Since $\pi(ab) = (ab)N = (aN)(bN) = \pi(a)\pi(b)$, $\pi$ is a homomorphism. We see $\ker \pi = \{a \in G \mid \pi(a) = N\} = \{a \in G \mid aN = N\} = \{a \in G \mid a \in N\} = N$. Conversely, consider a homomorphism $f:G \to H$ of groups with $\ker f = N$. If $n \in N$ and $a \in G$, then $f(ana^{-1}) = f(a)f(n)f(a^{-1}) = f(a)f(a)^{-1} = 1$. Thus $aNa^{-1} \subset N$ and $N$ is normal by Theorem \ref{thm2.3}, as desired. $\blacksquare$

The map $\pi$ in the proof of Proposition \ref{prop2.5} is called the $\textit{canonical epimorphism}$. We now prove some fundamental results about homomorphisms. For further proofs of Theorem \ref{thm2.6}(i), see \cite[Theorem 14.11]{fraleigh}, \cite[Theorem 3.16]{dummit}; for (ii), see \cite[Corollary I.5.8]{hungerford}. For Theorem \ref{thm2.7}(i,ii), see \cite[Theorems I.8.10-11]{hungerford}.
\begin{theorem}\label{thm2.6}
Let $f:G \to H$ be a homomorphism of groups. \\
(i) (First Isomorphism Theorem) The map $\varphi:G/ \ker f \to \Im f$ by $a(\ker f) \mapsto f(a)$ is an isomorphism. \\
(ii) If $f$ is an isomorphism, $N \lhd G$, $M \lhd H$, and $f(N) = M$, then the map $G/N \to H/M$ by $aN \mapsto f(a)M$ is an isomorphism.
\end{theorem}
\noindent $\textbf{Proof}.$ [(i)] Let $N = \ker f$. If $aN = bN$, then $b = an$ for some $n \in N$, and $f(b) = f(an) = f(a)f(n) = f(a)$, so $\varphi(aN) = f(a) = f(b) = \varphi(bN)$ and $\varphi$ is well-defined. Since $\varphi((aN)(bN)) = \varphi(abN) = f(ab) = f(a)f(b) = \varphi(aN)\varphi(bN)$, $\varphi$ is a homomorphism. Clearly $\Im \varphi = \Im f$. Since $aN \in \ker \varphi$ if and only if $f(a) = 1$,  $\ker \varphi = \{aN \mid a \in N\} = N = 1_{G/N}$ and $\varphi$ is injective by (\ref{l42}). Thus $\varphi$ is an isomorphism. [(ii)] Consider the composition $g:G \xrightarrow{f} H \to \xrightarrow{\pi} H/M$. Then $\ker g = \{a \in G \mid f(a) \in M\} = f^{-1}(M)$. Since $f(N) = M$ and $f$ bijective, $N = f^{-1}(M)$. Thus by (i) the map $aN \mapsto g(a) = f(a)M$ is an isomorphism, as desired. $\blacksquare$
\begin{theorem}\label{thm2.7}
Let $\{f_{i}:G_{i} \to H_{i} \mid i \in \textit{I}\}$ be a family of homomorphisms of groups and let $\{N_{i} \mid i \in \textit{I}\}$ be a family of groups with $N_{i} \lhd G_{i}$ for each $i \in \textit{I}$. \\
(i) Let $f = \prod^{w} f_{i}$ be the map $\prod^{w}_{i \in \textit{I}} G_{i} \to \prod^{w}_{i \in \textit{i}} H_{i}$ by $\{a_{i}\}_{i \in \textit{I}} \mapsto \{f_{i}(a_{i})\}_{i \in \textit{I}}$. Then $f$ is a homomorphism of groups, $\ker f = \prod^{w}_{i \in \textit{I}} \ker f_{i}$, and $\Im f = \prod^{w}_{i \in \textit{I}} \Im f_{i}$. Consequently $f$ is a monomorphism [resp. epimorphism] if and only if each $f_{i}$ is. \\
(ii) $\prod^{w}_{i \in \textit{I}} N_{i}$ is a normal subgroup of $\prod^{w}_{i \in \textit{I}} G_{i}$ and $$\prod^{w}_{i \in \textit{I}} G_{i}/ \prod^{w}_{i \in \textit{I}} N_{i} \cong  \prod^{w}_{i \in \textit{I}} G_{i}/N_{i}.$$
\end{theorem}
\noindent $\textbf{Proof}.$ [(i)] For ease of notation, we drop the ``$i \in \textit{I}$" from $\{a_{i}\}_{i \in \textit{I}}$. Since identities map to identities with homomorphisms, it follows that $\Im f \subset \prod^{w}_{i \in \textit{i}} H_{i}$. Let $\{a_{i}\},\{b_{i}\} \in \prod^{w}_{i \in \textit{I}} G_{i}$. Then $f(\{a_{i}\}\{b_{i}\}) = f(\{a_{i}b_{i}\}) = \{f_{i}(a_i b_i)\} = \{f_i(a_i)\}\{f_i(b_i)\} = f(\{a_i\})f(\{b_i\})$ and $f$ is a homomorphism. Observe $\{x_{i}\} \in \ker f$ if and only if $f_{i}(x_{i}) = 1_{H_i}$ for each $i \in \textit{I}$ and thus $\ker f =  \prod^{w}_{i \in \textit{I}} \ker f_{i}$. Similarly, if $f(\{x_i\}) = \{f_i(x_i)\} \in \Im f$, then for each $i \in \textit{I}$, $f_i(x_i) \in \Im f_{i}$, and if $\{f_i(y_i)\} \in  \prod^{w}_{i \in \textit{I}} \Im f_{i}$, then $\{z_{i}\} \in \prod^{w}_{i \in \textit{I}} G_{i}$ with $z_{i} = y_{i}$ if $f_i(y_i) \neq 1$ and $z_{i} = 1$ otherwise has $\{f_i(y_i)\} = f(\{z_i\}) \in \Im f$ so $\Im f = \prod^{w}_{i \in \textit{I}} \Im f_{i}$. Thus if $f$ is a monomorphism [resp. epimorphism], we must have $\ker f_i = \{1_{H_i}\}$ [resp. $\Im f_i = H_i$], or otherwise $\ker f \neq \{1\}$ [resp. $\Im f \subsetneq  \prod^{w}_{i \in \textit{I}} H_{i}$]. Conversely, if each $f_{i}$ ($i \in \textit{I}$) is a monomorphism [resp. epimorphism], we have $\ker f = \prod^{w}_{i \in \textit{I}} \ker f_{i} = \prod^{w}_{i \in \textit{I}} \{1_{H_{i}}\} = \{1\}$ [resp. $\Im f = \prod^{w}_{i \in \textit{I}} \Im f_{i} =  \prod^{w}_{i \in \textit{I}} H_{i}$] and $f$ is a monomorphism [resp. epimorphism].

[(ii)] Consider the family of maps $\{\pi_{i}:G_{i} \to G_{i}/N_{i} \mid i \in \textit{I}\}$ with $\pi_{i}$ canonical epimorphism. By [(i)] the map $\prod \pi_{i}: \prod^{w}_{i \in \textit{I}} G_{i} \to \prod^{w}_{i \in \textit{I}} G_{i}/N_{i}$ is an epimorphism with kernel $\prod^{w}_{i \in \textit{I}} N_{i}$. Therefore by Theorem \ref{thm2.6}, $\prod^{w}_{i \in \textit{I}} G_{i}/\prod^{w}_{i \in \textit{I}} N_{i} \cong \prod^{w}_{i \in \textit{I}} G_{i}/N_{i}$, as desired. $\blacksquare$

We now define the generating set of a subset of a group $G$. Naturally for $a \in G$ and $n \in \mbb{N}$, $a^{0} = 1$, $a^{n} = a \cdots a$ (n times) if $n \in \mbb{Z}^{+}$, and $a^{-n} = \qty(a^{-1})^{n}$. It is then easy to verify the usual exponent properties $a^{n+m} = a^n(a^m)$ and $(a^n)^m = a^{nm}$ for $a \in G, m,n \in \mbb{Z}$. Given a subset $X \subset G$, the set $\{a_1^{n_1} \cdots a_k^{n_k} \mid a_{i} \in X, n_{i} \in \mbb{Z}, 1 \leq i \leq n\}$ naturally forms a subgroup of $G$ with inverses $\qty(a_1^{n_1} \cdots a_k^{n_k})^{-1} = a_k^{-n_k} \cdots a_1^{-n_1}$.
\begin{definition}\label{defsec2.6}
Let $G$ be a group and $X \subset G$. Then the subgroup $\langle X \rangle$ of $G$ consisting of all finite products $a_1^{n_1} \cdots a_k^{n_k}$ with $a_i \in X, n_{i} \in \mbb{Z}$ ($1 \leq i \leq n$) is called the subgroup of $G$ generated by $X$, with $\langle \emptyset \rangle \coloneqq \{1\}$. If $G = \langle X \rangle$ and $X$ is finite, $G$ is called finitely generated. Each subgroup $\langle x \rangle$ of $G$ with $x \in G$ is called the cyclic subgroup generated by $x$. If $G = \langle a \rangle$ for some $a \in G$, $G$ is called cyclic.
\end{definition}
\begin{example}\label{ex16}
The group $(\mbb{Z}, +)$ is cyclic with generator $1 \in \mbb{Z}$. Clearly, any finite group is also finitely generated. The group $(\mbb{Q}, +)$ is not finitely generated. To prove this, suppose $\mbb{Q} = \langle \frac{a_1}{b_n},...,\frac{a_n}{b_n} \rangle$. Then it follows that every $	q \in \mbb{Q}$ can be written as $q = \frac{z}{b_1 \cdots b_n}$, with $z \in \mbb{Z}$. But $\frac{1}{2 b_1 \cdots b_n} \in \mbb{Q}$ cannot be written as such multiple, a contradiction. Thus, $\mbb{Q}$ is not finitely generated. This shows that some infinite sets like $\mbb{Z}$ may be finitely generated, but other, more complicated structures, like $\mbb{Q}$, are not, highlighting the complexity and richness of infinity in group theory.
\end{example}
For a finite group $G$ and $a \in G$, there must exist distinct $n, m \in \mbb{N}$ such that $a^{n} = a^{m}$, for otherwise $|\{a^{n}|n \in \mbb{N}\}| = \bm{\al}_{0} > |G|$. Thus there exists $n \in \mbb{Z}^{+}$ such that $a^{n} = aa^{n-1} = 1$ and $a^{-1} = a^{n-1}$. It then follows that for finite groups, only natural powers in generating sets need be considered. In comparison, if $G$ is infinite, each integer power of $a \in G$ may very well be different, and thus all integer powers of $a$ must be considered, highlighting how infinite groups exhibit fundamentially different generative behavior from finite ones. 

A $\textit{basis}$ of an abelian group $F$ is a subset $X$ of $F$ such that $F = \langle X \rangle$ and for distinct $x_{1},...,x_{k} \in X$ and $n_{i} \in \mbb{Z}$ ($1 \leq i \leq k$) \begin{align}\label{l43}
n_1x_1+...+n_kx_k = 0 \implies n_i = 0 \hspace{2mm} \forall i \in \{1,...,k\}.
\end{align}
The sum $\sum_{i = 1}^{n} n_{i}x_{k}$ is called a $\textit{linear combination}$ of $x_1,...,x_k$. Readers familiar with linear algebra should be careful to assume theorems about bases of vector spaces hold for bases of abelian groups (e.g., see \cite[Exercise II.1.2]{hungerford}). Given $\{G_i \mid i \in \textit{I}\}$ a family of groups, the family of maps $\{\iota_k:G_{k} \to \prod^{w}_{i \in \textit{I}}\}$ given by $\iota_k(a) = \{a_i\}_{i \in \textit{I}}$, where $a_i = 1$ if $i \neq k$ and $a_k = a$ for each $k \in \textit{I}$ are called the $\textit{canonical injections}$. It is easy to verify each $\iota_{k}$ is a monomorphism of groups. We now prove an important proposition. For a further proof of Proposition \ref{prop2.8}(i), see \cite[Theorem I.3.2]{hungerford}; for (ii), see is \cite[Theorem I.8.6]{hungerford}.
\begin{proposition}\label{prop2.8}
Let $G$ be a group and $\{N_{i} \mid i \in \textit{I}\}$ a family of normal subgroups. \\
(i) If $G = \langle a \rangle$ with $a \in G$, then if $|G|$ is infinite, $G$ is isomorphic to the additive group $\mbb{Z}$, and if $|G| = m$ ($m \in \mbb{Z}^{+}$), $G$ is isomorphic to the additive group $\mbb{Z}/m\mbb{Z}$. \\
(ii) If $G = \langle \bigcup_{i \in \textit{I}} N_{i} \rangle$ and for each $k \in \textit{I}$, $N_{k} \cap \langle \bigcup_{i \neq k} N_{i} \rangle = \langle 1 \rangle$, then $G \cong \prod^{w}_{i \in \textit{I}} N_{i}$.
\end{proposition}
\noindent $\textbf{Proof}.$ [(i)] Consider the map $\varphi: \mbb{Z} \to G$ by $n \mapsto a^{n}$. Clearly $\varphi$ is well-defined and surjective. Notice $\varphi(n+l) = a^{n+l} = a^{n}a^{l} = \varphi(n) \varphi(l)$ and thus $\varphi$ is an epimorphism. Suppose $G$ is infinite and $a^k = 1$. If $k \neq 0$, then by the division algorithm \cite[Theorem 1.2]{niven}, for all $a^{c} \in G$ there exists $x,r \in \mbb{Z}$ with $0 \leq r < k$ and $a^{c} = a^{kx+r} = a^{r}$ and $G = \{1,...,a^{k}\}$, a contradiction. Thus $\ker \varphi = \{0\}$ and $\mbb{Z} \cong G$ by (\ref{l42}). If $|G| = m$, let $A = \{z \in \mbb{Z}^{+} \mid a^{z} = 1 \}$. By above $A \neq \emptyset$, so let $y = \min A$. Similarly, by the division algorithm again, if $y<m$ or $y>m$, it follows that $|G| < m$ or $|G| > m$, respectively. Thus $y = m$ and $a^m = 1$. Then immediately  $m \mbb{Z} \subset \ker \varphi$, and given $\alpha \in \ker \varphi$, by the minimality of $m$ and the division algorithm, $m| \alpha$ and $\alpha \in m \mbb{Z}$. Thus $\ker \varphi = m \mbb{Z}$ and $\mbb{Z}/m \mbb{Z} \cong G$ by Theorem \ref{thm2.6}(i).

 [(ii)] Let $\{a_{i}\} \in \prod^{w}_{i \in \textit{I}} N_{i}$ and $I_{0}$ be the finite set $\{i \in \textit{I}| a_{i} \neq 1\}$. Let $a \in N_{i}$ and $b \in N_{j}$ with $i \neq j$ ($i,j \in \textit{I}$). Since $N_{i}, N_{j} \lhd G$, $bab^{-1} \in N_{i}$ and $ab^{-1}a^{-1} \in N_{j}$ so $(bab^{-1})a^{-1} = b(ab^{-1}a^{-1}) \in N_{i} \cap N_{j} = \langle 1 \rangle$, which implies $ba = ab$. Thus it follows that $\prod_{i \in \textit{I}_{0}} a_{i}$ is a well-defined element of $G$. Consider the map $\varphi: \prod^{w}_{i \in \textit{I}} N_{i} \to G_{i}$ given by $\{1_{N_{I}}\} \mapsto 1$ and $\{a_{i}\} \mapsto \prod_{i \in \textit{I}_{0}} a_{i}$ otherwise. Since elements in distinct normal subgroups from the family commute, it follows that $\varphi$ is a homomorphism. By definition of $\varphi$ and $\iota_{i}$, we see $\varphi \iota_{i}(a_i) = a_i$ for $a_i \in N_{i}$.  Since $G$ is generated by the subgroups $N_{i}$ and elements from distinct $N_{i}$ and $N_{j}$ commute, given $a \in G$ we can write $a$ as a product $\prod_{i \in \textit{J}_{0}} a_{i}$, where $a_{i} \in N_{i}$ and $\textit{J}_{0}$ a finite subset of $\textit{I}$. Thus we have $\prod_{i \in \textit{J}_{0}} \iota_{i}(a_{i}) \in \prod^{w}_{i \in \textit{I}} N_{i}$ and $\varphi(\prod_{i \in \textit{J}_{0}} \iota_i(a_i)) = \prod_{i \in J_{0}} \varphi \iota_i(a_i) = \prod_{i \in J_{0}} a_{i} = a$, and $\varphi$ is an epimorphism.  Suppose $\varphi(\{c_{i}\}) = \prod_{i \in K_{0}} c_{i} = 1$. Let $K_{0} = \{x_{1},...,x_{n}\}$. Then $c_{x_{1}} \cdots c_{x_{n}} = 1$ with $c_{i} \in N_{i}$. Hence $c_{x_{1}}^{-1} = c_{x_{2}} \cdots c_{x_{n}} \in N_{1} \cap \langle \bigcup_{i \neq x_{1}} N_{i} \rangle = \langle 1 \rangle$ and $c_{x_{1}} = 1$. Repeating the process shows $c_i = 1$ for all $i \in \textit{I}$. Hence $\varphi$ is a monomorphism, as desired. $\blacksquare$
 
The difference in structure between the direct product and the weak direct product is one rooted in the difference between finite and infinite structures. The map $\varphi$ in the proof of Proposition \ref{prop2.8}(ii) could not be well-defined in the direct product, as defining the product of an infinite number of group elements requires additional structure. In contrast, the switch to the weak direct product, where only finitely many coordinates are non-identity, along with the crucial condition that each $N_i$ was normal (ensuring commutativity) induced a natural well-defined homomorphism. Thus Proposition \ref{prop2.8}(ii) shows how the presence of infinity in a group structure can cause the structure to fall apart, introducing a clear division between the finite and the infinite in group theory.

Proposition \ref{prop2.8}(i) further highlights a contrast between finite and infinite cyclic groups. In the finite case, a cyclic groups $G$ was isomorphic to $\mbb{Z}/m\mbb{Z}$, whereas in the infinite case $\mbb{Z}$. This classification fully characterizes cyclic groups up to isomorphism and hence displays the connection finite cyclic groups have to the finite structure $\mbb{Z}/m\mbb{Z}$ and infinite cyclic groups have to the infinite structure $\mbb{Z}$. This moreover reflects a kind of continuity between the two as $\mbb{Z}/m\mbb{Z}$ may be viewed as a truncation of $\mbb{Z}$, and in the limit as $m \to \infty$, we have $\mbb{Z} \cong \mbb{Z}/\{0\}$.

Infinite sets tend to appear when studying free abelian groups, and later on in this section, $R$-modules over division rings. Proposition \ref{prop2.8}(ii) showed that given the appropiate structural conditions on a family of normal subgroups, the group was isomorphic to the weak direct product of the subgroups. In algebra, isomorphic can essentially be thought to mean the structures are the same except for a relabeling of their elements under the isomorphism. Thus, we naturally have the following definition
\begin{definition}\label{defsec2.7}
Let $\{N_{i} \mid i \in \textit{I}\}$ be a family of normal subgroups of a group $G$ such that $G = \langle \bigcup_{i \in \textit{I}} N_{i} \rangle$ and for each $k \in \textit{I}$, $N_{k} \cap  \langle \bigcup_{i \neq k} N_{i} \rangle = \langle 1 \rangle$. Then $G$ is called the internal weak direct product of the family $\{N_i \mid i \in \textit{I}\}$, denoted $G = \prod^{w}_{i \in \textit{I}} N_{i}$. If $G$ is abelian, $G$ is called the internal direct sum of $\{N_i \mid i \in \textit{I}\}$ and denoted $G = \bigoplus_{i \in \textit{I}} N_{i}$.
\end{definition}
\begin{theorem}\label{thm2.9}
The following conditions on an abelian group $F$ are equivalent. \\
(i) F has a nonempty basis \\
(ii) F is the internal direct sum of a family of infinite cyclic subgroups. \\
(iii) F is (isomorphic to) a direct sum of copies of $\mbb{Z}$.
\end{theorem}
\noindent $\textbf{Proof}.$ [(i) $\Rightarrow$ (ii)] If $X$ is a basis of $F$, then for each $x \in X$, $nx = 0$ if and only if $n = 0$. It follows that each cyclic subgroup $\langle x \rangle$ is infinite. Since $F = \langle X \rangle$, we also have $F = \langle \bigcup_{x \in X} \langle x \rangle \rangle$. If $z \in X$ is such that $\langle z \rangle \cap \langle \bigcup_{x \in X \setminus \{z\}} \langle x \rangle \rangle \neq 0$, then for some nonzero $n \in \mbb{Z}$, $n_1x_1+...+n_{k}x_k+nz = 0$, contradicting hypothesis that $X$ is a basis. Thus $\langle z \rangle \cap \langle \bigcup_{x \in X \setminus \{z\}} \langle x \rangle \rangle \neq 0$ for every $z \in X$ and $F = \bigoplus_{x \in X} \langle x \rangle$. [(ii) $\Rightarrow$ (iii)] If $F = \bigoplus_{x \in X} \langle x \rangle$, then by Theorem \ref{thm2.7}(i) and Proposition \ref{prop2.8}(i,ii), we have $F \cong \bigoplus_{x \in X} \langle x \rangle \cong \bigoplus_{x \in X}  \mbb{Z}$. [(iii) $\Rightarrow$ (i)] Suppose $F \cong \bigoplus \mbb{Z}$ with the index over $X$. For each $x \in X$, let $\theta_{x}$ be the element $\{u_{i}\}$ of $\bigoplus  \mbb{Z}$ with $u_i = 0$ if $i \neq x$ and $u_x = 1$ and consider $A = \{\theta_{x} \mid x \in X\}$. Given $\{c_{x}\} \in \bigoplus \mbb{Z}$, let $X_0 = \{x \in X|c_x \neq 0\}$ which is finite. Then $\sum_{x \in X_0} c_x \theta_x = \{c_x\}$ and $\bigoplus \mbb{Z} = \langle A \rangle$. Let $Y \subset X$ finite and $d_x \in \mbb{Z}$ for each $x \in Y$. By the definition of $\theta_x$, it follows that if $\sum_{x \in Y} d_x \theta_x = 0$, $d_x = 0$ for each $x$ and thus $A$ is a basis of $\bigoplus \mbb{Z}$.  Let $f:\bigoplus \mbb{Z} \to F$ be an isomorphism and consider $B = \{f(\theta_x) \mid x \in X\}$. Since $A$ generates $\bigoplus \mbb{Z}$, $F = \langle B \rangle$, and if $\sum_{x \in Y} d_x f(\theta_x) = 0$, then by injectivity $\sum_{x \in Y} d_x \theta_{x} = 0$ so $d_x = 0$ for each $x$ and $B$ is a basis of $F$, as desired. $\blacksquare$

For a further proof of Theorem \ref{thm2.9}, see \cite[Theorem II.1.1]{hungerford}.  An abelian group satisfying the equivalent conditions of Theorem \ref{thm2.9} is called a $\textit{free-abelian group}$ (on $X$), with the trivial group $\{0\}$ free-abelian on $\emptyset$. The cardinality of a basis of a free-abelian group is called its $\textit{rank}$ (well-defined by Theorem \ref{thm2.11}). Theorem \ref{thm2.9} shows that (non-trivial) free-abelian groups inherently have infinity embedded within them, for they are isomorphic to infinite cyclic groups. We now come to a result that directly uses infinite cardinals. For further proofs and discussion of Lemma \ref{lem2.10}(i), see \cite[Theorem 0.8.12]{hungerford}; for (ii), see \cite[Theorem 0.16]{rotman}, \cite[Exercise 0.8.12]{hungerford}. For Theorem \ref{thm2.11}, see \cite[Theorem II.1.2]{hungerford} and \cite[Corollary 5.10]{rotman}.
\begin{lemma}\label{lem2.10}
Let $\{A_{i}|i \in \textit{I}\}$ be a family of sets such that $|A_i| \leq \alpha$ for each $i \in \textit{I}$ and $A$ a set such that $|A| = \beta$. Then \\
(i)  $\abs{\bigcup_{n \in \mbb{Z}^{+}} A^n} = \bm{\al}_{0} \beta;$ \\
(ii) $\abs{\bigcup_{i \in \textit{I}} A_i} \leq |\textit{I}| \alpha$.
\end{lemma}
\noindent $\textbf{Proof}.$ [(i)] If $A = \emptyset$, the result is immediate, so assume $A \neq \emptyset$. If $A$ is infinite, then $|A^n| = |A|$ for each $n \in \mbb{Z}^{+}$ (Theorem \ref{thm1.9}). Thus for each $n$ let $f_n:A^n \to A$ be bijective and the map $\bigcup_{n \in \mbb{Z}^{+}} A^n \to \mbb{Z}^{+} \times A$ by $u \mapsto (n,f_n(u))$ ($u \in A^n$) is bijective and $\abs{\bigcup_{n \in \mbb{Z}^{+}} A^n} = |\mbb{Z}^{+}||A| = \bm{\al}_{0} \beta$. If $A$ is finite, for each $n \in \mbb{Z}^{+}$, choose $a_{n} \in A^n$ ($\text{AC}_{\omega}$). Then the map $n \mapsto a_{n}$ is injective so $\bm{\al}_{0}\beta = |\mbb{Z}^{+}| \leq \abs{\bigcup_{n \in \mbb{Z}^{+}} A^n}$ (Theorem \ref{thm1.9}). Moreover, since each $A^n$ is finite, let $g_n:A^n \to \mbb{Z}^{+}$ be an injection and thus the map $\bigcup_{n \in \mbb{Z}^{+}} A^n \to \mbb{Z}^{+} \times \mbb{Z}^{+}$ by $u \mapsto (n,g_n(u))$ ($u \in A^n$) is injective so $\abs{\bigcup_{n \in \mbb{Z}^{+}} A^n} \leq |\mbb{Z}^{+} \times \mbb{Z}^{+}| = \bm{\al}_{0}\beta$ (Theorem \ref{thm1.9}). The result then follows from Schröder–Bernstein (Lemma \ref{lem1.4}). [(ii)] Let $A$ be such that $|A| = \alpha$. If $A = \emptyset$ or $I = \emptyset$, the results immediate, so assume $A \neq \emptyset$ and $I \neq \emptyset$. For each $i \in \textit{I}$, let $f_{i}:A_i \to A$ be an injection. For each $a \in \bigcup_{i \in \textit{I}} A_{i}$, choose one $i \in \textit{I}$ such that $a \in A_{i}$ (AC). Consider the map $f:\bigcup_{i \in \textit{I}} A_i \to \textit{I} \times A$ by $a \mapsto (i, f_{i}(a))$, where $i \in \textit{I}$ is the $i$ we chose such that $a \in A_{i}$. By our choice, $f$ is well defined, and since each $f_{i}$ is injective, so is $f$. Thus $\abs{\bigcup_{i \in \textit{I}}} \leq |\textit{I} \times A| = |\textit{I}|\alpha$, as desired. $\blacksquare$

Lemma \ref{lem2.10}(i) displays a relationship between a set with cardinal number $\beta$ and unions of power of that set, in that it equals $\bm{\al}_{0} \beta$. In the case where $A$ is finite and nonempty, Theorem \ref{thm1.9} shows the unions cardinality is simply $\bm{\al}_{0}$, but the infinite case could be a greater cardinal, and shows how countable unions of sets relate to infinity. This highlights a central theme of infinity, that operations like union and product behave differently depending on whether the underlying sets are finite or infinite. Then (ii) demonstrates how unions of a family of sets all bounded above by some cardinal $\alpha$  are bounded above by the cardinality of the index times $\alpha$, and shows how a result that would trivially hold in a finite case naturally extends to infinity. It shows how in this case, cardinal inequalities remain coherent when extended to infinity families of sets and infinite cardinalities, preserving bounds in a way that becomes crucial in 
\begin{theorem}\label{thm2.11}
Any two bases of a free abelian group $F$ have the same cardinality.
\end{theorem}
\noindent $\textbf{Proof}.$ If $F$ has $\emptyset$ as a basis, then $F = \{0\}$ and clearly all other bases are empty and hence the same cardinality. Otherwise, suppose $F$ has basis $X$ with $|X| = n$ ($n \in \mbb{Z}^{+}$) and then $F \cong \bigoplus_{i = 1}^{n} \mbb{Z}$ by the proof of Theorem \ref{thm2.9}. Let $f:F \xrightarrow{\cong} \bigoplus_{i = 1}^{n} \mbb{Z}$. We see $f(2F) = 2f(F) = 2\bigoplus_{i = 1}^{n} \mbb{Z}$ and thus $2F \cong \bigoplus_{i = 1}^{n} 2\mbb{Z}$, so $F/2F \cong \bigoplus_{i = 1}^{n} \mbb{Z}/2\mbb{Z}$ by Theorem \ref{thm2.7}(ii) and $|F/2F| = 2^{n}$. If $Y$ is another basis of $F$ (which need not be finite) and $r$ any positive integer such that $r \leq |Y|$, then there exists a monomorphism $g:\bigoplus_{i = 1}^{r} \mbb{Z} \to F$ (since $F$ is isomorphic to $|Y|$ copies of $\mbb{Z}$) and similarly a monomorphism $\bigoplus_{i = 1}^{r} \mbb{Z}/2\mbb{Z} \to F/2F$. Thus $2^r \leq |F/2F|$, so $2^r \leq 2^n$ and $r \leq n$. It follows that $|Y| = m$ is finite so by above $|F/2F| = 2^m$. Therefore, $|X| = n = m = |Y|$.

If one basis of $F$ is infinite, so are all other bases by above. Thus it suffices to show for any infinite basis $X$ of $F$, $|X| = |F|$. Clearly $|X| \leq |F|$. Let $S = \bigcup_{n \in \mbb{Z}^{+}} X^{n}$, and for each $s = (x_1,...,x_n) \in S$ let $G_{s}$ be the group generated by $\langle x_1,...,x_n \rangle$. Let $y_1,...,y_t$ ($t \leq n$) be the distinct elements of $\{x_1,...,x_n\}$. Since $X$ is a basis, it follows that the map $G_s \to \bigoplus_{i = 1}^{t} \mbb{Z} y_i$ by $\sum_{i = 1}^{t} n_iy_i \mapsto (n_1y_1,...,n_ty_t)$ is a well-defined isomorphism and $G_s \cong \bigoplus_{i = 1}^{t} \mbb{Z} y_i$. It then follows that $|G_s| = |\mbb{Z}^t| = \bm{\al}_0$ (Theorem \ref{thm1.9}). Since $F = \bigcup_{s \in S} G_s$, by Lemma \ref{lem2.10}(ii) $|F| = \abs{\bigcup_{s \in S} G_s} \leq |S| \bm{\al}_0$. By Lemma \ref{lem2.10}(i) and Theorem \ref{thm1.9}, $|S| = \abs{\bigcup_{n \in \mbb{Z}^{+}} X^{n}} = |X|$ and thus $|F| \leq |X| \bm{\al}_0 = |X|$. Hence $|F| = |X|$ by Lemma \ref{lem1.4}, as desired. $\blacksquare$

The presence of infinity in the proof of Theorem \ref{thm2.11} is not just in the case where every basis of $F$ is infinite, but also when $F$ has a finite basis $X$. In the proof, we cannot assume another basis must also be finite, hence we considered integers at most the cardinality $|Y|$ rather than $|Y|$ itself. The infinite case then relies on the fact that the unions of the $G_s$'s ($s \in S$) is bounded above by $|S|\bm{\al}_0$ and that the cardinality of the unions of the $X^n$'s ($n \in \mbb{Z}^{+}$) equals the cardinality of $X$, inducing both Lemma \ref{lem2.10}(i) and (ii) to then create a cardinality squeeze and invoking Schröder–Bernstein (Lemma \ref{lem1.4}), establishing a relationship between how infinity interacts between set-theoretic mathematics to algebraic mathematics. The proposition below is stated in  \cite[Theorem II.1.4]{hungerford}, but the proof here follows a different strategy.
\begin{proposition}\label{prop2.12}
Every abelian group $G$ is the homomorphic image of a free abelian group of rank $|X|$, where $X$ is a set of generators of $G$.
\end{proposition}
$\textbf{Proof}.$ Let $F = \bigoplus_{x \in X} \mbb{Z}x$. By the proof of Theorem \ref{thm2.9}, $\{\theta_{x} \mid x \in X\}$ is a basis of $F$ with rank $|X|$. Then since $\{\theta_{x} \mid x \in X\}$ is a basis, the map $f:F \to G$ by $\sum_{i = 1}^{k} n_k \theta_{x_1} \mapsto \sum_{i = 1}^{k} n_k x_k$ ($x_i \in X, n_i \in \mbb{Z}$) is a well-defined epimorphism and $\Im f = G$, as desired. $\blacksquare$

Proposition \ref{prop2.12} shows how infinity relates to abelian groups in general. While not all (non-trivial) abelian groups are free-abelian, and thus related to infinity by being the (internal) direct sum of a family of infinite cyclic subgroups, they are the homomorphic image of such a group, hence establishing a natural connection between infinity and more general abelian groups. We now present a crucial result, the $\textit{Fundamental Theorem of Finitely Generated Abelian Groups}$. Only the existence statement is proven here (Theorem \ref{thm2.14}), but without too much more theory (see \cite[Section II.2]{hungerford}), uniqueness may be proved as well. For further discussion of Lemma \ref{lem2.13}, see \cite[Theorem II.1.6]{hungerford}; for Theorem \ref{thm2.14}, see \cite[Theorem 5.3]{dummit} and  \cite[Theorem II.2.1]{hungerford}.
\begin{lemma}\label{lem2.13}
If $F$ is a free abelian group of finite rank $n \in \mbb{Z}^{+}$ and $G$ is a non-trivial subgroup of $F$, then there exists a basis $\{x_1,...,x_n\}$ of $F$, and positive integers $r, d_1,..,d_r$ ($1 \leq r \leq n$) such that $d_1 \mid d_2 \mid \cdots \mid d_r$ and $G$ is free abelian with basis $\{d_1x_1,...,d_rx_r\}$.
\end{lemma}
\noindent $\textbf{Proof}.$ If $F$ has rank 1, let $\{x_1\}$ be a basis and $F = \langle x_1 \rangle \cong \mbb{Z}$ (Proposition \ref{prop2.8}(i)). If $G \subset F$ ($G \neq \{0\}$), let $d_1$ be the smallest positive integer such that $d_1x_1 \in G$. Given $cx_1 \in G$, write $c = kd_1+r$ ($k,r \in \mbb{Z}$, $1 \leq r < d_1$). Then $rx_1 \in G$ and thus by minimality $r = 0$, whence $G = \langle d_1 x_1 \rangle \cong \mbb{Z}$ (Proposition \ref{prop2.8}(i)). Suppose, by way of induction, the result holds for all ranks less than $n$ ($n > 1$). Let $S$ be the set of all integers $s \neq 0$ such that there exists a basis $\{y_1,...,y_n\}$ of $F$ and an element in $G$ of the form $sy_1+\sum_{i = 2}^{n} k_iy_i$ ($k_i \in \mbb{Z}$). Since $\{y_2,y_1,y_3,...,y_n\}$ is also a basis, $k_2 \in S$, and similarly $k_j \in S$ for $j = 2,3,...,n$. Since $G \neq \{0\}$, $S \neq \emptyset$, and by the well-ordering principle, there exists a least $d_1 \in \mbb{Z}^{+}$ (since $G$ has additive inverses) such that for some basis $X = \{z_1,...,z_n\}$ of $F$ and $u \in G$, $u = d_1z_1+\sum_{i = 2}^{n} k_iz_i$. Write $k_i = d_1q_i+r_i$ ($0 \leq r_i < d_1$) for each $2 \leq i \leq n$, and $u = d_1(z_1+\sum_{i = 2}^{n} q_iz_i)+\sum_{i = 2}^{n} r_iz_i$. Let $x_1 = z_1+\sum_{i = 2}^{n} q_iz_i$ and $Y = \{x_1,z_2,...,z_n\}$. Since $z_1 = x_1-\sum_{i = 2}^{n} q_iz_i$, $Y$ generates $F$, and since $X$ is a basis, it then follows that $Y$ is. By the minimality of $d_1$, since $Y$ is a basis in any order, $r_2=r_3= \cdots = r_n = 0$ and thus $u = d_1x_1 \in G$. 

Let $H = \langle z_2,z_3,...,z_n \rangle$. Then $H$ is free-abelian with rank $n-1$ and $F = \langle x_1 \rangle \bigoplus H$. Since $\{x_1,z_2,...,z_n\}$ is a basis of $F$, $\langle u \rangle \cap (G \cap H) = \{0\}$. If $v = c_1x_1+\sum_{i = 2}^{n} c_iz_i \in G$ ($c_i \in \mbb{Z}$), write $c_1 = d_1s_1+r_1$ with $0 \leq r_1 < d_1$. Thus $v-s_1u = r_1x_1+\sum_{i = 2}^{n} c_iz_i$ and by the minimality of $d_1$, $r_1 = 0$ and $\sum_{i = 2}^{n} c_i z_i \in G \cap H$.  Thus $v = s_1u+\sum_{i = 2}^{n} c_i z_i$ and $G = \langle \langle u \rangle \cup (G \cap H) \rangle$. Thus by Definition \ref{defsec2.7}, $G = \langle u \rangle \bigoplus (G \cap H)$.

If $G \cap H = \{0\}$, then $G = \langle d_1x_1 \rangle$ and we are done. Otherwise, by induction hypothesis there exists a basis $\{x_2,...,x_n\}$ of $H$ and positive integers $r, d_2,...,d_r$ such that $d_2 \mid d_3 \mid \cdots d_r$ and $G \cap H$ is free abelian with basis $\{d_2x_2,...,d_rx_r\}$. Since $F = \langle x_1 \rangle \bigoplus H$ and $G = \langle d_1x_1 \rangle \bigoplus (G \cap H)$, it follows that $\{x_1,x_2,...,x_n\}$ is a basis of $F$ and $\{d_1x_1,...,d_rx_r\}$ of $G$ (since $x_1 =  z_1+\sum_{i = 2}^{n} q_iz_i$). Thus, to complete the proof, we show $d_1 \mid d_2$. Write $d_2 = qd_1+r_0$ with $0 \leq r_0 < d_1$. Since $\{x_1,..,x_n\}$ is a basis of $F$, it follows that $\{x_2,x_1+qx_2,...,x_n\}$ is also a basis and $r_0x_2+d_1(x_1+qx_2) = d_1x_1+d_2x_2 \in G$, the minimality of $d_1$ in $S$ implies that $r_0 = 0$ and $d_1 \mid d_2$, as desired. $\blacksquare$

Lemma \ref{lem2.13} shows how (non-trivial) subgroups of free-abelian groups relate to infinity, primarily the fact that they are free-abelian themselves, and hence internal direct sums of infinite cyclic groups, or isomorphic to copies of the infinite abelian group $\mbb{Z}$ (Theorem \ref{thm2.9}).
\begin{theorem}\label{thm2.14}
Every finitely generated abelian group $G$ is isomorphic to a finite direct sum of (possibly empty) additive groups $\mbb{Z}/m_k\mbb{Z}$ ($1 \leq k \leq t$), where $m_1,...,m_t \in \mbb{Z}^{+}$ ($m_1 > 1$) satisfy $m_1 \mid m_2 \mid \cdots m_t$, and a (possibly trivial) free abelian group $F$ of finite rank.
\end{theorem}
\noindent $\textbf{Proof}.$ If $G = \{0\}$, then $G$ is trivially free-abelian. Otherwise, if $G$ is generated by $n \in \mbb{Z}^{+}$ elements, there is a free abelian group $S$ and an epimorphism $f:S \to G$ (Proposition \ref{prop2.12}). If $f$ is an  isomorphism, $G \cong S$. If not, then by Lemma \ref{lem2.13}, there is a basis $\{x_1,...,x_n\}$ of $S$ and positive integers $r,d_1,...,d_r$ ($1 \leq r \leq n$), $d_1 \mid d_2 \mid \cdots d_r$, and $\{d_1x_1,...,d_rx_r\}$ is a basis of $K = \ker f$. Thus $S = \bigoplus_{i = 1}^{n} \langle x_i \rangle$ and $K = \bigoplus_{i = 1}^{r} \langle d_ix_i \rangle$, where $\langle x_i \rangle \cong \mbb{Z}$ (Proposition \ref{prop2.8}(i)) and under the same isomorphism $\langle d_i x_i \rangle \cong d_i \mbb{Z}$. For $r+1 \leq i \leq n$, let $d_i = 0$, and then by Theorem \ref{thm2.6}(i,ii), Theorem \ref{thm2.7}(ii), $$G \cong S/K = \bigoplus_{i = 1}^{n} \langle x_i \rangle \Big/ \bigoplus_{i = 1}^{n} \langle d_ix_i \rangle \cong \bigoplus_{i = 1}^{n} \langle x_i \rangle/ \langle d_ix_i \rangle \cong \bigoplus_{i = 1}^{n} \mbb{Z}/d_i \mbb{Z}.$$
If $d_i = 1$, then $\mbb{Z}/d_i\mbb{Z} = \mbb{Z}/\mbb{Z} = \{0\}$, likewise if $d_i = 0$, $\mbb{Z}/d_i\mbb{Z} = \mbb{Z}/\{0\} \cong \mbb{Z}$. Let $m_1,...,m_t$ be the $d_i$ in order such that $d_i \neq 0,1$, and $s$ the number of $d_i$ with $d_i = 0$. Then we have $$G \cong \mbb{Z}/m_1 \mbb{Z} \oplus \cdots \oplus \mbb{Z}/m_t \mbb{Z} \oplus F.$$
where $m_1 > 1$, $m_1 \mid m_2 \mid \cdots m_t$ and $F$ is free abelian of rank $s$, as desired. $\blacksquare$

In the case where $G$ is finite, the free abelian group must be trivial, and thus we have \begin{align}\label{test}
G \cong \mbb{Z}/m_1 \mbb{Z} \oplus \cdots \oplus \mbb{Z}/m_t \mbb{Z}
\end{align}
The decomposition in (\ref{test}) is called the $\textit{fundamental theorem of finite abelian groups}$. It can be proven independently to Theorem \ref{thm2.14} as well using $p$-groups (see \cite{dummit, navarro, rotman}). The fundamental theorems of finite and finitely generated abelian groups illustrate the structural differences between finite and infinite finitely generated groups. In the purely finite case, a classification in terms of finitely many additive quotient groups $\mbb{Z}/m_i \mbb{Z}$. In contrast, infinity introduces a (non-trivial) free-abelian group into the decomposition, highlighting that the fundamental distinction between finite and infinite finitely generated abelian groups lies in the presence of a free abelian component. This structure is inherently infinite due to Theorem \ref{thm2.9}.

As shown above, the theory of free-abelian groups in general sheds light on the role of infinity in abelian group theory. Theorems \ref{thm2.9} and \ref{thm2.11} illustrate how free abelian groups themselves embody infinity, while Proposition \ref{prop2.12} shows how general abelian groups relate. Lemma \ref{lem2.13} highlights how nontrivial subgroups of free-abelian groups relate to infinity, and Theorem \ref{thm2.14} establishes the connection to infinity within finitely generated abelian groups.
\subsection{Ring and Modules}
The study of rings and  $R$-modules shows once again structural differences between the finite and the infinite in algebra. Groups form a structure much simpler than that of rings, which add an additional operation, whereas modules form a generalization of abelian groups (which are modules over $\mbb{Z}$). Naturally, we begin by defining a ring.
\begin{definition}\label{def3.1}
A ring $R$ is a non-empty set, together with two binary operations (denoted as addition and multiplication here) with the following properties. \\
(i) ($R$, $+$) is an abelian group \\
(ii) ($R$, $\cdot$) is a monoid \\
(iii) $a(b+c) = ab+ac$ and $(a+b)c = ac+bc$ for all $a,b,c \in R$ \\
A ring $R$ with $ab = ba$ for all $a,b \in R$ is called a commutative ring.
\end{definition} 
\begin{example}\label{ex17}
The set $\mbb{Z}$ of integers is a commutative ring under the usual definitions of addition and multiplication, and similarly for $\mbb{Z}/n\mbb{Z}$ with multiplication defined by $(a+n\mbb{Z})(b+n\mbb{Z}) = ab+n\mbb{Z}$. The set $M_{n}(\mbb{R})$ of all $n$-by-$n$ matrices with entries in $\mbb{R}$ is a ring under the usual addition and multiplication and is non-commutative for $n \geq 2$.
\end{example}
This definition aligns with \cite{lang, jacob1}, although some authors \cite{dummit, hungerford} do not require that rings have a multiplicative identity. Let $R$ be a ring. A $\textit{subring}$ of $R$ is a subset $S$ of $R$ that's a ring under the same addition and multiplication as $R$, with the same multiplicative identity.  A non-zero element $a \in R$ is called a $\textit{zero-divisor}$ if there exists $b \in R$ ($b \neq 0$) such that $ab = 0$ or $ba = 0$. An element $c \in R$ is left [resp. right] invertible if there exists $x \in R$ [resp. $y \in R$] such that $xc = 1$ [resp. $cy = 1$]. If $c \in R$ is both left and right invertible, $c$ is called a $\textit{unit}$. 

An $\textit{integral domain}$ is a commutative ring $R$ with $1 \neq 0$ and no zero divisors. A ring $D$ with $1 \neq 0$ such that every element is a unit is called a $\textit{division ring}$. A commutative division ring is called a $\textit{field}$. We now prove some basic results about rings. For further proofs and discussion of the following theorem, see \cite[Section 2.1]{jacob1}, \cite[Theorem III.1.2]{hungerford}.
\begin{theorem}\label{thm2.15}
Let $R$ be a ring. Then \\
(i) $0a = a0 = 0$ for all $a \in R$; \\
(ii) $(-a)b = a(-b) = -(ab)$ for all $a,b \in R$; \\
(iii) $(-a)(-b) = ab$; \\
(iv) 
$$\qty(\sum_{i = 1}^{n} a_i)\qty(\sum_{j = 1}^{m} b_j) = \sum_{i = 1}^{n} \sum_{j = 1}^{m} a_i b_j \hspace{4mm} \text{for all} \hspace{3mm} a_i, b_j \in R;$$
(v) If $R$ is an integral domain, then left and right cancelation hold. That is, if $ac = bc$ or $ca = cb$ with $a,b,c \in R$ and $c \neq 0$, then $a = b$.
\end{theorem}
\noindent $\textbf{Proof}.$ Let $a,b \in R$. [(i)] $0a = (0+0)a = 0a+0a$, thus $0a = 0$. Similarly, $a0 = 0$. [(ii)] $ab+(-a)b = (a+(-a))b = 0b = 0$, thus $(-a)b = -(ab)$. Similarly, $a(-b) = -(ab)$.  [(iii)] By (ii), $(-a)(-b) = -(a(-b)) = -(-(ab)) = ab$. [(iv)] We proceed by induction on $m$. If $m = 1$, the result follows immediately by definition. By way of induction, assume $m > 1$ and the result holds for $m-1$. Then \begin{align*}
\qty(\sum_{i = 1}^{n} a_i)\qty(\sum_{j = 1}^{m} b_j)  & = \qty(\sum_{i = 1}^{n} a_i)\qty(\sum_{j = 1}^{m-1} b_j+b_m) = \sum_{i = 1}^{n} \sum_{j = 1}^{m-1} a_ib_j+\sum_{i =1}^{n} a_ib_m \\  & = \sum_{i = 1}^{n} \sum_{j = 1}^{m} a_i b_j,
\end{align*}
(v) If $ac = bc$, then $(a-b)c = 0$. If $a \neq b$, then $a-b \neq 0$ and $c \neq 0$ is a zero divisor, thus $a = b$. Similarly, $ca = cb$ implies $c(a-b) = 0$ and $a = b$, as desired. $\blacksquare$

Let $R$ be a ring and $I$ a subset of $R$ with $(I, +)$ a subgroup of $(R, +)$. Then $I$ is a $\textit{left ideal}$ [resp. $\textit{right ideal}$] of $R$ if $r \in R$ and $x \in I$ implies that $rx \in I$ [resp. $xr \in I$]. If $I$ is both a left and right ideal, $I$ is called an $\textit{ideal}$ of $R$. We now prove a useful lemma about ideals
\begin{lemma}\label{lem2.16}
Let $I$ be a nonempty subset of a ring $R$. Then $I$ is an ideal of $R$ if and only if $a-b \in I$ for every $a,b \in I$ and $ra, ar \in I$ for every $a \in I, r \in R$.
\end{lemma}
\noindent $\textbf{Proof}.$ Suppose $I$ is an ideal. Then $a-b \in I$ for every $a,b \in I$ since $I$ is a subgroup of $R$. The second property is immediate by definition. Conversely, if $a-b \in I$ for every $a,b \in I$, then $a-a = 0 \in I$, $a+b = a-(-b) \in I$, and $0-a = -a \in I$ for every $a,b \in I$. Thus $I$ is a subgroup, and by the second property $I$ is an ideal, as desired. $\blacksquare$ 

An ideal [resp. left ideal] $M$ in a ring $R$ is called $\textit{maximal}$ if $M \neq R$ and every ideal [resp. left ideal] $N$ such that $M \subset N \subset R$, either $N = M$ or $N = R$. Further proofs and discussion of the following result can be found in \cite{hungerford, conrad, macdonald}. The result was originally proven by Krull \cite{krull}, and is equivalent to the axiom of choice \cite{hodges}. Note that the statement holds for left and right ideals as well (as in \cite{hungerford}).
\begin{theorem}\label{thm2.17}
In a nonzero ring, maximal ideals always exist.
\end{theorem}
\noindent $\textbf{Proof}.$ Since $\{0\}$ is (trivially) an ideal and $\{0\} \neq R$, it suffices to prove that every ideal ($\neq R$) in $R$ is contained in a maximal ideal of $R$. Given an ideal $A$ in $R$ ($A \neq R$) let $\mathcal{S}$ be the set of all ideals $B$ in $R$ such that $A \subset B \neq R$. $\mathcal{S} \neq \emptyset$ since $A \in \mathcal{S}$. Partially order $\mathcal{S}$ by set-theoretic inclusion, and let $C = \{C_i| i \in \textit{I}\}$ be a chain in $\mathcal{S}$. Consider $D = \bigcup_{i \in \textit{I}} C_i$. If $a,b \in D$, then $a \in C_i$, $b \in C_j$ for some $i,j \in \textit{I}$. Without loss of generality, assume $C_i \subset C_j$. Then $a,b \in C_j$. Since $C_j$ is an ideal, it follows that $D$ is (Lemma \ref{lem2.16}). Thus, by Zorn's Lemma $\mathcal{S}$ has a maximal element $M$. If $N$ is an ideal such that $N \in \mathcal{S}$ and $A \subset M \subset N \subset R$, then by definition of maximal element, $N = M$, and if $N \notin \mathcal{S}$, by definition of $\mathcal{S}$, $N = R$. Thus $M$ is a maximal ideal, as desired. $\blacksquare$

The proof of Theorem \ref{thm2.17} induces infinity through the use of Zorn's Lemma, which is equivalent to the axiom of choice. The axiom of choice says that any set $X$ has a \textit{choice function}, that is a function with domain $X \setminus \{\emptyset\}$ such that $f(y) \in y$ for all $y \in X \setminus \{\emptyset\}$ \cite{shoenfieldaxioms}. This axiom is particularly relevent in the context of infinite sets, as it ensures such a choice function holds for infinite sets, where choosing an element from each element of $X$ may not be constructively possible. The fact that maximal ideals exist for any non-zero ring is a crucial result, and highlights how infinity intertwines deeply in the theory of rings. Theorem \ref{thm2.17} thus serves as a concrete example where the set-theoretic axioms pertaining to infinity profoundly underpin algebraic structures. For further proofs and discussion of the following theorem, see \cite[Corollary II.3]{dummit} and \cite[Exercise III.1.6]{hungerford}.
\begin{theorem}\label{thm2.18}
Every finite integral domain is a field.
\end{theorem}
\noindent $\textbf{Proof}.$ Let $R$ be a finite integral domain and $a \in R$ ($a \neq 0$). Consider the map $f:R \to R$ by $x \mapsto ax$. By Theorem \ref{thm2.15}(v), $f$ is well-defined and injective. Let $k = |\Im f| \leq |R| = n$. If $f$ is not surjective, take $x_0 \notin R$ and $a_0 \in R \setminus \Im f$. Extend $f$ to $f'$ by mapping $\{x_0\} \mapsto a_0$. Then $f'$ is injective and $|R \cup \{x_0\}| = n+1 > k$, contradicting the pigeonhole principle. Thus $f$ is bijective. In particular, there exists $b \in R$ such that $ab = 1$, and $a$ is a unit. Since $a$ was arbitrary, $R$ is a field, as desired. $\blacksquare$

Naturally, the finite condition in Theorem \ref{thm2.18} is absolutely crucial. The ring $\mbb{Z}$ of integers, for example, is an infinite integral domain that is not a field. The proof above relies on the fact that an injective function $f:A \to A$ is a bijection, which is only true for finite sets (e.g., the map $f:\mbb{N} \to \mbb{N}$ by $n \mapsto n+1$ fails). Thus, Theorem \ref{thm2.18} highlights a disconnect between the finite and the infinite in ring theory. We now introduce a useful ring
\begin{definition}[\cite{hungerford}, Theorem III.5.1]\label{def2.27}
Let $R$ be a ring. The set $R[x]$ of all sequences of  $(a_0,a_1,...)$ of elements of $R$ such that $a_i = 0$ for all but a finite number of indices $i$ along with addition and multiplication defined by $(a_0,a_1,...)+(b_0,b_1,...) = (a_0+b_0,a_1+b_1,...)$
and $(a_0,a_1,...)(b_0,b_1,...) = (c_0,c_1,...)$, where $c_n = \sum_{i = 0}^{n} a_{n-i}b_i = \sum_{k+j = n} a_kb_j$ is called the $\textit{ring of polynomials}$ over $R$.
\end{definition}
We now verify $R[x]$ is a ring for a ring $R$. Let $a = (a_0,a_1,...),b = (b_0,b_1,...), c = (c_0,c_1,...) \in R[x]$. Clearly addition is commutative, associative, and the additive identity is $0 = (0,0,...)$, with additive inverse of $a$ being $(-a_0,-a_1,...)$. Thus $R[x]$ is an abelian group. The multiplicative identity is naturally $1 = (1_R,0,0,...)$ for at each $n$ $a_n = \sum_{i = 0}^{n} a_{n-i}1_i = \sum_{i = 1}^{n} 1_{n-i} a_i$, and it is easy to see distributivity (left and right) holds. We see \begin{align*}
((ab)c)_{n} & = \sum_{k+j = n} (ab)_k c_j = \sum_{k+j = n} \qty(\sum_{p +q = k} a_p b_q)c_j = \sum_{k+j = n} \sum_{p+q = k} a_pb_qc_j = \sum_{p+q+j = n} a_pb_qc_j \\ & = \sum_{p+m = n} a_p \qty(\sum_{q+j = m} b_q c_j) \\ & = \sum_{p+m = n} a_p(bc)_m \\ & = (a(bc))_n
\end{align*}
Hence multiplication is associative and $R[x]$ a ring.Naturally, if $R$ is commutative, so is $R[x]$. 

A \textit{ring-homomorphism} is a group homomorphism under addition and a monoid homomorphism under multiplication. Consider the map $f:R \to R[x]$ by $r \mapsto (r,0,0,...)$. Clearly $f$ is well-defined and injective. Moreover, if $a,b \in R$, $f(a+b) = (a+b,0,0,...) = (a,0,0,...)+(b,0,0,...) = f(a)+f(b)$, $f(ab) = (ab,0,0,...) = (a,0,0,...)(b,0,0,...) = f(a)f(b)$, and $f(1_{R}) = (1_{R},0,0,...) = 1_{R[x]}$, so $f$ is a monomorphism of rings. In light of this monomorphism, we identity $R$ with its isomorphic image in $R[x]$ and write $(r,0,0,...)$ as simply $r$ ($r \in R$). Notice $r(a_0,a_1,...) = (ra_0,ra_1,...)$ for $(a_0,a_1,...) \in R[x]$. For further discussion of the following proposition, see \cite[Section 2.10]{jacob1}, \cite[Theorem III.5.2]{hungerford}, and \cite[Section II.3]{lang}.
\begin{proposition}\label{prop2.28}
Let $R$ be a ring. Define $x^n = (0,0,...,0,1_{R},0,...)$ where $1_{R}$ is the (n+1)th coordinate ($n \in \mbb{N}$). Then for every $f \in R[x]$ there exists unique $n \in \mbb{N}$ and elements $a_0,...,a_n \in R$ such that $f = \sum_{i = 0}^{n} a_ix^{i}$ in the sense that if $f = \sum_{i = 0}^{m} b_i x^i$ ($b_i \in R$), then $a_i = b_i$ for $0 \leq i \leq n$, and if $m > n$, then $b_i = 0$ for all $n<i \leq m$.
\end{proposition}
\noindent $\textbf{Proof}$. If $f = (a_0,a_1,...) \in R[x]$, there must exist a largest index $n$ such that $a_n \neq 0$ (otherwise $f \notin R[x]$). Thus $f = \sum_{i = 0}^{n} a_ix^{i}$. If $f = \sum_{i = 0}^{m} b_ix^i$ and $m < n$, then the nth coordinate of the two expressions differ, so $m \geq n$. Similarly, If $a_i \neq b_i$ for some $1 \leq i \leq n$ or $m > n$ and $b_i \neq 0$ for some $n<i \leq m$, the ith coordinate of the two expressions differ, a contradiction. Thus $b_i = a_i$ for $1 \leq i \leq n$ and if $m > n$ then $b_i = 0$ for all $n<i \leq m$, as desired. $\blacksquare$

In light of Proposition \ref{prop2.28}, we write a polynomial $f \in R[x]$ in a familiar form as a sum $f = \sum_{i = 0}^{n} a_ix^{i}$.

Given a ring $R$, the set of all sequences of elements of $R$ under the same addition and multiplication as Definition \ref{def2.27} forms a generaliztion of $R[x]$, denoted $R[[x]]$, and called the \textit{ring of formal power series}. The ring of formal power series, introducing infinity, allows for representations of more general functions. For example, in the ring $\mbb{R}$ of real numbers, it is well know that $f(x) = e^x$ has representation (see \cite[Chapter 20]{spivak}) \begin{align}\label{l45}
e^x = \sum_{n = 0}^{\infty} \frac{x^n}{n!}.
\end{align}
Thus, the ring of formal power series $R[[x]]$ is an example of a natural extension from the finite in $R[x]$, where only finitely many coordinate in $f \in R[x]$ are non-zero, to the infinite in $R[[x]]$, whose representations extend naturally from finite sums of $\sum_{i = 0}^{n} a_i x^i$, to infinite sums $\sum_{i = 0}^{\infty} a_ix^{i}$. This extension is conceptually similar to the extension from weak direct products to direct products. That extension from the finite to the infinite is natural in a similar way, however as we saw in Proposition \ref{prop2.8}, there are some consequences of extending to potentially infinitely many non-zero coordinates. There are two major consequences with extending to formal power series. For one, a function representation in $R[[x]]$ may have multiple representations. For example, as in (\ref{l45}) $$e^x = \sum_{n = 0}^{\infty} \frac{x^n}{n!} = e\sum_{n = 0}^{\infty} \frac{(x-1)^n}{n!}.$$
For two, such representation need not converge for all $x \in R$. For example, the representation $$\frac{1}{1-x} = \sum_{n = 0}^{\infty} x^n.$$
it well known to converge only for $|x| < 1$ (see \cite[Section 9.9]{larson}). These two examples thus highlight that, in an analogous way to Proposition \ref{prop2.8}, the structural differences between polynomial rings and formal power series illustrate a division between the finite and the infinite in ring theory.

Given a (commutative) ring $R$ and $X \subset R$, the ideal \textit{generated by $X$} is defined as $(X) = RX = \{\sum_{i = 1}^{n} r_ix_i|r_i \in R, x_i \in X\}$, with the empty set generating the zero ideal (for general definitions, see \cite[Section 7.4]{dummit}). It is immediate that $(X)$ is indeed an ideal (Lemma \ref{lem2.16}). The ideal $(a)$ for $a \in R$ is called a \textit{principal ideal}. If every ideal in $R$ is principal, $R$ is called a \textit{principal ideal ring}. If $R$ is an integral domain, $R$ is called a \textit{principal ideal domain}.
\begin{example}\label{ex2.200}
The set $\mbb{Z}$ of integers is a principal ideal domain. First off, clearly it is an integral domain. Let $\textit{I}$ be an ideal of $\mbb{Z}$. If $\textit{I} = \{0\}$, then $\textit{I} = (0)$. If not, by the well-ordering principle let $a$ be the least positive integer in $\textit{I}$. Since $\textit{I}$ is an ideal, $(a) \subset \textit{I}$. Moreover, given $b \in \textit{I}$, by the division algorithm, $b = qa+r$ with $q,r \in \mbb{Z}$ and $0 \leq r < a$. Since $r = b-qa \in \textit{I}$, by the minimality of $a$ we have $r = 0$ and $b = qa \in (a)$, whence $\textit{I} = (a)$ and $\mbb{Z}$ a principal ideal domain. 
\end{example} 
For further proofs and discussion of the following lemma, see \cite[Lemma 6.18]{rotman}, and for the theorem, see \cite[Proposition 6.1-2]{macdonald}; \cite[Proposition 6.38]{rotman}.
\begin{lemma}\label{lemma2.29}
Let $R$ be a commutative ring and $I_1 \subset I_2 \subset...$ an ascending chain of ideals in $R$. Then $J = \bigcup_{n \in \mbb{Z}^{+}} I_n$ is an ideal in $R$
\end{lemma}
\noindent $\textbf{Proof}.$ If $a \in J$, then $a \in I_n$ for some $n \in \mbb{Z}^{+}$. Thus since $I_n$ is an ideal, $ra \in I_n \subset J$ for all $r \in R$. Similarly, if $a,b \in J$, $a,b \in I_{k}$ for some $k \in \mbb{Z}^{+}$. Thus $a-b \in I_{k} \subset J$ and $J$ is an ideal by Lemma \ref{lem2.16}. $\blacksquare$
\begin{theorem}\label{thm2.30}
Let $R$ be a commutative ring. Then each of the following conditions are equivalent. \\
(i) Every ascending chain of ideals $I_1 \subset I_2 \subset...$ in $R$ is stationary, that is there exists $N \in \mbb{Z}^{+}$ such that $I_N = I_{N+1}=I_{N+2}...$. \\
(ii) Every non-empty family of ideals in $R$ has a maximal element. \\
(iii) Every ideal in $R$ is finitely generated. 
\end{theorem}
\noindent $\textbf{Proof}$. [(i) $\Rightarrow$ (ii)] Suppose, by way of contradiction, $\mathcal{F}$ is a family of ideals with no maximal element. Choose $I_1 \in \mathcal{F}$. Since $I_1$ is not maximal, there exists $I_2 \in \mathcal{F}$ such that $I_1 \subsetneq I_2$. Continuing this way, we construct a chain of ideals $I_1 \subsetneq I_2 \subsetneq I_3 \subsetneq...$ that is not stationary, contradicting (i). [(ii) $\Rightarrow$ (iii)] Let $I$ be an ideal in $R$ and $\mathcal{F}$ the family of finitely generated ideals of $R$ contained in $I$. $\mathcal{F} \neq \emptyset$ since $\{0\} \in \mathcal{F}$. Thus $\mathcal{F}$ has maximal element $M$. Notice $M \subset I$ because $M \in \mathcal{F}$. If $M \neq I$, there exists $x \in I$ such that $x \notin M$. The ideal $J = \{m+rx \mid m \in M, r \in R\}$ is finitely generated since $M$ is, and moreover $J \in \mathcal{F}$ with $M \subsetneq J$, contradicting the maximality of $M$. Thus $M = I$ and $I$ is finitely generated. [(iii) $\Rightarrow$ (i)] Let $I_1 \subset I_2 \subset ...$ be a chain of ideals in $R$. By Lemma \ref{lemma2.29} $J = \bigcup_{n \in \mbb{Z}^{+}} I_n$ is an ideal. Thus, $J = (a_1,...,a_k)$ for some $a_i \in J$ ($1 \leq i \leq k$). Since for each $i$, $a_i \in I_{n_i}$ for some $n_i$, let $N = \max\{n_1,...,n_k\}$. Thus it follows that $J \subset I_{N}$ and $I_n \subset J \subset I_{N}$ for every $n \in \mbb{Z}^{+}$. Thus $I_{N} = I_{N+1} =...$ and the chain is stationary, as desired. $\blacksquare$

A commutative ring satisfying any of the equivalent conditions of Theorem \ref{thm2.30} is called a \textit{Noetherian Ring}. The condition (i) in Theorem \ref{thm2.30} is called the \textit{ascending chain condition} (ACC).

Noetherian rings, by ACC, do not permit an infinitely ascending chain of ideals. The finiteness moreover manifests in the crucial equivalent condition that each ideal is finitely generated. We saw above that relaxing finiteness conditions, such as in direct products or formal power series, often causes properties true in the finite case to fall apart. In this sense, Noetherian rings serves as a bridge to keep the structure of ideals contained and well-behaved. They allow infinite rings through this bridge, but only when they satisfy the necessary properties to ensure the rings well-behaved nature. Therefore, in an analogous way to weak direct products and polynomial rings, Noetherian rings highlight the brittle nature of infinity, and provide a means to contain it.

Given a ring $R$ and a non-zero polynomial $f = a_nx^n+...+a_0 \in R[x]$ ($a_n \neq 0$), we say $\deg f = n$ (uniquely defined by Proposition \ref{prop2.28}). For technical reasons we define $\deg(0) \coloneqq -\infty$, along with the properties $(-\infty)< m$, $(-\infty)+m = -\infty = m+(-\infty)$, and $(-\infty)+(-\infty) = -\infty$ for every integer $m$ \cite{axlerlinear, hungerford}. The result below was originally proven by Hilbert \cite{hilbert}, and for further proofs, see \cite[Theorem VIII.4.9]{hungerford}, \cite[Theorem 6.42]{rotman}, and \cite[Theorem IV.4.1]{lang}.
\begin{theorem}[Hilbert's Basis Theorem]\label{thm2.31}
If $R$ is a commutative Noetherian Ring, then so is $R[x]$
\end{theorem}
\noindent \textbf{Proof}. We show every ideal $I \subset R[x]$ is finitely generated (Theorem \ref{thm2.30}). If $I = \{0\}$ the result is immediate, assume $I \neq \{0\}$. Let $f_1$ be a polynomial of least (non-negative) degree in $I$ and leading coefficient $a_1$. Recursively define $f_{m+1}$ ($m \in \mbb{Z}^{+}$) to be a polynomial of least degree in $I \setminus (f_1,...,f_m)$ with leading coefficient $a_{m+1}$, and if $I = (f_1,...,f_j)$ for some $j \in \mbb{Z}^{+}$, stop the process and we are done. Since $R$ is Noetherian, the ascending chain of ideals $$(a_1) \subset (a_1,a_2) \subset (a_1,a_2,a_3) \subset...$$ is stationary and thus $(a_1,a_2,...) = (a_1,...,a_k)$ for some $k \in \mbb{Z}^{+}$. We claim $I = (f_1,...,f_k)$. If not, consider $f_{k+1}$ and $a_{k+1} = \sum_{i = 1}^{k} c_ia_i$ for some $c_i \in R$. Since $\deg f_{k+1}$ is greater than or equal to $\deg f_i$ for all $i = 1,...,k$, $$g = \sum_{i = 1}^{k} c_i f_i x^{\deg f_{k+1}-\deg f_{i}} \in (f_1,...,f_k).$$
is a well-defined polynomial with the same degree and leading coefficient as $f_{k+1}$. If $f_{k+1}-g \in (f_1,...,f_k)$, then $f_{k+1} \in (f_1,...,f_k)$, a contradiction. Thus $f_{k+1}-g \in I \setminus (f_1,...,f_k)$ with degree strictly less than $f_{k+1}$, again a contradiction. Thus $I = (f_1,...,f_k)$, as desired. $\blacksquare$

Hilbert's basis theorem shows that if $R$ is Noetherian, ascending ideals in $R[x]$ are stationary. This prevents unbounded complexity in the structure of ideals within $R[x]$. $R[x]$ is inherently infinite due to it containing infinitely many monomials $r_k x^k \in R[x]$ ($k \in \mbb{N}$). For example, the ascending chain of ideals $$(r_1x) \subset (r_1x, r_2x^2) \subset (r_1x, r_2 x^2, r_3 x^3) \subset...$$
at first glance appears to grow indefinitely, for each step introduces a new generator. However, in the case when $R$ is Noetherian, such a chain stabilizes to some $(r_1 x,..., r_k x^k)$. This illustrates a fundamental connection between infinity and finiteness in algebra, that while $R[x]$ seems to allow infinitely growing chains due to its infinitely many monomials, the presence of a Noetherian ring $R$ ensures such complexity is ultimately governed by finitely generated ideals. For further proofs and discussion of the following theorem, see \cite[Theorem 3.3]{matsumura}, \cite[Proposition VIII.4.10]{hungerford}, and \cite[Theorem 7.5]{macdonald}.
\begin{theorem}\label{thm2.32}
If $R$ is a commutative Noetherian ring, then so is $R[[x]]$.
\end{theorem}
\noindent \textbf{Proof}. Let $I$ be an ideal in $R[[x]]$. We show $I$ is finitely generated (Theorem \ref{thm2.30}). Let $B = R[x]$ and for each $r \in \mbb{N}$ define $I(r) = \{a_r \in R|f = a_rx^{r}+a_{r+1}x^{r+1}+... \in I \cap x^{r}B\}$. Since $I$ is an ideal, so is $I(r)$, and since for each $f \in I \cap x^rB$, $xf \in I \cap x^{r+1}B$, we have $$I(0) \subset I(1) \subset I(2) \subset...$$
Since $R$ is Noetherian, the chain is stationary and there exists $s \in \mbb{N}$ such that $I(s) = I(s+1)=...$. Moreover, each $I(r)$ is finitely generated, so let $I(r) = (a_{r_1},...,a_{r_{k_r}})$ for each $r \in \mbb{N}$ and $g_{r_{i}} \in I \cap x^{r}B$ a polynomial with $a_{r_{i}}$ the coefficient of $x^r$ for each $1 \leq i \leq k_r$. Given $f \in I$, there exists $g_0 = \sum_{i = 1}^{k_0} c^{0}_{i}g_{0_{i}}$ ($c^{0}_{i} \in R$) such that $f-g_0 \in I \cap xB$. Similarly, there exists $g_1 =  \sum_{i = 1}^{k_1} c^{1}_{i}g_{1_{i}}$ ($c^{1}_{i} \in R$) such that $f-g_0-g_1 \in I \cap x^2B$. Repeating this process, we see $$f-g_0-g_1- \cdots -g_s \in I \cap x^{s+1}B.$$
Since $I(s) = I(s+1)$, we can write $g_{s+1} = \sum_{i = 1}^{k_s} c^{s+1}_{i}xg_{s_{i}}$ ($c^{s+1}_{s_{i}} \in R$) such that $$f-g_0-g_1- \cdots-g_{s+1} \in I \cap x^{s+2}B.$$
Proceed this way to obtain $g_{s+2},...$ ($\text{AC}_{\omega}$). For each $r < s$, each $g_r$ is a linear combination of $g_{r_{i}}$ ($1 \leq i \leq k_{r}$) with coefficients in $R$, and for $r \geq s$, we have $g_{r} = \sum_{i = 1}^{k_s} c^{r}_{i} x^{r-s} g_{s_{i}}$. Then let $$h_{i} = \sum_{r = s}^{\infty} c^{r}_{i} x^{r-s} \in B.$$ and $$f = g_0+ \cdots+g_{s-1}+ \sum_{i = 1}^{k_s} h_{i} g_{s_{i}}.$$
Thus $I$ is finitely generated and $B = R[[x]]$ is Noetherian, as desired. $\blacksquare$

Theorem \ref{thm2.32} shows that in the presence of a Noetherian ring $R$, $R[[x]]$ is also Noetherian. The structure $R[[x]]$ is inherently infinite due to its elements being the sum of infinitely many monomials $a_i x^i$. Despite this infinite structure, each ideal of elements of $R[[x]]$ is finitely generated. As with Theorem \ref{thm2.31}, this ensures the infinite complexity of the structure of ideals in $R[[x]]$ is ultimately governed by finitely generated ideals and establishes a connection between infinity and finiteness in ring theory.

The study of $R$-modules over rings illustrates more structural connections and differences between the finite and the infinite in the study of algebra. Naturally, we begin with the definition of an $R$-module.
\begin{definition}\label{def2.33}
Let $R$ be a ring. A left $R$-module is an abelian group $(A, +)$ together with a function $R \times A \to A$ called an action (whose image of $(r,a)$ is denoted $ra$) such that for all $r,s \in R$ and  $a,b \in A$ \\
(i) $r(a+b) = ra+rb$ \\
(ii) $(r+s)a = ra+sa$ \\
(iii) $r(sa) = (rs)a$ \\
(iv) $1_{R} a = a$ for all $a \in A$ \\
If $R$ is a field, then an $R$-module is called a vector space. Elements of a vector space are called vectors. A right $R$-module is defined analogously via a function $A \times R \to A$ denoted $(a,r) \mapsto ar$.
\end{definition}
\begin{example}\label{ex18}
Any abelian group is immediately a $\mbb{Z}$-module by definition. The set $$\mbb{R}^{n} = \left\{(a_1,...,a_n) \;\middle|\; a_1,...,a_n \in \mbb{R}\right\}$$ is a vector space over the real numbers $\mbb{R}$ with addition coordinate wise and $r(a_1,...,a_n) = (ra_1,...,ra_n)$ ($r \in \mbb{R}$). Any ring $R$ can be trivially made into a module over itself with an action of the multiplication in the ring, and likewise for any ideal in of $R$.
\end{example}
In the case when rings are defined without $1$, a module is defined only satisfying (i)-(iii), and rings with $1$ satisfying either satisfying (iv) or those satisfying (iv) called unitary $R$-modules (see \cite[Section 10.1]{dummit} and \cite[Section IV.1]{hungerford}). The study of vector spaces is called \textit{linear algebra}. Most definitions define vector spaces over fields \cite{roman, axlerlinear, dummit, lang}, but since much of the theory of linear algebra remains the same over non-commutative division rings, some authors define vector spaces over general division rings \cite{hungerford}. We will only consider left $R$-modules unless stated otherwise, and hence simply call them $R$-modules. Most results about left $R$-modules hold, mutatis mutandis, for right $R$-modules \cite[Section IV.1]{hungerford}.

If $R$ is commutative, every left $R$-module can be made into a right $R$-module by defining $ar = ra$ for $r \in R, a \in A$ (commutativity is needed to satisfy (iii)), and thus left and right modules coincide. We now prove some basic properties of modules.
\begin{proposition}\label{prop2.34}
Let $R$ be a ring and $A$ and $R$-module with $0_{A}$ additive identity on $A$ and $0_{R}$ on $R$. Then for all $r \in R, a \in A, n \in \mbb{Z}$ \\
(i) $r0_{A} = 0_{A}$ and $0_{R}a = 0_{A};$ \\
(ii) $(-r)a = -(ra) = r(-a)$ and $n(ra) = r(na)$.
\end{proposition}
\noindent \textbf{Proof}. [(i)] Notice $r0_{A} = r(0_{A}+0_{A}) = r0_{A}+r0_{A}$ and $r0_{A} = 0_{A}$. Similarly, $0_{R}a = (0_{R}+0_{R})a$ and $0_{R}a = 0_{A}$. [(ii)] Notice $(-r)a+ra = (-r+r)a = 0_{A}$ so $(-r)a = -(ra)$ and $r(-a)+ra = r(-a+a) = 0_{A}$ so $r(-a) = -ra$. Similarly, if $n \in \mbb{N}$, then $n(ra) = r(na)$ by Definition \ref{def2.33}(ii) and [(i)]. If $n < 0$, then $n(ra) = -((-n)(ra)) = -(r((-n)a)) = -(r(-na)) = ra$, as desired. $\blacksquare$

For modules, homomorphisms are defined slightly differently as 
\begin{definition}\label{def2.34}
Let $A$ and $B$ be modules over a ring $R$. A function $f:A \to B$ is an $R$-module homomorphism if $f$ is a group homomorphism such that for all $a \in A, r \in R$, $$f(ra) = rf(a).$$
If $f$ is injective [resp. surjective], then $f$ is called a monomorphism [resp. epimorphism] of modules. If $R$ is a field, $f$ is called a linear map.
\end{definition}
\begin{example}\label{ex19}
The map $f: \mbb{Z} \to \mbb{Z}$ by $n \mapsto 2n$ is a monomorphism of abelian groups, and hence of $\mbb{Z}$ modules. The map $T:\mbb{R}^{2} \to \mbb{R}$ by $(x,y) \mapsto x+y$ is a vector space homomorphism (over $\mbb{R}$), and more generally the map $f:A \oplus A \to A$ by $(a,b) \mapsto a+b$ with $A$ an $R$-module is an $R$-module homomorphism.
\end{example}
A \textit{submodule} of an $R$-module $A$ is a subgroup $B$ of $A$ such that $rb \in B$ for all $r \in R, b \in B$. A submodule of a vector space is called a \textit{subspace}. It is easy to verify that submodules of an $R$-module are themselves modules. Given a subset $X \subset A$, the module \textit{generated} by $X$ is defined as $RX = \{\sum_{i = 1}^{s} r_ia_i \mid s \in \mbb{Z}^{+}; a_i \in X; r_i \in R\}$, with the empty set generating the zero module. For further proofs and discussion of the following theorem, see \cite[Proposition 10.3]{dummit} and  \cite[Theorem IV.1.6]{hungerford}.
\begin{theorem}\label{thm2.36}
Let $R$ be a ring and $B$ a submodule of an $R$-module $A$. Then the quotient group $A/B$ is an $R$-module with action $r(a+B) = ra+B$ for all $r \in R, a \in A$. The map $\pi:A \to A/B$ by $a \mapsto a+B$ is an $R$-module epimorphism with kernel $B$.
\end{theorem}
\noindent \textbf{Proof}. Since $A$ is abelian, $B$ is a normal subgroup and $A/B$ is a well-defined, abelian group (Theorem \ref{thm2.4}). If $a+B = a'+B$, then $a-a' \in B$ so $ra-ra' = r(a-a') \in B$ for all $r \in R$ since $B$ is a submodule.Thus $(ra-ra')+B = B$ and $ra+B = ra'+B$, so the action is well-defined. We then see $r((a+B)+(b+B)) = r(a+b)+B = (ra+B)+(rb+B) = r(a+B)+r(b+B)$ and similarly $(r+s)(a+B) = r(a+B)+s(a+B)$ for all $r,s \in R, a,b \in A$. Moreover, $r(s(a+B)) = r(sa+B) = (rs)a+B = (rs)(a+B)$ and $1_{R}(a+B) = 1_{R}a+B = a+B$, so $A/B$ is an $R$-module. Since $\pi(ra) = ra+B = r(a+B) = r\pi(a)$ and $\pi$ is a group epimorphism with kernel $B$, $\pi$ is an $R$-module epimorphism with kernel $B$, as desired. $\blacksquare$

In light of Theorem \ref{thm2.36}, it is easy to see the analogue of the First Isomorphism Theorem (Theorem \ref{thm2.6}(i)) for modules. A subset $X$ of an $R$-module $A$ is said to be \textit{linearly independent} if for distinct $x_1,...,x_n \in X$ and $r_i \in R$ ($1 \leq i \leq n$) \begin{align}\label{l46}
r_1x_1+r_2x_2+...+r_nx_n = 0 \implies r_i = 0 \hspace{2mm} \forall i \in \{1,...,n\}.
\end{align}
A set that is not linearly independent is called \textit{linearly dependent}. If $A$ is generated by $X$, then $X$ is said to \textit{span} $A$. More generally, the set of all \textit{linear combinations} $\sum_{i = 1}^{n} r_ix_i$ ($r_i \in R, x_i \in X, n \in \mbb{Z}^{+}$) of $X$ is called the span of $X$. If $X$ spans $A$ and is linearly independent, $X$ is called a \textit{basis} of $A$. Notice the empty set is (vacuously) linearly independent and by definition is a basis of the zero module. A vector space $V$ over a field $F$ is called \textit{finite-dimensional} if there exists a finite spanning set of vectors $\{v_1,...,v_n\} \subset V$. For further proofs of the theorem below, see \cite[Theorem 3.78]{rotman} and \cite[Theorem 2.21]{axlerlinear}.
\begin{theorem}\label{thm2.37}
Every finite dimensional vector space $V$ over a field $F$ has a basis.
\end{theorem}
\noindent \textbf{Proof} Suppose $B = \{v_1,...,v_n\}$ (ordered) spans $V$. We will reduce $B$ to form a basis of $V$. If $v_1 = 0$, remove $v_1$ from $B$, otherwise leave $B$ unchanged. Repeating the process, if $j > 1$ and $v_j$ is in $\text{span}(v_1,...,v_{j-1})$, remove $v_j$ from $B$. If not, leave $v_j$ in $B$. Stopping the process after step $n$, we get a new set $B$, say $B = \{w_1,...,w_m\}$. A moments thought shows $B$ still spans $V$ since the only vectors removed were in the span of the previous vectors. If $B$ is linearly dependent, suppose $\sum_{i = 1}^{m} r_i w_i = 0 $ ($r_i \in F$) with some $r_i \neq 0$. Let $j$ be the largest element of $\{1,...,m\}$ such that $r_j \neq 0$. Then $\sum_{i = 1}^{j} r_i w_i = 0$ and (Theorem \ref{thm2.15}(ii)) \begin{align}\label{l47}
w_j = -r_j^{-1}r_1w_1- \cdots-r_j^{-1} r_{j-1} w_{j-1}
\end{align}
which contradicts the process above. Thus $B$ is a basis of $V$, as desired. $\blacksquare$

It is a well-known result in linear algebra that every vector space has a basis. Theorem \ref{thm2.37} shows that this is true for finite dimensional vector spaces, but the general case is in fact equivalent to Zorn's Lemma \cite{grabczewski}, and hence the axiom of choice. For further proofs of the following lemma, see \cite[Lemma IV.2.3]{hungerford}; for the theorem, see \cite[Theorem III.5.1]{lang} and \cite[Theorem IV.2.4]{hungerford}.
\begin{lemma}\label{lem2.38}
Every maximal linearly independent subset $X$ of a vector space $V$  over a field $F$ is a basis of $V$.
\end{lemma}
\noindent \textbf{Proof}. Let $W$ be the subspace of $V$ spanned by $X$. Since $X$ is linearly independent, $X$ is a basis of $W$. If $W = V$, we are done. If not, there exists $a \in V$ such that $a \notin W$ and $a \neq 0$. Consider $X \cup \{a\}$. If $ra+\sum_{i = 1}^{n} r_i x_i = 0$ ($r, r_i \in F, x_i \in X$) with $r \neq 0$, then $a = -r^{-1}r_1x_1- \cdots -r^{-1}r_nx_n \in W$, a contradiction. Hence $r = 0$ and since $X$ is linearly independent, $r_i = 0$ for all $1 \leq i \leq n$. Thus $X \cup \{a\}$ is linearly independent, contradicting the maximality of $X$. Therefore $W = V$ and $X$ is a basis of $V$, as desired. $\blacksquare$
\begin{theorem}\label{thm2.39}
Every vector space $V$ over a field $F$ has a basis. More generally, every linearly independent subset of $V$ is contained in a basis of $V$.
\end{theorem}
\noindent \textbf{Proof}. Since the empty set is (vacuously) a linearly independent subset of $V$, the first statement immediately follows from the second. Thus, let $X$ be a linearly independent subset of $V$ and $\mathcal{S}$ the set of all linearly independent subsets of $V$ that contain $X$. Since $X \in \mathcal{S}$, $\mathcal{S} \neq \emptyset$. Partially order $\mathcal{S}$ by set-theoretic inclusion and let $\{C_i \mid i \in \textit{I}\}$ be a chain in $\mathcal{S}$. Consider $C = \bigcup_{i \in \textit{I}} C_i$. If $x_1,...,c_n \in C$, then it follows that $x_1,..,x_n \in C_i$ for some $i \in \textit{I}$ since $\{C_i \mid i \in \textit{I}\}$ is a chain. Since $C_i$ is linearly independent, so is $C$. Thus by Zorn's Lemma $\mathcal{S}$ has a maximal element $M$, which is necessarily a maximal linearly independent subset of $V$, and hence a basis of $V$ by Lemma \ref{lem2.38}, as desired. $\blacksquare$
\begin{example}\label{ex20}
Unfortunately, unlike with vector spaces (Theorem \ref{thm2.39}), not every $R$-module has a basis. For example, take $\mbb{Q}$ over $\mbb{Z}$. Since $\mbb{Q}$ is not finitely generated (Example \ref{ex16}), any basis must be infinite. But given $\frac ab, \frac cd \in \mbb{Q} \setminus \{0\}$, oberserve $bc (\frac ab)-ad (\frac cd) = 0$ with $bc, ad \neq 0$, showing linear dependence. Thus $\mbb{Q}$ cannot have a basis.
\end{example}
The difference between the proofs of the existence of bases for finite and (potentially) infinite dimensional vector spaces now becomes evident. The finite dimensional case relies purely on finite cardinalities and finite processes. In contrast, the introduction of infinity changes the statement to require, and in fact be equivalent to, the axiom of choice. Much of the theory of linear algebra needed for applied mathematics can be restricted to finite dimensional vector spaces. After all, as Axler writes, ``Linear algebra is the study of linear maps on \textit{finite-dimensional} vector spaces" \cite[Section 1.A]{axlerlinear}. That said, for the purposes of pure mathematics, and even some applied mathematics, the statement that \textit{any} vector space has a basis is crucial, for which we must assume AC, inherently related to infinity by choice functions on infinite sets, and thus (as in Theorem \ref{thm2.37} and Theorem \ref{thm2.39}) displays disconnects between infinity and finiteness in the theory of $R$-modules. Thus, while the existence of bases in finite-dimensional cases follows from concrte, finite processes, the generalization is deeply intertwined within powerful set-theoretic assumptions. We now prove two more fundamental results from linear algebra in both finite and infinite dimensional cases. For further proofs and discussion of Lemma \ref{lem2.40}, see \cite[Theorem 1.10]{roman} and \cite[Theorem 2.23]{axlerlinear}; for Theorem \ref{thm2.41}, see \cite[Corollary 1.11]{roman} and \cite[Theorem 2.35]{axlerlinear}. For further discussion of Theorem \ref{thm2.42}, see \cite[Theorem III.5.2]{lang} and \cite[Theorem 1.12]{roman}.
\begin{lemma}\label{lem2.40}
Let $V$ be a finite-dimensional vector space over a field $F$. Then the cardinality of a linearly independent set is less than or equal to the cardinality of any spanning set.
\end{lemma}
\noindent \textbf{Proof}. Let $L$ be a linearly independent set and $B = \{w_1,...,w_n\}$ (in order) a finite spanning set. Suppose, by way of contradiction, $|L| > n$. Take $\{u_1,...,u_{n+1}\} \subset L$. Since $B$ spans $V$ and $u_1 \neq 0$, $u_1 = a_1w_1+...+a_nw_n$ ($a_i \in F$) with some $a_i \neq 0$. Thus it follows that $w_i \in \text{span}(u_1,w_1,...,w_{i-1},w_{i+1},...,w_n)$ and the set $\{u_1,w_1,...,w_{i-1},w_{i+1},...,w_n\}$ (in order) spans $V$. Replace $B$ with this set. Repeating the process, assume $j > 1$ and we completed the process up to $j-1$. Adjoin $u_j$ to $B$ right after $u_1,...,u_{j-1}$, creating a linearly dependent set. Then as in (\ref{l47}), some vector must be in the span of the previous vectors. It cannot be one of the $u_i's$ since $L$ is linearly independent, so it must be some $w_k$. Thus removing $w_k$ from $B$ still leaves a spanning set. After step $n$, we are left with $\{u_1,...,u_n\}$ as a spanning set, so $\{u_1,...,u_{n+1}\}$ is linearly dependent, a contradiction. Thus $|L| \leq n$. Now let $M$ be an arbitrary spanning set of $V$. If $|M|$ is infinite, clearly $|L| \leq |M|$. If $|M|$ is finite, the argument above shows $|L| \leq |M|$, as desired. $\blacksquare$
\begin{theorem}\label{thm2.41}
Any two bases of a finite-dimensional vector space $V$ have the same cardinality.
\end{theorem}
\noindent \textbf{Proof}. Let $B_1$ and $B_2$ be bases of $V$. Then viewing $B_1$ as linearly independent and $B_2$ as spanning, we have $|B_1| \leq |B_2|$ by Lemma \ref{lem2.40}. Similarly, viewing $B_1$ as spanning and $B_2$ as linearly independent, $|B_2| \leq |B_1|$. Thus $|B_1| = |B_2|$ by Schröder–Bernstein (Lemma \ref{lem1.4}), as desired, $\blacksquare$.
\begin{theorem}\label{thm2.42}
Every basis of a vector space $V$ over a field $F$ has the same cardinality. 
\end{theorem}
\noindent \textbf{Proof}.  If $V$ has a finite basis, it is finite-dimensional, and thus all bases have the same cardinality by Theorem \ref{thm2.41}. So assume every basis of $V$ is infinite. Let $X = \{x_i \mid i \in \textit{I}\}$ and $Y$ be bases of $V$. Then any vector $y \in Y$ can be written as a unique finite linear combination of vectors in $X$, with all coefficients non-zero, say $y = \sum_{i \in U_{y}} c_i x_i$ ($c_i \in F, x_i \in X$). Since $Y$ is a basis, $\bigcup_{y \in Y} U_{y} = \textit{I}$. Since $|U_{y}| < \bm{\al}_{0}$ for all $y \in Y$, by Lemma \ref{lem2.10}(ii) and Theorem \ref{thm1.9} $$|X| = |\textit{I}| = \abs{\bigcup_{y \in Y} U_{y}} \leq \bm{\al}_{0}|Y| = |Y|.$$
By symmetry, we then see $|Y| \leq |X|$ as well, thus by Schröder–Bernstein (Lemma \ref{lem1.4}) $|X| = |Y|$, as desired. $\blacksquare$

In light of Theorem \ref{thm2.42}, we call the cardinality of a basis of a vector space $V$ over a field $F$ its \textit{dimension}, which is then well-defined, and denote it by $\dim_{F} V$. Once again we observe the differences between the finite-dimensional case and the infinite-dimensional case. The finite case relies heavily on the assumption that there exists a finite spanning set, to then show any linearly independent set must be finite with cardinality less than or equal to the cardinality of any spanning set. In contrast, the general case shifts to only considering infinite bases, since if one basis was finite, then by definition $V$ must be finite-dimensional as well. The proof in that case relied heavily on properties of infinite cardinals, thus showing how the finite case extends to the infinite case by extending finite cardinals to infinite cardinals, and displaying a connection between them. Thus, the notion of dimension is highlighted in both theorems, but only through a deeper interaction with infinity. We now come to one of the most fundamental results from linear algebra, the \textit{rank-nullity theorem}, or as it's sometimes called, the \textit{fundamental theorem of linear maps} \cite{axlerlinear}. For further proofs of Theorem \ref{thm2.43}, see \cite[Theorem 3.22]{axlerlinear}; and for Theorem \ref{thm2.44}, see \cite[Theorem III.5.3]{lang} and \cite[Corollary IV.2.14]{hungerford}
\begin{theorem}\label{thm2.43}
Suppose $V,W$ are vector spaces over field $F$ and $T:V \to W$ is a linear map with $V$ finite-dimensional. Then $\Im T$ is finite-dimensional and $$\dim_{F} V = \dim_{F} \ker T+\dim_{F} \Im T$$
\end{theorem}
\noindent \textbf{Proof}. First off, since $T$ is a group homomorphism, $\ker T \leq V$ and $\Im T \leq W$ (as groups). Moreover, if $a \in \ker T$ and $b \in \Im T$, $T(ra) = rT(a) = r(0) = 0$ and $b = T(c)$ ($c \in V$) so $rb = rT(c) = T(rc) \in \Im T$ ($r \in R$). Thus $\ker T$ and $\Im T$ are subspaces of $V$ and $W$ respectively. If $\ker T$ were infinite dimensional, given any finite linearly independent subset of $\ker T$, there exists some $v \in \ker T$ not in its span, and adjoining $v$ still gives a linearly independent set. Repeating this process can give a linearly independent set of any finite cardinality, contradicting Lemma \ref{lem2.40}. Thus $\ker T$ is finite-dimensional, say $\{u_1,...,u_m\}$ (in order) is a basis (by Lemma \ref{lem2.40} a basis of a finite-dimensional space is finite). Since $V$ is finite-dimensional let $\{w_1,...,w_k\}$ be a basis of $V$ and \begin{align}\label{l48}
\{u_1,...,u_m,w_1,...,w_k\}
\end{align}
(in order) spans $V$. Applying the procedure in the proof of Theorem \ref{thm2.37} to reduce (\ref{l48}) to a basis, none of the $u_i's$ are removed since $\{u_1,...,u_m\}$ is linearly independent. Thus we are left with a basis $\{u_1,...,u_m,v_1,...,v_n\}$ of $V$ and $\dim_{F} V = m+n$. We claim $\{T(v_1),...,T(v_n)\}$ is a basis of $\Im T$. Given $v \in V$, write $v = \sum_{i = 1}^{m} a_i u_i+\sum_{j = 1}^{n} b_j v_j$ ($a_i,b_j \in F$). Apply $T$ to both sides for $T(v) = \sum_{j = 1}^{n} b_j T(v_j)$. Hence $\{T(v_1),...,T(v_n)\}$ spans $\Im T$. If $\sum_{i = 1}^{n} c_i T(v_i) = 0$ ($c_i \in F$), then $T(\sum_{i = 1}^{n} c_i v_i) = 0$ and $\sum_{i =1}^{n} c_iv_i \in \ker T$. Thus since $\{u_1,...,u_m\}$ spans $\ker T$, we have $$c_1v_1+...+c_nv_n = d_1u_1+...+d_mu_m$$
for some $d_i \in F$. Since $\{u_1,...,u_m,v_1,...,v_n\}$ is a basis, all the $c_i's$ and $d_j's$ must be zero, thus $\{T(v_1),...,T(v_n)\}$ is linearly independent and hence a basis. Therefore, $$\dim_{F} V = m+n = \dim_{F} \ker T+\dim_{F} \Im T,$$
as desired. $\blacksquare$
\begin{theorem}\label{thm2.44}
Suppose $V,W$ are vector spaces over $F$ and $T:V \to W$ a linear map. Then $$\dim_{F} V = \dim_{F} \ker T+\dim_{F} \Im T$$
\end{theorem}
\noindent \textbf{Proof}. Let $U =\{u_i \mid i \in \textit{I}\}$ be a basis of $\ker T$ (Theorem \ref{thm2.39}). Then by Theorem \ref{thm2.39}, $U$ is a subset of a basis $W$ of $V$. Consider $W \setminus U = \{v_j \mid j \in \textit{J}\}$. We claim $\{T(v_j) \mid j \in \textit{J}\}$ is a basis of $\Im T$. Given $v \in V$, we can write $v$ as $v = \sum_{i = 1}^{m} a_i u_i+\sum_{j = 1}^{n} b_j v_j$ ($a_i,b_j \in F, u_i \in U, v_j \in W \setminus U$). Thus $T(v) = \sum_{j = 1}^{n} b_j T(v_j)$ and $\{T(v_j) \mid j \in \textit{J}\}$ spans $\Im T$. Similarly, if $\sum_{i = 1}^{n} c_iT(v_i) = 0$ ($c_i \in F$), then $T(\sum_{i = 1}^{n} c_i v_i) = 0$ and $\sum_{i =1}^{n} c_iv_i \in \ker T$. Thus since $U$ spans $\ker T$, we have $$c_1v_1+...+c_nv_n = d_1u_1+...+d_mu_m$$
for some $d_i \in F, u_i \in U$. Then since $W$ is a basis, all the $c_i's$ and $d_j's$ must be zero, thus $\{T(v_j) \mid j \in \textit{J}\}$ is linearly independent and hence a basis. Therefore, $$\dim_{F} V = |U|+|W \setminus U| = \dim_{F} \ker T+\dim_{F} \Im T,$$
as desired. $\blacksquare$

Here we see a contrast from the first two examples above. The proofs in the finite-dimensional and the general case are almost identical. The only real difference is careful wording to ensure no assumption that $V$ or $\ker T$ are finite dimensional. The differences lies in the assumptions rather than the proof itself. Since $V$ may not be finite dimensional, and hence so may $\ker T$, to ensure $V$ or $\ker T$ have a basis relies on Zorn's Lemma, just like in Theorem \ref{thm2.39}. In contrast, the fact that $V$ or $\ker T$ have a basis when $V$ is finite-dimensional follows from Theorem \ref{thm2.37}, which requires only finite processes without need for the axiom of choice. In fact, much of the theory of finite-dimensional vector spaces, as seen above, can be viewed separately than the theory of infinite-dimensional vector spaces. In particular, the existence of bases for finite-dimensional spaces, equal cardinality of bases, and the rank-nullity theorem can all be proved independently of the axiom of choice or infinite cardinals, whereas the infinite-dimensional proofs either rely on the axiom of choice or infinite cardinals. The contrast of the two theorems then demonstrates how infinity can quietly permeate in results, and require foundational assumptions masked in the finite case that only become apparent in extensions to infinite cardinals.

More generally speaking, an $R$-module with a basis is called a \textit{free $R$-module}. Similarly to groups, the \textit{direct product} [resp. \textit{direct sum}] of a family of modules $\{A_i \mid i \in \textit{I}\}$, denoted $\prod_{i \in \textit{I}} A_i$ [resp. $\bigoplus_{i \in \textit{I}} A_{i}$], is the direct product [resp. direct sum] of $\{A_i \mid i \in \textit{I}\}$ as groups along with action $r\{a_i\} = \{ra_i\}$. The maps $\pi_{i}:\prod_{i \in \textit{I}} A_i \to A_{i}$ [resp. $\iota_{i}: A_{i} \to \prod_{i \in \textit{I}} A_i$] of group canonical projections [resp. injections] then naturally form $R$-module homomorphisms. It is then easy to see the counterpart of Proposition \ref{prop2.8}(ii) for modules, and thus we adopt the notation $A = \bigoplus_{i \in \textit{I}} A_{i}$ for the \textit{internal direct sum} of the family of modules $\{A_i \mid i \in \textit{I}\}$ (i.e., $A$ satisfies the counterpart of Proposition \ref{prop2.8}(ii) for modules). Just as free abelian groups ($\mbb{Z}$-modules) decompose into direct summands of $\mbb{Z}$ (Theorem \ref{thm2.9}), free $R$-modules decompose into direct summands of $R$. For further discussion of the following theorem, see \cite[Theorem IV.2.1]{hungerford}
\begin{theorem}\label{thm2.100}
Let $R$ be a ring. The following conditions on an $R$-module $F$ are equivalent; \\
(i) $F$ has a nonempty basis; \\
(ii) $F$ is the internal direct sum of a family of cyclic $R$-modules, each of which is isomorphic as a left $R$-module to $R$; \\
(iii) $F$ is $R$-module isomorphic to a direct sum of copies of the left $R$-module $R$.
\end{theorem}
\noindent \textbf{Proof}. [(i) $\Rightarrow$ (ii)] Let $X$ be a basis of $F$ and $x \in X$. Clearly $f:R \to Rx$ by $r \mapsto rx$ is an epimorphism and $\ker f = \{0\}$ by linear independence, whence $R \cong Rx$. Clearly $F = \langle \bigcup_{x \in X} Rx \rangle$ since $X$ is a basis. If $z \in Rx \cap \langle \bigcup_{y \neq x} Ry \rangle$, then $z = rx = \sum_{i = 1}^{n} r_i x_i$ ($r_i \in R, x_i \in X \setminus \{x\}$) and $rx-\sum_{i = 1}^{n} r_i x_i = 0$, so $r = r_1 = \cdots = r_n = 0$ and $z = 0$. Thus $F = \bigoplus_{x \in X} Rx$. 

[(ii) $\Rightarrow$ (iii)] Suppose $F = \bigoplus_{x \in X} \langle x \rangle$ with $\langle x \rangle \cong R$ for each $x \in X$. For each $x \in X$, let $\psi_{x}: \langle x \rangle \to R$ be an isomorphism (AC). Then the map $f:\bigoplus_{x \in X} \langle x \rangle \to \bigoplus_{x \in X} R$ by $\{a_x\}_{x \in X} \mapsto \{\psi_{x}(a_x)\}_{x \in X}$ is already seen to be a group isomorphism by Theorem \ref{thm2.7}, and moreover $f(r\{a_x\}) = f(\{ra_x\}) = \{\psi_{x}(ra_x)\} = \{r \psi_{x}(a_x)\} = r\{\psi_{x}(a_x)\} = rf(\{a_x\})$ (where $x \in X$ was dropped for ease of notation) so $f$ is a $R$-module isomorphism and $F \cong \bigoplus_{x \in X} \langle x \rangle \cong \bigoplus_{x \in X} R$ by the analogue of Proposition \ref{prop2.8}(ii) for modules.
 
[(iii) $\Rightarrow$ (i)] Suppose $F \cong \bigoplus R$ indexed by set $X$. For each $x \in X$, let $\theta_{x}$ be the element $\{r_i\}$ of $\sum R$ with $r_i = 0$ for $i \neq x$ and $r_x = 1$. Identically to Theorem \ref{thm2.9} (mutandis mutandis) we show $\{\theta_{x} \mid x\in X\}$ is a basis of $\bigoplus R$ and under the isomorphism $F$ obtains a non-empty basis, as desired. $\blacksquare$

Although Theorem \ref{thm2.100} closely parallels Theorem \ref{thm2.9}, there are some distinctions worth emphasizing. In the case when $R$ is finite, the presence of a free $R$-module alone no longer offers an immediate association to infinity. However, the structure of free $R$-modules still opens the door to infinite behavior. In any infinite ring, the connection is analogous to Theorem \ref{thm2.9}. Moreover, through any infinite basis, the free module is isomorphic to a direct sum of infinitely many cyclic modules isomorphic to $R$, again echoing Theorem \ref{thm2.9}. Thus, free $R$-modules, even when $R$ is finite, are deeply rooted in a relationship to infinite structures.

Theorem \ref{thm2.39} shows that every vector space over a field $F$ is a free $F$-module. A ring $R$ such that every basis of a free $R$-module has the same cardinality is said to have the \textit{invariant dimension property}. Theorem \ref{thm2.42} shows that all fields have the invariant dimension property. In fact, as one may verify, none of the proofs above required the commutativity under multiplication of a field, hence it also holds that every division ring has the invariant dimension property (for alternate proofs, see \cite[Section IV.2]{hungerford}). Theorem \ref{thm2.11} shows that the ring $\mbb{Z}$ of integers has the invariant dimension property. In fact, it can be proven that, just like the special cases of free-abelian groups and $R$-modules over fields, that every \textit{infinite} basis of a free $R$-module over any ring has the same cardinality (see \cite[Theorem IV.2.6]{hungerford}). Unfortunately, it is not true in general that linearly independent sets can be extended to a basis or spanning sets reduced to a basis with free $R$-modules, for example with free $\mbb{Z}$-modules (see \cite[Exercise II.1.2]{hungerford}).

We now examine modules more generally. We begin by defining an important type of $R$-module.
\begin{definition}\label{def1.100}
A triangle diagram of functions
$$\begin{tikzcd}
A \arrow[r, "f"] \arrow[d, "g"'] & B \\
C \arrow[ur, "h"']
\end{tikzcd}$$
is called commutative if $f = hg$. A square diagram of functions
$$\begin{tikzcd}
A \arrow[r, "\varphi"] \arrow [d, "\zeta"] & B \arrow[d, "\psi"] \\
 C \arrow[r, "\sigma"] & D
\end{tikzcd}$$
is called commutative if $\psi \varphi  = \sigma \zeta$. A function diagram is called commutative if every square and triangle in the diagram is commutative.
\end{definition}
\begin{definition}\label{def2.45}
Let $R$ be a ring. An $R$-module $P$ is called projective if for every diagram of $R$-module homomorphisms
$$\begin{tikzcd}
& P \arrow[d, "f"] \\
A \arrow[r, "g"] & B \arrow[r] & 0
\end{tikzcd}$$
with $g$ an epimorphism, there exists an $R$-module homomorphism $h:P \to A$ such that the diagram 
$$\begin{tikzcd}
& P \arrow[d, "f"] \arrow[dl, dashed, "h" sloped]  \\
A \arrow[r, "g"] & B \arrow[r] & 0
\end{tikzcd}$$
is commutative (i.e., $f = gh$).
\end{definition}
The bottom row (i.e., $A \xrightarrow{g} B \xrightarrow{} 0$) in Definition \ref{def2.45} is an example of an \textit{exact sequence}, which we now define.
\begin{definition}[\cite{hungerford}, Definition IV.1.16]\label{def2.46}
A pair of module homomorphisms $A \xrightarrow{f} B \xrightarrow{g} C$ is said to be exact at $B$ provided $\Im f = \ker g$. A finite sequence of module homomorphisms, $A_0 \xrightarrow{f_1} A_1 \xrightarrow{f_2} A_2 \xrightarrow{f_3} \cdots \xrightarrow{f_{n-1}} A_{n-1} \xrightarrow{f_n} A_n$ is exact provided $\Im f_i = \ker f_{i+1}$ for $i = 1,2,...,n-1$. An infinite sequence $\cdots \xrightarrow{f_{i-1}} A_{i-1} \xrightarrow{f_i} A_{i} \xrightarrow{f_{i+1}} A_{i+1} \xrightarrow{f_{i+2}} \cdots$ is exact provided $\Im f_{i} = \ker f_{i+1}$ for all $i \in \mbb{Z}$.
\end{definition}
The bottom row $A \xrightarrow{g} B \xrightarrow{} 0$ in the definition of a projective module is exact because $\Im g = B$ and the kernel of $B \to 0$ is all of $B$. For further proofs and discussion of the following theorem, see \cite[Theorem IV.3.2]{hungerford}. 
\begin{theorem}\label{thm2.47}
Every free module $F$ is projective
\end{theorem}
\noindent \textbf{Proof}. Let $M$ be a basis of $F$. Suppose $f: F \to B$ is a homomorphism and $g:A \to B$ an epimorphism (of modules). For each $m \in M$, the element $f(m) \in B$ has a non-empty pre-image in $A$ under $g$ since $g$ is surjective. Thus for each $m \in M$, choose $a_m \in A$ such that $g(a_m) = f(m)$ (AC). Define a homomorphism $h:F \to A$ by $\sum_{i = 1}^{n} r_i m_{i} \mapsto \sum_{i = 1}^{n} r_i a_{m_i}$ ($r_i \in R, m_i \in M$). Each $x \in F$ can be expressed uniquely as such sum $x =  \sum_{i = 1}^{n} r_i m_{i}$ (since $F$ is a basis) and $g(h(x)) = g\qty(\sum_{i = 1}^{n} r_i a_{m_i}) = \sum_{i = 1}^{n} r_i g(a_{m_i}) = \sum_{i = 1}^{n} r_if(m) = f(x).$ Thus $f = gh$ and the diagram $$\begin{tikzcd}
& F \arrow[d, "f"] \arrow[dl, dashed, "h" sloped]  \\
A \arrow[r, "g"] & B \arrow[r] & 0
\end{tikzcd}$$
is commutative, as desired. $\blacksquare$

Theorem \ref{thm2.47} shows every free $R$-module is projective, and one may be tempted to infer from this a deep relationship between free $R$-modules, projective modules, and infinity. However, the fact that every free $R$-module is an $R$-module is trivial, but doesn't imply that free modules and general modules share meaningful structural properties. Therefore, a deeper connection must be made to claim infinity is embedded in projective modules. 

In the case of a general ring $R$ and projective $R$-module $P$, we can say that there exists a submodule $K$ of $R$ such that there exists a free $R$-module $F$ with $F \cong K \oplus P$, which indeed highlights a much deeper connection between general free-modules, projective modules, and infinity. Even further, over specific classes of rings, such as principal ideal domains, we can show they are actually equivalent. This equivalence allows one to directly transfer properties of free modules—including those shaped by infinite bases or decompositions—to the realm of projective modules. Therefore, these two results establish how free modules, projective modules, and infinity relate in both general and specialized settings. We now build towards the two results. 

An exact sequences of the form $0 \xrightarrow{} A  \xrightarrow{f} B  \xrightarrow{g} C  \xrightarrow{} 0$ (i.e., $f$ is a monomorphism, g is an epimorphism) is called a \textit{short exact sequence}. Two short exact sequences $0 \xrightarrow{} A  \xrightarrow{f} B  \xrightarrow{g} C  \xrightarrow{} 0$ and $0 \xrightarrow{} A'  \xrightarrow{f'} B  \xrightarrow{g'} C'  \xrightarrow{} 0$ are \textit{isomorphic} if there exists a commutative diagram $$\begin{tikzcd}
	0 & A & B & C & 0 \\
	0 & {A'} & {B'} & {C'} & 0
	\arrow[from=1-1, to=1-2]
	\arrow[from=1-2, to=1-3, "f"]
	\arrow[from=1-2, to=2-2, "\alpha"]
	\arrow[from=1-3, to=1-4, "g"]
	\arrow[from=1-3, to=2-3, "\beta"]
	\arrow[from=1-4, to=1-5]
	\arrow[from=1-4, to=2-4, "\gamma"]
	\arrow[from=2-1, to=2-2]
	\arrow[from=2-2, to=2-3, "f'"]
	\arrow[from=2-3, to=2-4, "g'"]
	\arrow[from=2-4, to=2-5]
\end{tikzcd}$$
such that $\alpha, \beta$, and $\gamma$ are all isomorphisms. It is not hard to verify that isomorphisms of short exact sequences form an equivalence relation \cite[Exercise IV.1.14]{hungerford}. A short exact sequences  $0 \xrightarrow{} A_1  \xrightarrow{f} B  \xrightarrow{g} A_2  \xrightarrow{} 0$ isomorphic to $0 \xrightarrow{} A_1  \xrightarrow{\iota_1} A_1 \oplus A_2  \xrightarrow{\pi_2} A_2  \xrightarrow{} 0$ is called \textit{split exact}. For further proofs and discussion of the following theorem, see \cite[Lemma IV.1.17]{hungerford}. For Theorem \ref{thm2.50}, see \cite[Proposition III.3.2]{lang} and \cite[Theorem IV.1.18]{hungerford}. 
\begin{theorem}[Short Five Lemma]\label{thm2.48}
Let $R$ be a ring and $$\begin{tikzcd}
	0 & A & B & C & 0 \\
	0 & {A'} & {B'} & {C'} & 0
	\arrow[from=1-1, to=1-2]
	\arrow[from=1-2, to=1-3, "f"]
	\arrow[from=1-2, to=2-2, "\alpha"]
	\arrow[from=1-3, to=1-4, "g"]
	\arrow[from=1-3, to=2-3, "\beta"]
	\arrow[from=1-4, to=1-5]
	\arrow[from=1-4, to=2-4, "\gamma"]
	\arrow[from=2-1, to=2-2]
	\arrow[from=2-2, to=2-3, "f'"]
	\arrow[from=2-3, to=2-4, "g'"]
	\arrow[from=2-4, to=2-5]
\end{tikzcd}$$
a commutative diagram of $R$-modules and $R$-module homomorphisms such that each row is short exact. Then \\
(i) $\alpha, \gamma$ monomorphisms $\implies$ $\beta$ is a monomorphism;  \\
(ii)  $\alpha, \gamma$ epimorphism $\implies$ $\beta$ is an epimorphism; \\
(iii)  $\alpha, \gamma$ isomorphism $\implies$ $\beta$ is an isomorphism.
\end{theorem}
\noindent \textbf{Proof}. [(i)] Let $b \in B$ and suppose $\beta(b) = 0$. By commutativity $$\gamma g(b) = g' \beta(b) = g'(0) = 0.$$
Thus $g(b) = 0$ since $\gamma$ is a monomorphism. Since the top row is exact, $b \in \ker g = \Im f$ and $b = f(a)$ for some $a \in A$. By commutativity, $$f' \alpha(a) = \beta f(a) = \beta(b) = 0.$$
Since the bottom row is exact, $f'$ is a monomorphism and $\alpha(a) = 0$, so $a = 0$. Thus $b = f(a) = f(0) = 0$ and $\beta$ is a monomorphism (\ref{l42}).

[(ii)] Let $b' \in B'$. Then $g'(b') \in C'$ and since $\gamma$ is an epimorphism $g'(b') = \gamma(c)$ for some $c \in C$. By the exactness of the top row, $g$ is an epimorphism and $c = g(b)$ for some $b \in B$. By commutativity, $$g' \beta(b) = \gamma g(b) = \gamma(c) = g'(b').$$
Thus $g'[\beta(b)-b'] = 0$ and $\beta(b)-b' \in \ker g' = \Im f'$. Let $f'(a') = \beta(b)-b'$ ($a' \in A'$). Since $\alpha$ is an epimorphism, let $\alpha(a) = a'$ ($a \in A$). By commutativity, $\beta f(a) = f' \alpha(a) = f'(a') = \beta(b)-b'$. Hence $$\beta[b-f(a)] = \beta(b)-\beta f(a) = \beta(b)-(\beta(b)-b') = b'$$
and $\beta$ is an epimorphism. [(iii)] Since $\alpha$ and $\gamma$ are both monomorphisms and epimorphisms, so is $\beta$ by (i) and (ii) and hence $\beta$ is an isomorphism, as desired. $\blacksquare$
\begin{lemma}\label{lem2.49}
Let $R$ be a ring and $0 \xrightarrow{} A_1  \xrightarrow{f} B  \xrightarrow{g} A_2  \xrightarrow{} 0$ a short exact sequence of $R$-module homomorphisms. Then $0 \xrightarrow{} A_1  \xrightarrow{f} B  \xrightarrow{g} A_2  \xrightarrow{} 0$ is split exact if and only if there exists an isomorphism $\varphi: A_1 \oplus A_2 \to B$ such that the diagram  
$$\begin{tikzcd}
	0 & A_1 & A_1 \oplus A_2 & A_2 & 0 \\
	0 & {A_1} & {B} & {A_2} & 0
	\arrow[from=1-1, to=1-2]
	\arrow[from=1-2, to=1-3, "\iota_1"]
	\arrow[from=1-2, to=2-2, "1_{A_1}"]
	\arrow[from=1-3, to=1-4, "\pi_2"]
	\arrow[from=1-3, to=2-3, "\varphi"]
	\arrow[from=1-4, to=1-5]
	\arrow[from=1-4, to=2-4, "1_{A_2}"]
	\arrow[from=2-1, to=2-2]
	\arrow[from=2-2, to=2-3, "f"]
	\arrow[from=2-3, to=2-4, "g"]
	\arrow[from=2-4, to=2-5]
\end{tikzcd}$$
is commutative.
\end{lemma}
\noindent \textbf{Proof}. If the given sequence is split exact, by definition there exists a commutative diagram $$\begin{tikzcd}
	0 & A_1 & A_1 \oplus A_2 & A_2 & 0 \\
	0 & {A_1} & {B} & {A_2} & 0
	\arrow[from=1-1, to=1-2]
	\arrow[from=1-2, to=1-3, "\iota_1"]
	\arrow[from=1-2, to=2-2, "\zeta"]
	\arrow[from=1-3, to=1-4, "\pi_2"]
	\arrow[from=1-3, to=2-3, "\psi"]
	\arrow[from=1-4, to=1-5]
	\arrow[from=1-4, to=2-4, "\sigma"]
	\arrow[from=2-1, to=2-2]
	\arrow[from=2-2, to=2-3, "f"]
	\arrow[from=2-3, to=2-4, "g"]
	\arrow[from=2-4, to=2-5]
\end{tikzcd}$$
with $\zeta, \psi, \sigma$ isomorphims. Let $\varphi = \psi (\zeta^{-1} \oplus \sigma^{-1})$, where $\zeta^{-1} \oplus \sigma^{-1}:A_1 \oplus A_2 \to A_1 \oplus A_2$ is defined by $(a_1, a_2) \mapsto (\zeta^{-1}(a_1), \sigma^{-1}(a_2))$ ($a_i \in A_i$). By Theorem \ref{thm2.48}, it suffices to show $g \varphi = \pi_{2}$ and $\varphi \iota_1 = f$. Notice $g \varphi(a_1,a_2) = g(\psi(\zeta^{-1}(a_1), \sigma^{-1}(a_2))) = \sigma \pi_2(\psi(\zeta^{-1}(a_1), \sigma^{-1}(a_2))) = \sigma(\sigma^{-1}(a_2)) = a_2 = \pi_2(a_1,a_2)$ and $\varphi \iota_1(a_1) = \varphi(a_1,0) = \psi(\zeta^{-1}(a_1),0) = \psi \iota_{1}(\zeta^{-1}(a_1)) = f \zeta(\zeta^{-1}(a_1)) = f(a_1)$ ($a_i \in A_i$). Conversely, if such a commutative diagram exists with $\varphi$ isomorphism, then the sequences are immediately isomorphic by definition, as desired. $\blacksquare$
\begin{theorem}\label{thm2.50}
Let $R$ be a ring and $0 \xrightarrow{} A_1  \xrightarrow{f} B  \xrightarrow{g} A_2  \xrightarrow{} 0$ a short exact sequence of $R$-module homomorphisms. Then the following conditions are equivalent. \\
(i) There is an $R$-module homomorphism $h:A_2 \to B$ with $gh = 1_{A_{2}};$ \\
(ii) There is an $R$-module homomorphism $k:B \to A_1$ with $kf = 1_{A_1};$ \\
(iii) The given sequence is split exact.
\end{theorem}
\noindent \textbf{Proof}. [(i) $\Rightarrow$ (iii)] It is a routine check that $\varphi:A_1 \oplus A_2 \to B$ given by $(a_1,a_2) \mapsto f(a_1)+h(a_2)$ is a homomorphism. Consider the diagram 
$$\begin{tikzcd}
	0 & A_1 & A_1 \oplus A_2 & A_2 & 0 \\
	0 & {A_1} & {B} & {A_2} & 0
	\arrow[from=1-1, to=1-2]
	\arrow[from=1-2, to=1-3, "\iota_1"]
	\arrow[from=1-2, to=2-2, "1_{A_1}"]
	\arrow[from=1-3, to=1-4, "\pi_2"]
	\arrow[from=1-3, to=2-3, "\varphi"]
	\arrow[from=1-4, to=1-5]
	\arrow[from=1-4, to=2-4, "1_{A_2}"]
	\arrow[from=2-1, to=2-2]
	\arrow[from=2-2, to=2-3, "f"]
	\arrow[from=2-3, to=2-4, "g"]
	\arrow[from=2-4, to=2-5].
\end{tikzcd}$$
Given $a_1 \in A_1$, we have $\varphi \iota_1(a_1) = \varphi(a_1,0) = f(a_1) = f1_{A_1}(a_1)$ and $\varphi \iota_1 = f1_{A_1}$. For each $(a_1,a_2) \in A_1 \oplus A_2$, since the bottom row is exact, we have $1_{A_2} \pi_2(a_1,a_2) = a_2 = gh(a_2)+gf(a_1) = g(f(a_1)+h(a_2)) = g\varphi(a_1,a_2)$. Thus the diagram is commutative. Then by Theorem \ref{thm2.48}, $\varphi$ is an isomorphism and the given sequence is split exact.

[(ii) $\Rightarrow$ (iii)] Consider the diagram $$\begin{tikzcd}
	0 & A_1 & B & A_2 & 0 \\
	0 & {A_1} & {A_1 \oplus A_2} & {A_2} & 0
	\arrow[from=1-1, to=1-2]
	\arrow[from=1-2, to=1-3, "f"]
	\arrow[from=1-2, to=2-2, "1_{A_1}"]
	\arrow[from=1-3, to=1-4, "g"]
	\arrow[from=1-3, to=2-3, "\psi"]
	\arrow[from=1-4, to=1-5]
	\arrow[from=1-4, to=2-4, "1_{A_2}"]
	\arrow[from=2-1, to=2-2]
	\arrow[from=2-2, to=2-3, "\iota_1"]
	\arrow[from=2-3, to=2-4, "\pi_2"]
	\arrow[from=2-4, to=2-5]
\end{tikzcd}$$
where $\psi$ is the homomorphism defined by $b \mapsto (k(b),g(b))$. Notice $\pi_2 \psi(b) = g(b) = 1_{A_2} g$ and $\psi f(a_1) = (kf(a_1), gf(a_1)) = (a_1,0) = \iota_1 1_{A_1}(a_1)$ for every $a_1 \in A_1, b \in B$, whence the diagram is commutative. Then by Theorem \ref{thm2.48}, $\psi$ is an isomorphism and the given sequence is split exact by Lemma \ref{lem2.49}.

[(iii) $\Rightarrow$ (i), (ii)] By Lemma \ref{lem2.49} consider the commutative diagram with exact rows and $\varphi$ an isomorphism. 
$$\begin{tikzcd}
	0 & A_1 & A_1 \oplus A_2 & A_2 & 0 \\
	0 & {A_1} & {B} & {A_2} & 0
	\arrow[from=1-1, to=1-2]
	\arrow[from=1-2, to=1-3, "\iota_1"]
	\arrow[from=1-2, to=2-2, "1_{A_1}"]
	\arrow[from=1-3, to=1-4, "\pi_2"]
	\arrow[from=1-3, to=2-3, "\varphi"]
	\arrow[from=1-4, to=1-5]
	\arrow[from=1-4, to=2-4, "1_{A_2}"]
	\arrow[from=2-1, to=2-2]
	\arrow[from=2-2, to=2-3, "f"]
	\arrow[from=2-3, to=2-4, "g"]
	\arrow[from=2-4, to=2-5].
\end{tikzcd}$$
Define $h = \varphi \iota_2$ and $k = \pi_{1} \varphi^{-1}$. By the commutativity of the diagram, $kf = (\pi_1 \varphi^{-1})(f 1_{A_1}) = (\pi_1 \varphi^{-1})(\varphi \iota_1) = \pi_1 \iota_1 = 1_{A_1}$ and similarly $gh = g(\varphi \iota_2) = (1_{A_2} \pi_2)\iota_2 = 1_{A_2}$. Thus (i) and (ii) follow, as desired. $\blacksquare$

In algebra, it is often important to characterize structures as we did with free-abelian groups (Theorem \ref{thm2.9}). We now characterize projective modules. For further proofs and discussion of Theorem \ref{thm2.52}, see \cite[Section III.4]{lang}, \cite[Proposition 7.54]{rotman}, and  \cite[Theorem IV.3.4]{hungerford}.
\begin{lemma}\label{lem2.51}
Let $R$ be a ring. Every $R$-module $A$ is the homomorphic image of a free $R$-module.
\end{lemma}
\noindent \textbf{Proof}. Let $X$ be a set of generators of $A$ and consider $\bigoplus_{x \in X} R$. By the proof of Theorem \ref{thm2.100}, $\bigoplus_{x \in X} R$ is free with basis $\{\theta_{x} \mid x \in X\}$. Then the map $f: \bigoplus_{x \in X} R \to A$ by $\sum_{i = 1}^{k} r_i \theta_{x_{i}} \mapsto \sum_{i = 1}^{k} r_i x_i$ is an epimorphism, as desired. $\blacksquare$

Lemma \ref{lem2.51}, in an analogous way to Proposition \ref{prop2.12}, shows a more general relationship between infinity and $R$-modules. Proposition \ref{prop2.12} addresses the specific case of $\mbb{Z}$-modules (abelian groups), and Lemma \ref{lem2.51} generalizes the result to arbitrary rings, showing each $R$-module is the homomorphic image of a free $R$-module (with a potential infinite basis). Although this connection is not deep since it only says that given an $R$-module $A$, each element $a \in A$ can be written as a finite linear combination of images of basis elements under an epimorphism $f$, it still highlights how infinite bases and images of free $R$-modules arise naturally in module theory.
\begin{theorem}\label{thm2.52}
Let $R$ be a ring. Then the following conditions on an $R$-module $P$ are equivalent. \\
(i) $P$ is projective; \\
(ii) every short exact sequence $0 \xrightarrow{} A \xrightarrow{f} B \xrightarrow{g} P \xrightarrow{} 0$ is split exact; \\
(iii) there exists a free module $F$ and an $R$-module $K$ such that $F \cong K \oplus P$.
\end{theorem}
\noindent \textbf{Proof}. [(i) $\Rightarrow$ (ii)] Consider the diagram
$$\begin{tikzcd}
	& P \\
	B & P & 0
	\arrow[from=1-2, to=2-2, "1_{P}"]
	\arrow[from=2-1, to=2-2, "g"]
	\arrow[from=2-2, to=2-3]
\end{tikzcd}$$
with the bottom row exact. Since $P$ is projective there is an $R$-module homomorphism $h:P \to B$ such that $gh = 1_{P}$. Thus the short exact sequence $0 \xrightarrow{} A \xrightarrow{f} B \xrightarrow{g} P \xrightarrow{} 0$ is split exact by Theorem \ref{thm2.50}.

[(ii) $\Rightarrow$ (iii)] By Lemma \ref{lem2.51}, there exists a free module $F$ and an epimorphism $g:F \to P$. If $K = \ker g$, then $0  \xrightarrow{} K \xrightarrow{\subset} F \xrightarrow{g} 0$ is exact and by hypothesis $F \cong K \oplus P$.

[(iii) $\Rightarrow$ (i)] Let $\psi$ be the composition $F \xrightarrow{\alpha} K \oplus P \xrightarrow{\pi_2} P$ and $\varphi$ the composition $P \xrightarrow{\iota_2} K \oplus P \xrightarrow{\alpha^{-1}} F$, where $\alpha:F \cong K \oplus P$. Given a diagram of $R$-module homomorphisms  
$$\begin{tikzcd}
& P \arrow[d, "f"] \\
A \arrow[r, "g"] & B \arrow[r] & 0
\end{tikzcd}$$
with bottom row exact, consider the diagram 
$$\begin{tikzcd}
& F \arrow[d, "f \psi"] \\
A \arrow[r, "g"] & B \arrow[r] & 0
\end{tikzcd}$$
Since $F$ is projective by Theorem \ref{thm2.47}, there is an $R$-module homomorphism $u:F \to A$ such that $gu = f \psi$. Let $h = u \varphi:P \to A$. Then $gh = gu \varphi = (f \psi)\varphi = f(\psi \varphi) = f 1_{P} = f$. Thus $P$ is projective, as desired. $\blacksquare$

As discussed above, Theorem \ref{thm2.52} highlights a connection between general projective $R$-modules, general free $R$-modules, and infinity. Additionally, in analyzing the proof, the implication [(ii) $\Rightarrow$ (iii)] relies heavily on Lemma \ref{lem2.51}, another result showing how more general modules relate to infinity. While it holds true that every free $R$-module is projective (Theorem \ref{thm2.47}), Example \ref{ex21} shows that the converse isn't always true. Theorem \ref{thm2.52} characterizes projectivity in terms of exact sequences and decomposition, along with establishing a strong link between projective modules and free modules, even in the more general case where they are not equivalent. This connection highlights how infinity manifests not just in terms of cardinality or dimension, but in algebraic structures and decompositions themselves.
\begin{example}\label{ex21}
Consider the ring $R = \mbb{Z}/6\mbb{Z}$. Clearly $R$ is a free module over itself. $\mbb{Z}/2\mbb{Z}$ and $\mbb{Z}/3\mbb{Z}$ are $R$-modules under action $r(a+n\mbb{Z}) = ra+n\mbb{Z}$ ($r \in R, n = 2,3$), and under isomorphism $a+6\mbb{Z} \mapsto (a+3\mbb{Z},a+2\mbb{Z})$, we have $R \cong \mbb{Z}/3\mbb{Z} \oplus  \mbb{Z}/2\mbb{Z}$. Thus by Theorem \ref{thm2.52}, $\mbb{Z}/2\mbb{Z}$ is projective, however since it is of cardinality $2$, it is clearly not free. Thus not every free $R$-module is projective in general, but as we will see, in specific cases, it is.
\end{example}
In Theorem \ref{thm2.47} we showed that every free $R$-module is projective. As discussed above, over principal ideal domains, the converse also holds, as we now show. Lemma \ref{lem2.53}(i) is in \cite[Exercise IV.1.15]{hungerford}; for further proofs and discussion of (ii), see \cite[Corollary 7.55]{rotman}; for (iii) see \cite[Theorem 7.9]{rotman}. For further discussion of Theorem \ref{thm2.54}, see \cite[Theorem IV.6.1]{hungerford} and \cite[Theorem 9.8]{rotman}. 
\begin{lemma}\label{lem2.53}
Let $R$ be a ring with $A$ and $B$ $R$-modules. \\
(i) If $f:A \to B$ and $g:B \to A$ are $R$-module homomorphisms such that $gf = 1_{A}$, then $B = \Im f \oplus \ker g;$ \\
(ii) If $A$ is a submodule of $B$ and $B/A$ projective, then there exists a submodule $C$ of $B$ such that $C \cong B/A$ and $B = A \oplus C;$ \\
(iii) (Second Isomorphism Theorem) If $A$ and $B$ are submodules of a module $C$, the map $\varphi:A \to (A+B)/B$ by $a \mapsto a+B$ is an epimorphism of modules. Consequently, $A/(A \cap B) \cong (A+B)/B$ ($A+B = \{a+b \mid a \in A, b \in B\}$). 
\end{lemma}
\noindent \textbf{Proof}. [(i)] Given $b \in B$, write $b = f(g(b))+(b-f(g(b))$. Clearly $f(g(b)) \in \Im f$ and notice $g(b-f(g(b))) = g(b)-gf(g(b)) = g(b)-g(b) = 0$ thus $b-f(g(b)) \in \ker g$ and $B = \Im f+\ker g$. If $c = f(a) \in \Im f$ and $g(c) = 0$, then $g(f(a)) = a = 0$ and $c = f(0) = 0$ so $\Im f \cap \ker g = \{0\}$ and $B = \Im f \oplus \ker g$. [(ii)] Since the short exact sequence $0 \xrightarrow{} A \xrightarrow{\subset} B \xrightarrow{\pi} B/A \xrightarrow{} 0$ is split exact, there exists homomorphism $\psi:B/A \to B$ such that $\pi \psi = 1_{B/A}$ (Theorem \ref{thm2.50}). Let $C = \Im \psi \subset B$ and by (i) $B = \ker \pi \oplus \Im \psi = A \oplus C$. Clearly $\psi:B/A \to C$ is an epimorphism, and if $\psi(x) = \psi(y)$, then $x = \pi \psi(x) = \pi \psi(y) = y$ so $\psi$ is an isomorphism and $C \cong B/A$. [(iii)] Clearly $\varphi$ is well-defined and since $\varphi(a+b) = (a+b)+B = (a+B)+(b+B) = \varphi(a)+\varphi(b)$ and $\varphi(ra) = ra+B = r(a+B) = r\varphi(a)$, $\varphi$ is a homomorphism. If $(a+b)+B \in (A+B)/B$, then $\varphi(a) = a+B = (a+b)+B$ and $\varphi$ is an epimorphism. Since $\ker \varphi = \{a \in A \mid a+B = 0\} = A \cap B$, by the first isomorphism theorem for modules, $A/(A \cap B) \cong (A+B)/B$, as desired. $\blacksquare$
\begin{theorem}\label{thm2.54}
If $R$ is a principal ideal domain, then every submodule $G$ of a free $R$-module $F$ is itself free.
\end{theorem}
\noindent \textbf{Proof} Let $X = \{x_i \mid i \in \textit{I}\}$ be a basis of $F$. Then $F = \bigoplus_{i \in \textit{I}} Rx_i$ with each $Rx_i$ isomorphic to $R$ (as a left $R$-module) by Theorem \ref{thm2.100}. Choose a well ordering $\leq$ of $\textit{I}$ (AC). Given $i \in \textit{I}$, denote the immediate successor of $i$ by $i+1$ (if it exists). Let $J = I \cup \{\alpha\}$ ($a \notin I$), with order $\leq$ on $\textit{I}$ and $i < \alpha$ for all $i \in \textit{I}$. Then $J$ is well ordered and every element of $\textit{I}$ has an immediate successor in $J$. For each $j \in J$, let $F_j$ be the submodule of $F$ generated by the set $\{x_i \mid i < j\}$.We claim the following: \begin{flushleft}
\hspace*{2em}(i) $j < k \Leftrightarrow F_j \subsetneq F_k;$ \\
\hspace*{2em}(ii) $\bigcup_{j \in J} F_j = F;$ \\
\hspace*{2em}(iii) For each $i \in \textit{I}$, $F_{i+1}/F_{i} \cong Rx_i \cong R.$ \\
For each $j \in J$, let $G_{j} = G \cap F_j$. \\
\hspace*{2em}(iv) $j < k \implies G_{j} \subset G_{k};$ \\
\hspace*{2em}(v) $\cup_{j \in J} G_j = G;$ \\
\hspace*{2em}(vi) For each $i \in \textit{I}$, $G_i = G_{i+1} \cap F_i.$
\end{flushleft}

\noindent [(i)] If $j < k$, clearly $F_j \subsetneq F_k$, and if $F_j \subsetneq F_k$, then every $i < j$ also has $i < k$ so it follows that $j < k$ ($j \neq k$ as $F_j$ proper subset). [(ii)] Since $F_{\alpha} = F$, clearly $\bigcup_{j \in J} F_j = F$. [(iii)] By Theorem \ref{thm2.100}, we have $F_{i+1} = \bigoplus_{k < i+1} R x_k$ and the epimorphism  $F_{i+1} = \bigoplus_{k < i+1} R x_k \xrightarrow{\pi_{i}} Rx_i$ has kernel $F_i$, whence by the first isomorphism theorem for modules, $F_{i+1}/F_i \cong Rx_i \cong R$. [(iv)] If $j < k$, then $F_j \subset F_k$, whence $G_j \subset G_k$. [(v)] $\bigcup_{j \in J} G_j = G \cap \bigcup_{j \in J} F_j = G \cap F = G$. [(vi)] Fix $i \in \textit{I}$. Then $G_{i+1} \cap F_i = (G \cap F_{i+1}) \cap F_i = G \cap (F_{i+1} \cap F_{i}) = G \cap F_i = G_i$.

By property (vi) and Lemma \ref{lem2.53}(iii) $G_{i+1}/G_i = G_{i+1}/(G_{i+1} \cap F_i) \cong (G_{i+1}+F_i)/F_i$ for each $i \in \textit{I}$. Notice $(G_{i+1}+F_i)/F_i$ is a submodule of $F_{i+1}/F_i$ and hence $G_{i+1}/G_i$ is isomorphic to a submodule of $R$ by (iii). But every submodule of $R$ is necessarily an ideal of $R$, and since $R$ is a principal ideal domain, of the form $(a) = Ra$ for some $a \in R$. If $a \neq 0$, since $R$ is an integral domain, the map $R \to Ra$ by $r \mapsto ra$ is an isomorphism and hence $G_{i+1}/G_{i}$ is free for each $i \in \textit{I}$ and has a basis of cardinality either $0$ or $1$ (Theorem \ref{thm2.100}). By Theorems \ref{thm2.47} and \ref{thm2.52} the short exact sequence $0 \xrightarrow{} G_i \xrightarrow{\subset} G_{i+1} \xrightarrow{\pi} G_{i+1}/G_i \xrightarrow{} 0$ is split exact for each $i \in \textit{I}$. If $G_{i+1} = G_{i}$ for some $i \in \textit{I}$, then $G_{i+1} = G_{i+1} \oplus 0 = G_{i+1} \oplus Rb_i$ with $b_i = 0$. If $G_{i+1} \neq G_{i}$, by Lemma \ref{lem2.53}(ii), each $G_{i+1}$ is an internal direct sum $G_{i+1} = G_{i} \oplus C$, where $C \subset G_{i+1}$ is a submodule and $C \cong G_{i+1}/G_{i} \cong Rc_i$ (where $c_i \in R$). Consider $\alpha:Rc_i \cong C$, and $C = R\alpha(c_i)$. Let $b_i = \alpha(c_i)$ and $G_{i+1} = G_{i} \oplus Rb_i$ with $b_i \in G_{i+1} \setminus G_i$ (since $G_i \cap Rb_i = \{0\}$). Define $B = \{b_i \mid b_i \neq 0\}$. We claim $B$ is a basis of $G$. 

Suppose $u = \sum_{j} r_j b_j = 0$ ($j \in \textit{I}; r_j \in R$; finite sum). By way of contradiction, let $k$ be the largest index such that $r_k \neq 0$. Then $u = \sum_{j < k} r_jb_j+r_kb_k \in G_{k} \oplus Rb_k = G_{k+1}$. Since $u = 0$, $r_k b_k = -\sum_{j<k} r_j b_k \in G_{k} \cap Rb_k = \{0\}$, so $r_k b_k = 0$. Then write $b_k = \sum_{i = 1}^{n} s_i x_i$ ($s_i \in R, x_i \in X$) and $r_kb_k = \sum_{i = 1}^{n} (r_k s_i) x_i = 0$. Since $X$ is a basis, $b_k \neq 0$, and $R$ is an integral domain, this forces $r_k = 0$, a contradiction. Thus $r_j = 0$ for all $j$ and $B$ is linearly independent. 

Finally, we prove $B$ spans $G$. By (v) it suffices to prove that for each $k \in J$ the subset $B_k = \{b_j \in B \mid j < k\}$ of $B$ spans $G_k$. We proceed by transfinite induction (Theorem \ref{thm1.26}). If $k = \min J$, then $B_k = \emptyset$, which spans $G_k = \{0\}$. Therefore, by way of induction, suppose $k \neq \min J$, $B_j$ spans $G_j$ for all $j < k$, and let $u \in G_k$. If $k = j+1$ for some $j \in \textit{I}$, then $G_k = G_{j+1} = G_j \oplus Rb_j$ and $u = v+rb_j$ with $v \in G_j$. By induction hypothesis, $v = \sum r_i b_i$ with $r_i \in R$ and $b_i \in B_j \subset B_k$. Thus $u = \sum r_i b_i +rb_j$ and $B_k$ spans $G_k$. If $k \neq j+1$ for all $j \in \textit{I}$, since $u \in G_k = G \cap F_k$, write $u$ as a finite sum $u = \sum r_jx_j$ with $j<k$ ($r_j \in R; x_j \in X$). If $u = 0$, clearly $u$ is in the span of $B_k$. Otherwise, let $t$ be the largest index such that $r_t \neq 0$. Then $u \in F_{t+1}$ with $t+1 < k$ by hypothesis. Thus, $u \in G \cap F_{t+1} = G_{t+1}$ with $t+1 < k$. By induction hypothesis $u$ is in the span of $B_{t+1} \subset B_k$. Hence $B_k$ spans $G_k$, as desired. $\blacksquare$
\begin{corollary}\label{cor2.55}
Let $R$ be a principal ideal domain. Then an $R$-module $P$ is projective if and only if $P$ is free.
\end{corollary}
\noindent \textbf{Proof} By Theorem \ref{thm2.47}, it suffices to show only that if $P$ is projective it is free. By Theorem \ref{thm2.52}, $F \cong K \oplus P$ for some $R$-module $K$ and free $R$-module $F$. Then $P \xrightarrow{\iota_2} K \oplus P \cong F$ is a monomorphism and $P$ is (isomorphic to) a submodule of $F$. Hence, by Theorem \ref{thm2.54}, $P$ is free, as desired. $\blacksquare$

Corollary \ref{cor2.55} shows that, under certain circumstances, projective $R$-modules are equivalent to free $R$-modules (e.g., $R = \mbb{Z}$ by Example \ref{ex2.200}). As we saw above (Theorem \ref{thm2.100}, Theorem \ref{thm2.9}), there are many instances in which free modules are related to infinity, either through infinite bases, or being the internal direct sum of infinite cyclic modules/groups, and thus Corollary \ref{cor2.55} shows a relationship between infinity and projective $R$-modules.

The proof of Theorem \ref{thm2.54}, which is crucial for Corollary \ref{cor2.55}, relies on the principal of transfinite induction. Regular induction is already related to infinity due to its reliance on the axiom of infinity and the well-ordering principal, but transfinite induction goes even further. It extends the notion of induction to all well-ordered sets, even those with uncountably many elements. In the case above, in a similar vein to the many results that rely on Zorn's Lemma above, to assume any set has a well-ordering, and therefore compatible with transfinite induction, requires the axiom of choice. Therefore, the use of transfinite induction in this context reinforces the deep interconnection between infinity, the axiom of choice, and the structure theory of modules.
\section*{Conclusion}
This work undertook a structural analysis of infinity across set theory and modern algebra, revealing both fundamental connections and differences between finite and infinite structures across the two disciplines, as well as demonstrating how infinity interacts between them. In set theory, the analysis of infinite cardinals and ordinals highlights the former's  counterintuitive structural and arithmetic properties, in contrast to the more intuitive structure of the latter's. The axiom of choice played a crucial role in many results throughout the paper, and underscores the importance of such an axiom in mathematics with infinity. Transitioning to algebra, the examination of infinite groups, rings, and $R$-modules demonstrated how the introduction of infinity can extend or alter structural properties observed in a finite setting, as illustrated by, for instance, the extensions from weak to general direct products and polynomial rings to formal power series. Moreover, the analysis above highlighted the distinct yet analogous theory the introduction of infinity can create, for example in the comparative theory of finite and infinite dimensional vector spaces, along with demonstrating how more general algebraic structures, such as finitely generated abelian groups or projective $R$-modules, embody infinity in their classifications and structural properties.
\section*{Acknowledgments}
The author would like to thank Professor Milos Podmanik and Alexander West for reviewing early manuscripts of this paper—the former for valuable literary suggestions, and the latter for mathematical corrections and insights. A great intellectual debt is owed to the many mathematicians who contributed to the development of set theory and algebra over the past two centuries, to whom a short list could not do justice. Without those who proved, slowly sharpened, and deepened many of the results above, a concise treatment such as this would not be possible. 
\newpage
\bibliographystyle{plain}
\bibliography{math.bib}

\begin{thebibliography}{10}

\bibitem{abbott}
S.~Abbott.
\newblock {\em Understanding Analysis}.
\newblock Undergraduate Texts in Mathematics. Springer New York, 2015.

\bibitem{macdonald}
M.~F. Atiyah and I.~G. Macdonald.
\newblock {\em Introduction to Commutative Algebra}.
\newblock Addison-Wesley Series in Mathematics. Addison-Wesley Publishing
  Company, 1969.

\bibitem{axlerlinear}
S.~Axler.
\newblock {\em Linear Algebra Done Right}.
\newblock Undergraduate Texts in Mathematics. Springer New York, 3rd edition,
  2015.

\bibitem{axlermeasuresupplement}
S.~Axler.
\newblock {\em Supplement: Measure, Integration \& Real Analysis}.
\newblock Graduate Texts in Mathematics. Springer International Publishing,
  2019.

\bibitem{cohen}
P.~J. Cohen.
\newblock {\em Set Theory and the Continuum Hypothesis}.
\newblock W. A. Benjamin, New York, 1966.

\bibitem{conrad}
K.~Conrad.
\newblock Zorn’s lemma and some applications.
\newblock {\em Expository papers}, 2016.

\bibitem{dummit}
D.~S. Dummit and R.~M. Foote.
\newblock {\em Abstract Algebra}.
\newblock John Wiley and Sons, 2004.

\bibitem{enderton}
H.~B. Enderton.
\newblock {\em Elements of set theory}.
\newblock Academic Press, 1977.

\bibitem{fraleigh}
J.~B. Fraleigh.
\newblock {\em A first course in abstract algebra}.
\newblock Addison-Wesley, 7th edition, 2003.

\bibitem{grabczewski}
K.~Grabczewski.
\newblock Equivalents of the axiom of choice.
\newblock volume~34 of {\em Studies in Logic and the Foundations of
  Mathematics}, pages 1--134. Elsevier, 2001.

\bibitem{godel2}
K.~Gödel.
\newblock What is cantor's continuum problem?
\newblock {\em The American Mathematical Monthly}, 54(9):515--525, 1947.

\bibitem{godel1}
K.~Gödel.
\newblock {\em Consistency of the Continuum Hypothesis. (AM-3)}.
\newblock Princeton University Press, 1968.

\bibitem{halmos}
P.~R. Halmos.
\newblock {\em Naive Set Theory}.
\newblock Undergraduate Texts in Mathematics. Springer New York, 1998.

\bibitem{hewitt}
E.~Hewitt and K.~Stromberg.
\newblock {\em Real and Abstract Analysis: A modern treatment of the theory of
  functions of a real variable}.
\newblock Springer Berlin Heidelberg, 2013.

\bibitem{hilbert}
D.~Hilbert.
\newblock Ueber die theorie der algebraischen formen.
\newblock {\em Mathematische Annalen}, 36:473--534, 1890.

\bibitem{hilbertlect}
D.~Hilbert.
\newblock {\em David Hilbert's Lectures on the Foundations of Arithmetic and
  Logic 1917-1933}.
\newblock Springer Berlin, Heidelberg, 2013.

\bibitem{hodges}
W.~Hodges.
\newblock Krull implies zorn.
\newblock {\em Journal of the London Mathematical Society}, s2-19:285--287,
  1979.

\bibitem{howard}
P.~Howard and J.~E. Rubin.
\newblock {\em Consequences of the Axiom of Choice}, volume~59 of {\em
  Mathematical Surveys and Monographs}.
\newblock American Mathematical Society, 1998.

\bibitem{jechintro}
K.~Hrbacek and T.~Jech.
\newblock {\em Introduction to set theory}.
\newblock Monographs and textbooks in pure and applied mathematics. Marcel
  Dekker, 3rd edition, 1999.

\bibitem{hungerford}
T.~W. Hungerford.
\newblock {\em Algebra}.
\newblock Graduate Texts in Mathematics. Springer New York, 2nd edition, 2003.

\bibitem{honig}
C.~S. Hönig.
\newblock Proof of the well-ordering of cardinal numbers.
\newblock {\em Proceedings of the American Mathematical Society}, 5:312, 1954.

\bibitem{jacob1}
N.~Jacobson.
\newblock {\em Basic algebra 1}.
\newblock W H Freeman and Company, 2nd edition, 1985.

\bibitem{jech}
T.~Jech.
\newblock {\em Set Theory: The Third Millennium Edition, revised and expanded}.
\newblock Springer Monographs in Mathematics. Springer Berlin Heidelberg, 2007.

\bibitem{kaplansky}
I.~Kaplansky.
\newblock {\em Set theory and metric spaces}.
\newblock Allyn and Bacon Series in Advanced Mathematics. Allyn and Bacon, 2nd
  edition, 1972.

\bibitem{krull}
W.~Krull.
\newblock Idealtheorie in ringen ohne endlichkeitsbedingung.
\newblock {\em Mathematische Annalen}, 101:729--744, 1929.

\bibitem{lang}
S.~Lang.
\newblock {\em Algebra}.
\newblock Graduate Texts in Mathematics. Springer New York, revised 3rd
  edition, 2002.

\bibitem{larson}
Larson and Edwards.
\newblock {\em Calculus}.
\newblock Cengage Learning, 9th edition, 2010.

\bibitem{matsumura}
H.~Matsumura.
\newblock {\em Commutative Ring Theory}.
\newblock Cambridge Studies in Advanced Mathematics. Cambridge University
  Press, 1989.

\bibitem{moore}
G.~H. Moore.
\newblock {\em Zermelo's Axiom of Choice: Its Origins, Development, and
  Influence}, volume~8 of {\em Studies in the History of Mathematics and
  Physical Sciences}.
\newblock Springer-Verlag New York, 1982.

\bibitem{munkres}
J.~R. Munkres.
\newblock {\em Topology: Pearson New International Edition}.
\newblock Pearson Education, 2nd edition, 2013.

\bibitem{navarro}
G.~Navarro.
\newblock On the fundamental theorem of finite abelian groups.
\newblock {\em The American Mathematical Monthly}, 110(2):153--154, 2003.

\bibitem{niven}
I.~Niven, H.~S. Zuckerman, and H.~L. Montgomery.
\newblock {\em An introduction to the theory of numbers}.
\newblock John Wiley and Sons, 5th edition, 1991.

\bibitem{potter}
M.~Potter.
\newblock {\em Set Theory and its Philosophy: A Critical Introduction}.
\newblock Oxford University Press, 2004.

\bibitem{roitman}
J.~Roitman.
\newblock {\em Introduction to Modern Set Theory}.
\newblock John Wiley \& Sons, 1990.

\bibitem{roman}
S.~Roman.
\newblock {\em Advanced Linear Algebra}.
\newblock Graduate Texts in Mathematics. Springer New York, 2nd edition, 2005.

\bibitem{rotman}
J.~J. Rotman.
\newblock {\em Advanced Modern Algebra}.
\newblock Prentice Hall, 2nd edition, 2003.

\bibitem{babyrudin}
W.~Rudin.
\newblock {\em Principles of Mathematical Analysis}.
\newblock International series in pure and applied mathematics. McGraw-Hill,
  3rd edition, 1976.

\bibitem{shoenfieldaxioms}
J.~R. Shoenfield.
\newblock Axioms of set theory.
\newblock volume~90 of {\em Studies in Logic and the Foundations of
  Mathematics}, pages 321--344. Elsevier, 1977.

\bibitem{spivak}
M.~Spivak.
\newblock {\em Calculus}.
\newblock Publish or Perish, 3rd edition, 1994.

\bibitem{stoll}
R.~R. Stoll.
\newblock {\em Set Theory and Logic}.
\newblock Dover Books on Mathematics. Dover Publications, 2012.

\bibitem{suppes}
P.~Suppes.
\newblock {\em Axiomatic Set Theory}.
\newblock D. Van Nostrand Company, 1960.

\bibitem{whitehead}
A.~Whitehead.
\newblock On cardinal numbers.
\newblock {\em American Journal of Mathematics}, 24(4):367--394, 1902.

\bibitem{ziemer}
W.~P. Ziemer.
\newblock {\em Modern real analysis}.
\newblock Graduate texts in mathematics. Springer, 2nd edition, 2017.

\end{thebibliography}
\end{document}